%% file: Colorful_Minors_without_EP_arxiv.tex
\documentclass[11pt]{article}

\input{formatting}

 \newboolean{appear}

\setboolean{appear}{true}

 \ifthenelse{\boolean{appear}}{%

 \newcommand{\SW}[1]{%
 \todo[linecolor=red,bordercolor=white,color=Magenta,
 backgroundcolor=CornflowerBlue!40,textcolor=BrightPink,
 size=\footnotesize]{\textsf{SW: #1}}}

 \newcommand{\SWin}[1]{%
 \todo[inline,bordercolor=white,color=Magenta,
 backgroundcolor=CornflowerBlue!40,textcolor=BrightPink,
 size=\footnotesize]{\textsf{SW: #1}}}

 \newcommand{\EP}[1]{%
 \todo[linecolor=red,bordercolor=white,color=DeepCarrotOrange,
 backgroundcolor=BananaYellow!40,textcolor=DeepCarrotOrange,
 size=\footnotesize]{\textsf{EP: #1}}}

 \newcommand{\EPin}[1]{%
 \todo[inline,bordercolor=white,color=DeepCarrotOrange,
 backgroundcolor=BananaYellow!40,textcolor=DeepCarrotOrange,
 size=\footnotesize]{\textsf{EP: #1}}}

 \newcommand{\sed}[1]{%
 \todo[linecolor=red,bordercolor=white,color=BrilliantRose,
 backgroundcolor=BrilliantRose!50,size=\footnotesize]{%
 \textsf{\footnotesize SED: #1}}}

 \newcommand{\VP}[1]{%
 \todo[linecolor=red,bordercolor=white,color=purple!40,
 backgroundcolor=cyan!20,textcolor=magenta,
 size=\footnotesize]{\sf Vagelis: #1}}

 \newcommand{\VPin}[1]{%
 \todo[inline,bordercolor=white,color=purple!40,
 backgroundcolor=cyan!20,textcolor=magenta,
 size=\footnotesize]{\sf Vagelis: #1}}

}{%

 \newcommand{\SW}[1]{}
 \newcommand{\SWin}[1]{}
 \newcommand{\EP}[1]{}
 \newcommand{\EPin}[1]{}
 \newcommand{\sed}[1]{}
 \newcommand{\VP}[1]{}
 \newcommand{\VPin}[1]{}

}

\ifthenelse{\boolean{appear}}{}{}

\title{Colorful Minors}

 \author{%
 Evangelos Protopapas\thanks{Faculty of Mathematics, Informatics and Mechanics, University of Warsaw, Poland.\\ Email: \texttt{eprotopapas@mimuw.edu.pl}}~$^{,}$\thanks{Supported by the ERC project BUKA (n°\! 101126229) and the French-German Collaboration ANR/DFG Project UTMA (ANR-20-CE92-0027).}\and 
 Dimitrios M. Thilikos\thanks{LIRMM, Univ Montpellier, CNRS, Montpellier, France.\\ Email: \texttt{sedthilk@thilikos.info}.}~$^{,}$\thanks{
Supported by the Franco-Norwegian project PHC AURORA 2024-25 (Projet n°\! 51260WL) and the French National Research Agency (ANR) under project GODASse ANR-24-CE48-4377 and under the France 2030 grant reference number ANR-24-RRII-0002 operated by the Inria Quadrant Program.}\and 
 Sebastian Wiederrecht\thanks{KAIST, School of Computing, Daejeon, South Korea. \\ Email:
 \texttt{wiederrecht@kaist.ac.kr}}~$^{,}$\thanks{Supported by the Institute for Basic Science (IBS-R029-C1)}}
\date{}

\begin{document}

\maketitle
\vspace{-15mm}

\begin{abstract}
\noindent We introduce the notion of \emph{colorful minors}, which generalizes the classical concept of rooted minors in graphs.
{A \emph{$q$-colorful graph} is defined as a pair $(G, \chi),$ where $G$ is a graph and $\chi$ assigns to each vertex a (possibly empty) subset of at most $q$ colors.}
The colorful minor relation enhances the classical minor relation by merging color sets at contracted edges and allowing the removal of colors from vertices.
This framework naturally models algorithmic problems involving graphs with (possibly overlapping) annotated vertex sets.
We develop a structural theory for colorful minors by establishing three core theorems characterizing $\mathcal{H}$-colorful minor-free graphs, where $\mathcal{H}$ consists either of a clique or a grid with all vertices assigned all colors, or of grids with colors segregated and ordered on the outer face. Our results reveal that when exclusion is imposed not only on graphs but also to the way colors are distributed in them, a more refined structural landscape appears.

On the algorithmic side, we deduce that colorful minor testing is fixed-parameter tractable.
Together with the fact that the colorful minor relation forms a well-quasi-order, this implies that every colorful minor-monotone parameter on colorful graphs admits a fixed-parameter algorithm.
Furthermore, we derive two algorithmic meta-theorems (AMTs) whose structural conditions are linked to extensions of treewidth and Hadwiger number on colorful graphs.
Our results suggest how known AMTs can be extended to incorporate not only the structure of the input graph but also the way the colored vertices are distributed in it.
\end{abstract}


\thispagestyle{empty}

\newpage
\thispagestyle{empty}


\tableofcontents
\thispagestyle{empty}
\newpage

\setcounter{page}{1}

\section{Introduction}

Over the last decades, a fruitful interplay between structural graph theory and theoretical computer science has been the search for structural properties of graphs that facilitate the design of efficient algorithms for otherwise \textsf{NP}-hard problems (see for example \cite{CyganFKLMPPS2015Parameterized}).
One of the most studied such parameters is a measure for the ``tree-likeness'' known in its most popular form as \textsl{treewidth} \cite{RobertsonS1986Graphb,ArnborgP1989Linear,BodlaenderK2007Combinatorial,BodlaenderCKN2015Deterministic}.
Treewidth plays an exceptional role not only in the world of (parameterized) algorithms, but also in the theory of graph minors by Robertson and Seymour \cite{RobertsonS1984Grapha,RobertsonS1986Grapha,RobertsonS1990Grapha,RobertsonS1991Graph}.
A central theorem of Robertson and Seymour is the \textsl{Grid Theorem} \cite{RobertsonS1986Grapha} that states that any minor-closed graph class $\Ccal$ has bounded treewidth if and only if $\Ccal$ does not contain all planar graphs.
This reveals a fundamental limitation in the use of treewidth (or any minor-monotone parameter really) as a parametrization for an \textsf{NP}-hard problem:
While there exists a wide range of problems that are tractable on classes of bounded treewidth, there exist many important problems that are \textsf{NP}-hard even on planar graphs \cite{GareyJ1977Rectilinear,Curticapean2016Counting}.
This means that, in the regime of minor-closed graph classes, such problems are tractable if and only if a planar graph is excluded, unless $\mathsf{P}=\mathsf{NP}$.

A particularly interesting class of problems tractable on graphs of bounded treewidth, but \textsf{NP}-hard on more general minor-closed graph classes, are certain problems whose instances come with specific vertex sets, typically called \textsl{roots}, \textsl{terminals}, or \textsl{annotated vertices}.
A prototypical example of such a problem is \textsc{Steiner Tree}.
The \textsc{Steiner Tree} problem receives a graph $G$ accompanied with a \textsl{special} set $X\subseteq V(G)$ of annotated vertices as input and asks for the minimum number of edges in a connected subgraph of $G$ containing all vertices of $X.$
The problem is known to be \textsf{NP}-complete on planar graphs, however it can be solved in polynomial time with respect to structural parameterizations that instead of (or in addition to) restricting the \emph{global structure} of the graph (like treewidth), focus on the structure \emph{relative to the annotation} \cite{DreyfusW1971Steiner,EricksonMV1987SendAndSplit,KisfaludiBakNvL2020NearlyETH}, thereby providing tractable algorithms \textsl{beyond} the regime of minor-closed classes of bounded treewidth.

This hints at the possibility of developing a general structural and algorithmic theory for annotated graphs, providing tools that supersede the limits imposed by classical graph minor theory.

Computational problems such as \textsc{Steiner Tree} can be naturally understood as problems on annotated graphs.
Recently, several positive algorithmic results for \textsc{Steiner Tree} on annotated graphs excluding certain types of so-called \textsl{rooted minors}\footnote{An annotated graph $(H,Y)$ is a \emph{rooted minor} of an annotated graph $(G,X)$ if there exists a collection $\{ J_v\}_{v\in V(H)}$ of pairwise vertex-disjoint connected subgraphs of $G$ such that for all $uv\in E(H)$ there is an edge between $J_u$ and $J_v,$ and for all $v \in Y$ we have $V(J_v ) \cap X \neq \emptyset.$} have emerged \cite{GroendlandNK2024Polynomial,JansenS2024SteinerTree} alongside purely structural results on rooted minors \cite{MarxSW2017Rooted,Hodor2024Quickly,Fiorini2025Face} and found various applications \cite{Hodor2024Quickly,AprileFJKSWY2025Integer}.
Moreover, \textsc{Steiner Tree} is just one example where a theory of rooted minors may open a new horizon of tractability.
There exists a variety of other important problems exhibiting similar behavior, among them \textsc{Multiway Cut} \cite{Marx12, DahlhausEtAl1992, PandeyL2025Planar}, \textsc{Disjoint paths} \cite{OkamuraSeymour1981, Schrijver2003}, \textsc{\#Defect Matching} \cite{Curticapean2016DefectMatching}, and others \cite{Frederickson1991Planar,Curticapean2016Counting,KrauthgamerLRr2019FlowCut}.
Moreover, annotation has recently been shown to also play a key role in structural extensions of \textsl{Courcelle’s Theorem} \cite{Courcelle1990Monadic, Courcelle1992Monadic, Courcelle1997Expression} by restricting quantification to sets of low \textsl{bidimensionality}\footnote{The \emph{bidimensionality} of a vertex set $X$ in a graph $G$, denoted by $\bidim(G,X)$, is the largest $k$ such that $(G,X)$ contains the $(k \times k)$-grid with all vertices annotated as a rooted minor \cite{ThilikosW2025graphminorsstructure}.} \cite{SauST25Parameterizing}. 

While the concept of rooted minors in annotated graphs appears to be relatively well explored, a general theory for deriving results such as those above still seems to be lacking.
Beyond the fact that many of the results mentioned are inspired by the theory of graph minors, no systematic extension of its powerful toolkit has yet emerged.
This is somewhat surprising given the significant structural depth and algorithmic potential of these concepts, as we elaborate below.

\subsection{Our results}\label{subsec_Results}

In this paper, we develop a comprehensive theory for a generalization of rooted minors in multi-annotated graphs, which we refer to as \textsl{colorful minors} of \textsl{colorful graphs}. Instead of a single annotated set of vertices $X,$ we consider $q$ sets $X_1, \dots, X_q$ of annotated vertices, which may be overlapping. We investigate both the algorithmic and structural properties of colorful graphs that exclude certain classes of \textsl{colorful minors} -- an appropriate adaptation of the concept of rooted minors to our setting.
Our results suggest that much of the combinatorial and algorithmic power of classical graph minor theory can be effectively extended to the colorful setting.

\smallskip
\textbf{From structure...} Our core contributions are a series of structural theorems for $(H,\psi)$-colorful minor-free $q$-colorful graphs, both for general $(H,\psi)$ and for the particular case where $H$ is planar -- under two distinct assumptions on the distribution of colors.
Our results use as a departure point the notion of the bidimensionality of a vertex set and related structural concepts \cite{PaulPTW2024Obstructions, PaulPTW2024obstructionsArXiV, PaulPTS2025LocalIndex}.

\smallskip
\textbf{...to algorithms.} Our next step is to apply our structural theorems to algorithmic settings, with the broader goal of developing an algorithmic graph structure theory applicable to problems involving vertex annotations.
To this end, we first observe that colorful graphs are well quasi-ordered under the colorful minor relation.
We also provide an algorithm that decides whether a $q$-colorful graph $(H, \psi)$ is a colorful minor of a $q$-colorful graph $(G, \chi)$ in time\footnote{By $g(n)\in\Ocal_{k}(f(n))$ we mean that there is a function $h\colon\Nbbb\to\Nbbb$ such that $g(n)\in\Ocal(h(k)\cdot f(n)).$ Also, given a graph $G,$ we use $|G|$ and $|\!|G|\!|$ for the number of its vertices and edges respectively.} $\Ocal_{|H|+q}(|\!|G|\!|^{1+o(1)}).$
These two facts imply that membership for any colorful minor-closed class of colorful graphs is decidable in polynomial time. 

Next, we show how classical algorithmic meta-theorems (AMTs) -- such as Courcelle's Theorem and the recent algorithmic meta-theorem of Sau, Stamoulis, and Thilikos \cite{SauST25Parameterizing} -- can be extended to colorful graphs, where the structural conditions are expressed via the exclusion of colorful minors.
This framework enables the resolution of problems on graphs containing arbitrarily large grid or clique minors while, up to now, 
the exclusion of a minor appeared as a necessary (sparsity) condition.
Our results {support} a key insight {emerging from the body of work mentioned above}:
\vspace{-8pt}
\begin{center}
\textsl{Problem complexity does not only depend on the global structure of the input graph,}\\
\textsl{but crucially on how t
he annotated sets are distributed within it.}
\end{center}
\vspace{-8pt}
The algorithmic results presented here are not intended to be exhaustive; rather, they serve to illustrate this principle and to lay the groundwork for future exploration of algorithmic problems involving vertex annotations.
In particular, we provide algorithmic applications for two of our three main structural results -- \zcref{thm_intro_ExcludeRainbowClique,thm_restrictiveTreewidthIntro} -- while leaving the investigation of algorithmic consequences of \zcref{thm_ExcludingRainbowGridIntro} to future work.

\paragraph{Colorful graphs.}

Let $q$ be a non-negative integer.
A \emph{$q$-colorful graph} is a pair $(G,\chi)$ where $G$ is a graph and\footnote{For $n,m\in\mathbb{Z},$ we denote $\{ x \mid 1\leq x \leq n\text{, }x\in\mathbb{Z}\}$ by $[n]$ and $\{ x \mid n \leq x \leq m\text{, }x\in\mathbb{Z} \}$ by $[n,m].$} $\chi\colon V(G)\rightarrow 2^{[q]}.$
For a vertex $v\in V(G)$ we call $\chi(v)$ the \emph{palette} of $v.$
For illustrations, we use \textcolor{HotMagenta}{magenta} for $1,$ \textcolor{CornflowerBlue}{blue} for $2$, \textcolor{ChromeYellow}{yellow} for $3,$ and \textcolor{AppleGreen}{green} for $4.$

Certainly, every vertex of $G$ receives a (possibly empty) subset of colors.
Given a vertex set $X\subseteq V(G),$ 
we define $\chi(X)\coloneqq \bigcup_{v \in X} \chi(v)$.
We also use notation $\chi(G)\coloneqq \chi(V(G))$. 
Similarly, for $I\subseteq [q]$ we denote by $\chi^{-1}(I)$ the set $\{ v\in V(G) \colon I\cap\chi(v)\neq\emptyset \}.$
When $I=\{ i\}$ is a singleton, we write $\chi^{-1}(i)$ instead of $\chi^{-1}(\{ i\}).$

A $q$-colorful graph $(G,\chi)$ where $\chi(G)\subsetneq[q]$ is called \emph{restricted} and a $q$-colorful graph $(G,\chi)$ with $\chi(G)= \emptyset$ is called \emph{empty}.
Note that any restricted $q$-colorful graph is also a $(q-1)$-colorful graph, so the notion of being restricted is highly dependent on the context given by the total set of colors.
Finally, a $q$-colorful graph $(G,\chi)$ is said to be \emph{rainbow}\footnote{Notice that -- in slight violation of the English language -- we use the word \say{rainbow} as an adjective here.} if $\chi(v)=[q]$ for all $v\in V(G).$ 

A colorful graph $(H,\psi)$ is a \emph{(colorful) subgraph} of a colorful graph $(G,\chi),$ denoted by $(H,\psi)\subseteq (G,\chi),$ if $H$ is a subgraph of $G$ and $\psi(v)\subseteq \chi(v)$ for all $v\in V(H).$
If $(G,\chi)$ is a colorful graph and $H\subseteq G$ is a subgraph of $G$ we write, in slight abuse of notation, $(H,\chi)$ for the $q$-colorful subgraph of $(G,\chi)$ where all vertices of $H$ receive precisely the colors they received in $(G,\chi).$

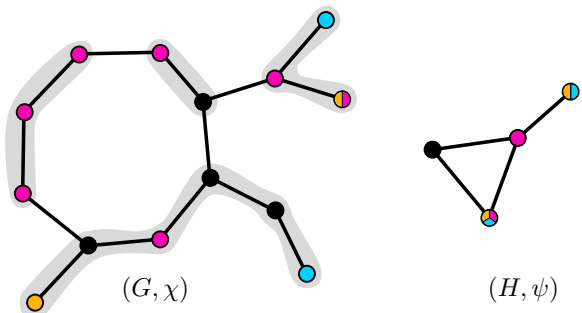
\begin{figure}[ht]
 \vspace{-8pt}
 \centering
 \scalebox{.89}{
 \begin{tikzpicture}

 \pgfdeclarelayer{background}
		\pgfdeclarelayer{foreground}
			
		\pgfsetlayers{background,main,foreground}
			
 \begin{pgfonlayer}{main}
 \node (M) [v:ghost] {};
 
 \node (MBottom) [v:ghost,position=270:20mm from M] {};

 \node (L) [v:ghost,position=180:20mm from MBottom] {$(G,\chi)$};
 \node (R) [v:ghost,position=0:35mm from MBottom] {$(H,\psi)$};

 \end{pgfonlayer}{main}

 \begin{pgfonlayer}{background}
 \pgftext{\includegraphics[width=10cm]{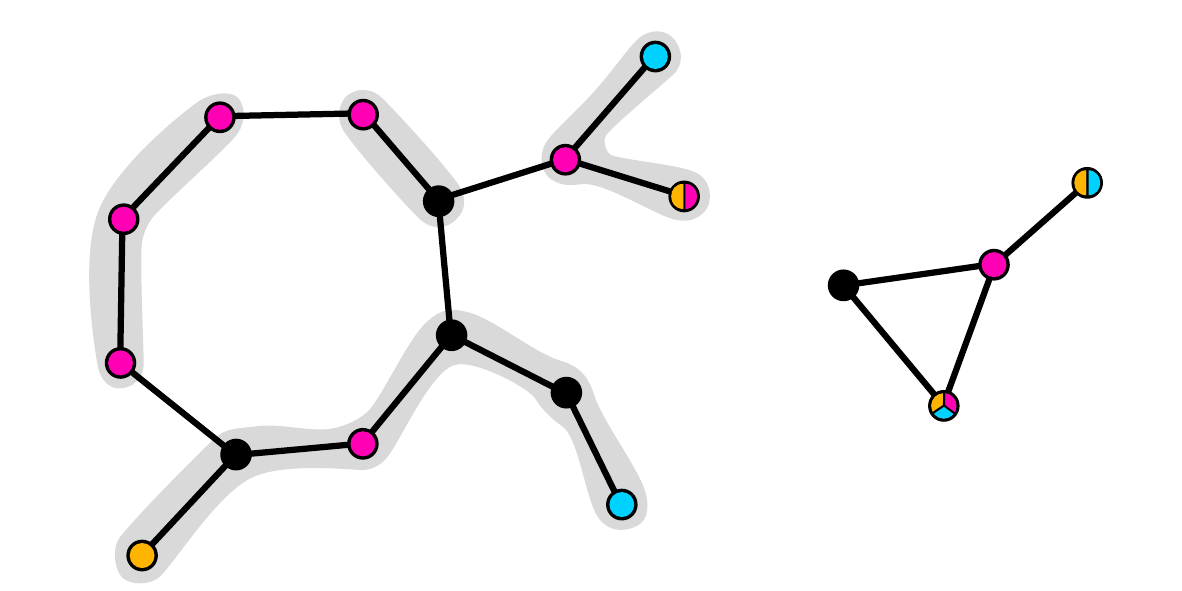}} at (M.center);
 \end{pgfonlayer}{background}
 
 \begin{pgfonlayer}{foreground}
 \end{pgfonlayer}{foreground}

 \end{tikzpicture}}
 \caption{Two colorful graphs $(G,\chi)$ and $(H,\psi)$ such that $(H,\psi)$ is a colorful minor of $(G,\chi).$
 The gray subgraphs of $(G,\chi)$ indicate the connected vertex sets that have to be contracted in order to form $(G,\psi).$ Notice that it is necessary to remove some colors from some of the vertices of $(G,\chi).$}
 \label{fig_ColorfulMinorIntro}
\end{figure}

A colorful graph $(H,\psi)$ is a \emph{colorful minor} of a colorful graph $(G,\chi)$ if $(H,\psi)$ can be obtained from $(G,\psi)$ by means of the following operations:
\begin{itemize}
 \item deleting an edge $e\in E(G),$
 \item deleting a vertex $v\in V(G)$ and restricting the domain of $\chi$ to $V(G)\setminus\{ v\},$
 \item for some $v\in V(G)$ and $i\in \chi(v)$ overwrite $\chi(v)\coloneqq \chi(v)\setminus\{ i\}$ (we refer to this operation as \emph{removing a color} (\emph{from $v$})), and
 \item contracting an edge $uv\in E(G),$ that is introducing a new vertex $x_{uv}$ with neighborhood $N_G(u)\cup N_G(v),$ setting $\chi(x_{uv})\coloneqq \chi(u) \cup \chi(v),$ and then deleting both $u$ and $v$ from the resulting colorful graph.
\end{itemize}
See \zcref{fig_ColorfulMinorIntro} for an example.
Notice that $1$-colorful graphs are exactly annotated graphs and the notions of colorful minors and rooted minors coincide on annotated graphs.

In \zcref{subsec_WQO} we derive the following, as a simple consequence of a deep theorem of Robertson and Seymour \cite{GraphMinorsXXIII}.

\begin{theorem}\label{thm_WQO}
For every non-negative integer $q,$ the class of all $q$-colorful graphs is well-quasi-ordered by the colorful minor relation.
\end{theorem}

\paragraph{Structural results.}

We prove a total of three structure theorems as follows.
The first two describe the structure of graphs excluding the $q$-colorful rainbow $t$-clique and the $q$-colorful rainbow $(k \times k)$-grid respectively as a colorful minor.
Moreover, we generalize a theorem of Marx, Seymour, and Wollan \cite{MarxSW2017Rooted}, that describes the structure of $1$-colorful graphs excluding a $(k \times k)$-grid where exactly the vertices of the first row receive color $1,$ to an arbitrary number of colors.

We begin with our simplest result -- excluding a rainbow clique -- which can be naturally phrased in terms of a modulator to restricted graphs.
For a vertex set $X\subseteq V(G)$ in a graph $G$ we define the \emph{torso} of $X$ in $G,$ denoted by $\torso(G,X),$ as the graph obtained from $G[X]$ by turning the neighborhood of $J$ in $X$ into a clique for every component $J$ of $G-X.$
In short, every $q$-colorful graph $(G,\chi)$ that excludes a rainbow clique as a colorful minor has a set $X$ whose torso excludes a minor and where every component of $(G-X,\chi)$ is free of some (not necessarily the same) color.
 
\begin{theorem}\label{thm_intro_ExcludeRainbowClique}
There exists a function $\rc\colon\mathbb{N}^2\to\mathbb{N}$ such that for all non-negative integers $q$ and $t,$ and all $q$-colorful graphs $(G,\chi)$ one of the following holds:
\begin{enumerate}
 \item $(G,\chi)$ contains a rainbow $t$-clique as a colorful minor, or
 \item there exists a set $X\subseteq V(G)$ such that the torso of $X$ in $G$ is $K_{\rc(q,t)}$-minor-free and for all components $J$ of $G-X,$ the $q$-colorful graph $(J,\chi)$ is restricted.
\end{enumerate}
Moreover, $\rc(q,t)\in\poly(qt)$ and there exists an algorithm that finds one of the two outcomes above in time $2^{\poly(qt)} \cdot |G|^{3}|\!|G|\!|\log|G|.$ 
\end{theorem} 
 
A \emph{tree-decomposition} for a graph $G$ is a pair $\Tcal=(T,\beta)$ where $T$ is a tree and $\beta\colon V(T)\to2^{V(G)},$ called the \emph{bags} of $\Tcal,$ assigns a set of vertices of $G$ to every node of $T$ such that $\bigcup_{t\in V(T)}\beta(t)=V(G),$ for every $e\in E(G)$ there is $t\in V(T)$ with $e\subseteq \beta(t),$ and for every $v\in V(G)$ the set $\{ t\in V(T) \mid v\in\beta(T)\}$ is connected.
The \emph{adhesion} of $\Tcal$ is $\max_{dt\in E(T)}|\beta(d)\cap \beta(t)|$ and the width of $\Tcal$ is defined as $\max_{t\in V(T)}|\beta(t)|-1.$
The \emph{treewidth} of a graph $G,$ denoted by $\tw(G),$ is the smallest integer $k$ such that there is a tree-decomposition of width at most $k$ for $G.$

\begin{figure}[ht]
 \vspace{-10pt}
 \centering
 \scalebox{1}{
 \begin{tikzpicture}

 \pgfdeclarelayer{background}
		\pgfdeclarelayer{foreground}
			
		\pgfsetlayers{background,main,foreground}
			
 \begin{pgfonlayer}{main}
 \node (M) [v:ghost] {};

 \end{pgfonlayer}{main}

 \begin{pgfonlayer}{background}
 \pgftext{\includegraphics[width=10.5cm]{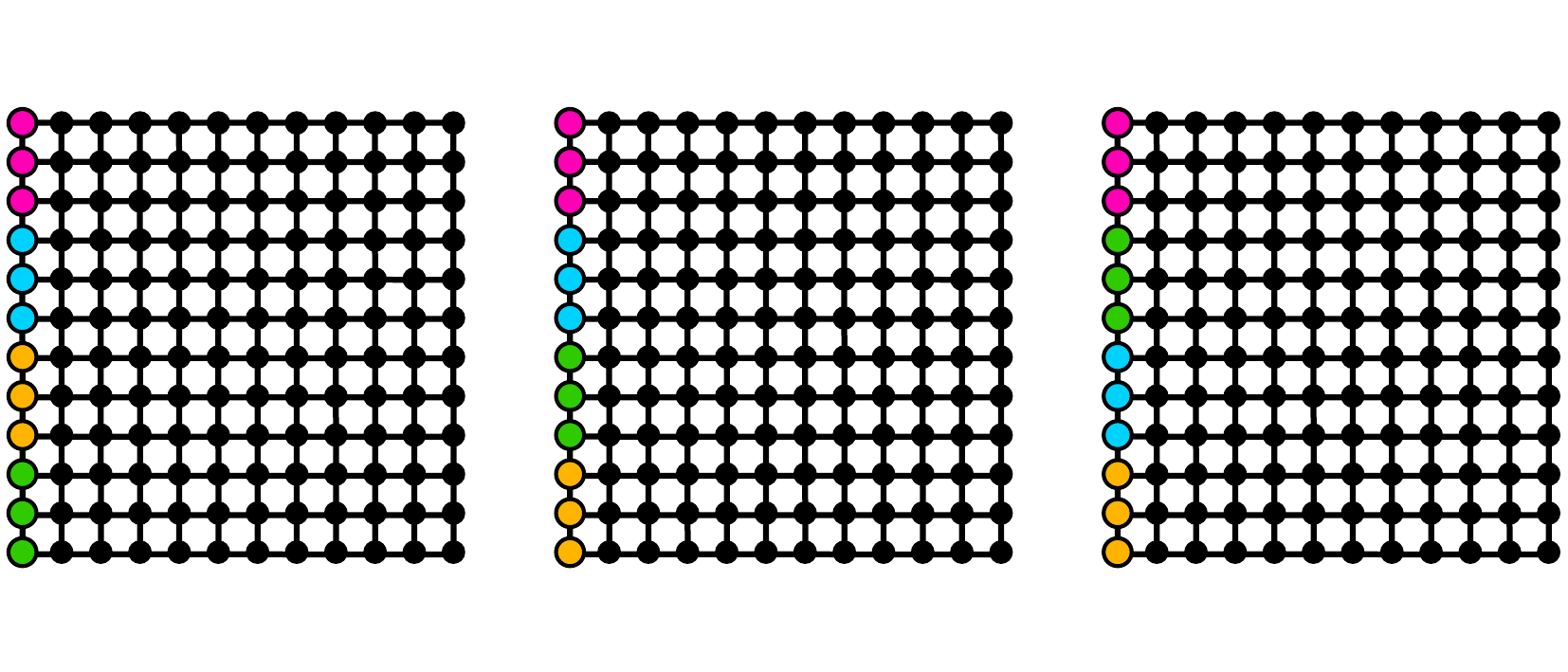}} at (M.center);
 \end{pgfonlayer}{background}
 
 \begin{pgfonlayer}{foreground}
 \end{pgfonlayer}{foreground}

 \end{tikzpicture}}
 \vspace{-14pt}
 \caption{Three non-isomorphic $(4,3)$-segregated grids.}
 \vspace{-4pt}
 \label{fig_SegregatedGrids}
\end{figure}

Our next theorem generalizes a theorem of Marx, Seymour, and Wollan \cite{MarxSW2017Rooted} (see also \cite{Hodor2024Quickly}).
Let $q\geq 1$ be an integer.
A \emph{$(q,k)$-segregated grid} is a $q$-colorful graph $(G,\chi)$ where $G$ is the $(qk\times qk)$-grid where the vertices of the first column can be numbered as $v_1,\dots,v_{qk}$ in order of their appearance, and we have $\chi(u)=\emptyset$ for all $u\in V(G)\setminus\{ v_1,\dots,v_{qk}\}$ and there exists a permutation $\pi$ of $[q]$ such that for every $i\in[q],$ $\chi(\{v_j \mid j\in[(i-1)q+1,iq] \})=\{ \pi(i)\}.$
In this situation, we say that $(G,\chi)$ \emph{realizes $\pi$}.
See \zcref{fig_SegregatedGrids} for an example.

\zcref{thm_intro_ExcludeRainbowClique} shows that every $(H,\psi)$-colorful minor-free colorful graph can be turned into a collection of restricted graphs by deleting an $H'$-minor-free subgraph.
Segregated grids capture exactly the situation where we require $H'$ to be a planar graph.
We may express this in form of a parameter as follows.
The \emph{torso treewidth} of a $q$-colorful graph $(G,\chi)$ is the smallest integer $k$ such that there exists a set $X\subseteq V(G)$ such that the treewidth of the torso of $X$ in $G$ is at most $k$ and for every component $J$ of $G-X$ the $q$-colorful graph $(J,\chi)$ is restricted. 

\begin{theorem}\label{thm_restrictiveTreewidthIntro}
There exists a function $\sg\colon\mathbb{N}^2\to\mathbb{N}$ such that for all non-negative integers $q$ and $k,$ and all $q$-colorful graphs $(G,\chi),$ $q\geq 1,$ one of the following holds:
\begin{enumerate}
\item $(G,\chi)$ contains as a colorful minor some $(q,k)$-segregated grid and the torso treewidth of $(G,\chi)$ is more than $k-1$, or 
\item the torso treewidth of $(G, \chi)$ is at most $\sg(q,k).$
\end{enumerate}
Moreover, $\sg(q,k)\in k^{2^{\Ocal(q)}}$ and there exists an algorithm that finds one of the two outcomes above in time $2^{k^{2^{\Ocal(q)}}} \cdot |G|^3|\!|G|\!|\log|G|.$
\end{theorem}

The second part of the first outcome is technically not necessary for the structural statement but serves as a pointer to this min-max relation.
Moreover, the min-max relation itself also supplies a guarantee for the approximation of torso treewidth via the algorithm from \zcref{thm_restrictiveTreewidthIntro}.

Notice that \zcref{thm_restrictiveTreewidthIntro} is optimal in the sense that the torso treewidth of segregated grids is unbounded (see \zcref{lemma_segregatedLowerbound}).
In other words, \zcref{thm_restrictiveTreewidthIntro} can be seen as a min-max duality theorem for segregated grids. We remark that for the case $q=1,$ the torso treewidth of $(G,\chi)$ has been independently defined by Jansen and Swennenhuis \cite{JansenS2024SteinerTree} and by Hodor, La, Micek, and Rambaud \cite{HodorLMR24quick} under different names.

It is easy to see that \zcref{thm_restrictiveTreewidthIntro} is far from a structure theorem for $q$-colorful graphs that exclude an arbitrary planar colorful minor.
This is in stark contrast to the $0$-colorful case provided by the Grid Theorem \cite{RobertsonS1986Grapha}.
It turns out that, in order to fully describe the structure of $q$-colorful graphs that exclude a fixed but arbitrary planar colorful minor, we need the full power of Robertson and Seymour's \emph{Graph Minor Structure Theorem} (GMST) \cite{RobertsonS2003GraphMinorsXVI,KawarabayashiTW2020Quickly,GorskySW2025Polynomial}.
The definitions necessary to state the GMST are rather involved and technical, for an in-depth explanation see \zcref{sec_RainbowGrid}.

We say that a graph $G$ has a \emph{$k$-near embedding} in a surface $\Sigma$ if there exists $A\subseteq V(G)$ with $|A|\leq k$ such that $G-A=G_0\cup G_1 \cup \dots \cup G_{\ell},$ $\ell\leq k,$ such that $G_0$ has an embedding into $\Sigma$ with $\ell$ pairwise vertex-disjoint faces $F_i$ such that for all $i\in[\ell],$ $V(G_0)\cap V(G_i)= V(F_i),$ the $G_i$'s are pairwise vertex-disjoint for positive $i,$ and for each $i\in[\ell],$ $G_{i}$ has a path-decomposition $(P_i,\beta_i)$\footnote{A \emph{path-decomposition} of a graph is a tree-decomposition $(T,\beta)$ where $T$ is a path.} of width at most $k$ such that $V(P_i)=V(F_i),$ the vertices of $P_{i}$ appear in agreement 
to their cyclic ordering on the boundary of $F_{i},$ and $v\in \beta_i(v),$ for all $v\in V(F_i).$
The set $A$ is called the \emph{apex set}, the graphs $G_i,$ $i\in[\ell],$ are called the \emph{vortices}, and the sets $V(G_i)\setminus V(F_i)$ are the \emph{interiors} of the vortices.

In the context of $q$-colorful graphs, we need some additional information regarding the colors.
Let $(T,\beta)$ be a tree decomposition of a $q$-colorful graph $(G,\chi).$
The \emph{colorful torso} $(G_t,\chi_t)$ of $(G,\chi)$ at a node $t\in V(T)$ is the colorful graph obtained from $(G[\beta(t)],\chi)$ by turning, for every $dt\in E(T),$ the set $\beta(d)\cap\beta(t)$ into a clique and adding $\chi(\bigcup_{h\in V(T_d)}(\beta(h)\setminus\beta(t)))$ to $\chi(v)$ for all $v\in\beta(d)\cap\beta(t)$ where $T_d$ is the unique component of $T-dt$ that contains $d.$ 

\begin{theorem}\label{thm_ExcludingRainbowGridIntro}

There exists a function $\rg\colon\mathbb{N}^2\to\mathbb{N}$ such that for all non-negative integers $q$ and $k,$ and all $q$-colorful graphs $(G,\chi)$ one of the following holds:
\begin{enumerate}
\item $(G,\chi)$ contains the $q$-colorful rainbow $(k \times k)$-grid as a colorful minor, or
\item $(G,\chi)$ has a tree-decomposition $(T,\beta)$ of adhesion at most $\rg(q,k)$ such that for all $t\in V(T),$ either $t$ is a leaf with unique neighbor $d$ and $(G[\beta(t)\setminus\beta(d)],\chi)$ is restricted, or the colorful torso $(G_t,\chi_t)$ of $G$ at $t$ has an $\rg(q,k)$-near embedding and there exists $\emptyset \neq I_t\subseteq [q]$ such that all vertices $v\in\beta(t)$ with $I_t\cap \chi_t(v)\neq\emptyset$ belong to the apex set or the interior of a vortex.
\end{enumerate}
Moreover, $\rg(q,k)\in 2^{k^{\Ocal(1)}2^{2^{\Ocal(q)}}}$ and there exists an algorithm that finds one of the two outcomes above in time $2^{2^{k^{\Ocal(1)}2^{2^{\Ocal(q)}}}} \cdot |G|^{3} |\!|G|\!| \log|G|.$
\end{theorem}

Notice that \zcref{thm_ExcludingRainbowGridIntro} can also be seen as a min-max duality as implied by results of Thilikos and Wiederrecht \cite{ThilikosW2025graphminorsstructure}.

To rephrase \zcref{thm_ExcludingRainbowGridIntro}, the absence of a rainbow grid as a colorful minor implies that, with respect to any large grid minor in $G$, and in the absence of a large clique minor—where a local version of \zcref{thm_excludeRainbowClique} would apply—there is always at least one color of small \textsl{bidimensionality}; moreover, all vertices of such colors can be confined to the apex set and the vortices of a near-embedding.

In particular, for $1$-colorful graphs, all colored vertices are confined entirely to the vortices and the apex set, while the surface part of the near-embedding contains no colored vertices. This appears to be a far-reaching generalization of the idea of covering all colored vertices in a planar graph by a bounded number of faces.

We believe that this provides a theoretical foundation for why many problems that are tractable on planar graphs when their annotated set is covered by few faces remain tractable on all colorful minor-closed classes of $1$-colorful graphs that exclude some fixed planar $1$-colorful graph. Moreover, it suggests that, in algorithm design, vortices in $k$-near embeddings may play a role analogous to that played by faces in planar graphs.
\paragraph{Towards a general framework for colorful parameters.}

With \zcref{thm_restrictiveTreewidthIntro} and \zcref{thm_ExcludingRainbowGridIntro} in hand, it becomes apparent that there are fundamentally distinct ways to extend a graph parameter to a colorful analogue, and that these choices can lead to remarkably different structural decompositions.
Which extension should be regarded as “canonical” appears to depend strongly on the intended application and on the structural/algorithmic consequences one wishes to capture.

For instance, from the perspective of structural decompositions, torso treewidth emerges as a particularly natural analogue of treewidth, mirroring the role played by treewidth in the uncolored setting.
In contrast, when viewed through the lens of obstructions, bidimensionality provides the appropriate analogue, as it yields the correct structure theorem for colorful graphs excluding an arbitrary colorful planar graph as a colorful minor, and potentially, the correct structural parameterization to obtain precise dichotomy complexity results for problems on colorful graphs that become \NP-hard in planar graphs, as previously discussed.

These observations point towards a broader, open-ended research direction, which is to better understand these phenomena by examining other classical graph parameters beyond treewidth, and by developing general principles for lifting such parameters to colorful analogues.
We return to this perspective in \zcref{sec_extend_to_colorful} of the conclusion, where we discuss this issue in more detail and propose a conjecture.

\paragraph{Algorithmic results.}
 
One of the ground-breaking results of Robertson and Seymour's Graph Minors Series is the \textsl{Graph Minor Algorithm} that allows to check if a graph $H$ is a minor of a graph $G$ in time $\Ocal_{|H|} (|G|^3)$ \cite{RobertsonS1995Graph} (see \cite{KawarabayashiW2010Shorter,KorhonenPS2024Minor} for improvements).
This algorithm also solves a rooted version of the \textsc{Minor Checking} problem which, as a special case, contains the so-called \textsc{$k$-Disjoint Paths} problem.

The \textsc{$q$-Colorful Minor Checking} problem has as input two $q$-colorful graphs $(G,\chi)$ and $(H,\psi),$ and the task is to determine whether $(G,\chi)$ contains $(H,\psi)$ as a colorful minor.

We present a simple reduction to the rooted minor checking problem which allows us to leverage the algorithm of Korhonen, Pilipczuk, and Stamoulis \cite{KorhonenPS2024Minor} in order to deduce a similar result for \textsc{$q$-Colorful Minor Checking}.

\begin{theorem}\label{thm_introColorfulMinors}
There exist an algorithm
that takes as input two $q$-colorful graphs $(H,\psi)$ and $(G,\chi)$ and solves \textsc{$q$-Colorful Minor Checking} in time $\Ocal_{q+|H|}(|\!|G|\!|^{1+o(1)}).$
\end{theorem}

\zcref{thm_WQO,thm_introColorfulMinors} imply that the recognition problem for any proper colorful minor-closed class can be solved in polynomial time. We present below a parametric variant of this, based on the general framework introduced by Fellows and Langston \cite{FellowsL1987Nonconstructive,FellowsL1988Nonconstructive} for the existence of fixed-parameter tractable algorithms. 

A \emph{($q$-colorful) graph parameter} is a {map} 
$\mathsf{p}$ assigning a non-negative integer to every ($q$-colorful) graph.
We say that a $q$-colorful graph parameter $\mathsf{p}$ is \emph{colorful minor-monotone} if for every $q$-colorful graph $(G,\chi)$ and every colorful minor $(H,\psi)$ of $(G,\chi)$ we have $\mathsf{p}(H,\psi) \leq \mathsf{p}(G,\chi).$ Again, combining \zcref{thm_WQO,thm_introColorfulMinors} we can deduce the following.

\begin{theorem}\label{thm_NonConstructiveColorfulParameters}
For every colorful minor-monotone $q$-colorful graph parameter $\mathsf{p},$ there exist an algorithm that, given a $q$-colorful graph $(G,\chi)$ and a 
non-negative integer $k,$ checks whether $\mathsf{p}(G,\chi)\leq k$ in 
$\Ocal_{q+k}(|\!|G|\!|^{1+o(1)})$ time.
\end{theorem}

We stress that the result of \zcref{thm_NonConstructiveColorfulParameters} is not constructive, i.e., it does not give any way to construct the claimed algorithm, as it is based on the finiteness of the obstruction set of the class 
$\{(G,\chi)\mid \mathsf{p}(G,\chi)\leq k\},$ which is only existentially provided by \zcref{thm_WQO}.

\paragraph{Meta-algorithmic Applications.}

Our two algorithmic results can be seen as extensions of Courcelle's Theorem 
and of a recent result due to Sau, Stamoulis, and Thilikos \cite{SauST25Parameterizing} to colorful graphs.
While the structural parameters in these theorems are treewidth, denoted by $\tw$ and Hadwiger number, denoted by $\hw$ (that is the maximum size of a clique minor) respectively,
we now introduce two parameters on colorful graphs that serve as our new combinatorial restrictions.

A key condition for applying our algorithmic results is the notion of \textsl{folio-representability} of a class of colorful graphs~$\Ccal$ (i.e., the \yes-instances of some algorithmic problem), formally defined in \zcref{sec_Applications}.

Informally, $\Ccal$ is \emph{folio-representable} if one can efficiently replace a subgraph~$P$ of a colorful graph~$(G, \chi)$ with colored vertices only within a small subset $B$ of the vertices of $G$ -- viewed as a \say{protuberance}\footnote{We use the term \textsl{protuberance} to distinguish it from the term \say{protrusion} introduced in \cite{BodlaenderFLPST2016Meta}, where similar replacement operations were applied.
Unlike protrusions, which are required to have bounded treewidth, protuberances are characterized by being \say{internally colorless}.} interacting with the rest of the graph via the \say{boundary} $B$ -- by a bounded-size subgraph~$P',$ such that the resulting graph~$(G', \chi)$ belongs to $\mathcal{C}$ if and only if $(G, \chi)$ does.

This replacement is guided by a notion of folio-equivalence: $P$ and $P'$ must have the same set of rooted topological minors (up to a given size), with respect to the specified boundary $B.$
Efficient computation of such replacements is enabled by results from~\cite{GroheKMW11Finding}, making folio-representability a natural combinatorial condition for the class $\Ccal.$

In what follows, whenever we refer to 1-colorful graphs we 
may alternatively use the notation $(G,X)$ instead of $(G,\chi)$ where $X\coloneqq \chi^{-1}(1)$.
Given a graph $G$ and a set $X\subseteq V(G)$, we define the \emph{monodimensionality} of $X$ in $G$ as the maximum order of a $(1,k)$-segregated grid colorful minor of $(G,X).$
Given a $q$-colorful graph $(G,\chi)$ and a color $i \in [q]$ we define the \emph{monodimensionality} of $i$ in $(G, \chi)$ as the monodimensionality of $\chi^{-1}(i)$ in $G.$
 
We now give the logical conditions under which our framework operates.
In our first result, we require the class $\Ccal$ to be definable in \textsl{Counting Monadic Second-Order Logic}\footnote{There exists a formula $\phi \in \mathsf{CMSO}$ such that $(G, \chi) \in \Ccal$ if and only if $(G, \chi)$ models $\phi.$} (CMSO). 
Courcelle's Theorem says that, for every graph class $\Gcal$ definable by some \CMSO-formula $\phi,$ membership in $\Gcal$ can be decided in time $\Ocal_{|\phi|+\tw(G)}(|G|).$ We prove the following analogue for colorful graphs.

\begin{theorem}
\label{thm_rwtMeta_Intro}
Let $q \in \mathbb{N},$ and let $\Ccal$ be a $q$-colorful graph problem that is both folio-representable and definable by some \CMSO-formula $\phi.$ Then, there exists an algorithm that, given a $q$-colorful graph $(G, \chi),$ decides whether $(G, \chi) \in \Ccal$ in time $\Ocal_{q + |\phi| + k}(|G|^{\Ocal(1)}),$ where $k$ is the maximum monodimensionality of a color in $(G, \chi)$.
\end{theorem}

From \zcref{thm_restrictiveTreewidthIntro}, the maximum monodimensionality of a color in a $q$-colorful graph $(G, \chi)$ is equivalent to the maximum torso treewidth of $(G, \chi^{-1}(i))$ over all colors $i \in [q].$
Notice that \zcref{thm_rwtMeta_Intro} can be restated by setting $k$ to be this equivalent quantity, which plays the same role for colorful graphs as treewidth does for graphs in Courcelle's Theorem.

For our second result, we consider a fragment of $\mathsf{CMSO}$ denoted as $\mathsf{CMSO/tw}+\mathsf{dp}.$
In this logic, quantification is restricted to sets of bounded bidimensionality, and a special predicate $\mathsf{dp}(t_1, \ldots, t_r, s_1, \ldots, s_r)$ is available, asserting the existence of $r$ vertex-disjoint paths connecting the pairs $(t_1, s_1), \ldots, (t_r, s_r).$
It was shown by Sau, Stamoulis, and Thilikos \cite{SauST25Parameterizing} that if a graph class $\Gcal$ is definable in $\mathsf{CMSO/tw}+\mathsf{dp},$ then membership in $\Gcal$ can be decided in time $\Ocal_{|\phi|+\mathsf{hw}(G)}(|G|^2)$.

The Hadwiger number of a color $i$ in a $q$-colorful graph $(G,\chi)$ is the maximum $k$ for which $(G,\chi^{-1}(i))$ contains the rainbow $k$-clique as a colorful minor.
We prove the following.

\begin{theorem}
\label{thm_rcMeta_Intro} 
Let $q \in \mathbb{N},$ and let $\Ccal$ be a $q$-colorful graph problem that is folio-representable and definable by some $\mathsf{CMSO/tw}+\mathsf{dp}$-formula $\phi.$ Then, there exists an algorithm that, given a $q$-colorful graph $(G, \chi),$ decides whether $(G, \chi) \in \Ccal$ in time $\Ocal_{q + |\phi| + k}(|G|^{\Ocal(1)})$, where $k$ is the maximum Hadwiger number of a color of $(G,\chi)$. 
\end{theorem}

We emphasize that the applicability of \zcref{thm_rwtMeta_Intro,thm_rcMeta_Intro}, based on the structural theorems \zcref{thm_intro_ExcludeRainbowClique,thm_restrictiveTreewidthIntro} respectively, exceeds the scope of Courcelle's Theorem and the results of \cite{SauST25Parameterizing}.
The generality of our extension is displayed by the fact that, for every $q≥1$, there are classes of $q$-colorful graphs $(G, \chi)$ for which the maximum monodimensionality/Hadwiger number of a color in $(G, \chi)$ is bounded while the treewidth/Hadwiger number of $G$ is unbounded.

This is a strong indication that, for problems on colorful graphs, the appropriate structural parameterizations toward deriving (meta-)algorithmic results should take into account the way the colors are distributed in the inputs.

We derive \zcref{thm_rwtMeta_Intro,thm_rcMeta_Intro} in \zcref{sec_Applications} as special cases of a more general result applicable to optimization problems (\zcref{thm_reduceToStars}).
In fact, in \zcref{sec_Applications}, we see these theorems as corollaries of broader statements (\zcref{the_MSOL} and \zcref{the_MSOL_tw}) formulated in terms of optimization problems where we additionally display the parameter dependencies hidden in the $\Ocal_{q+|\phi|+k}(|G|^{\Ocal(1)})$-notation.

\subsection{Brief insight into our approach}

The general approach to proving Robertson-Seymour-style structure theorems typically consists of establishing a local structure theorem with respect to a wall and then applying a well-known technique, originating from the work of Robertson and Seymour
\cite{RobertsonS1991Graph}, that turns the local statement into a global one based on tree-decompositions.
Our structure theorems \ref{thm_intro_ExcludeRainbowClique}, \ref{thm_restrictiveTreewidthIntro}, and \ref{thm_ExcludingRainbowGridIntro} follow this general paradigm.

In the uncolored setting, proving a local structure theorem with respect to a wall usually proceeds as follows: Starting from a large wall $W$ in the graph $G$, one either finds
\begin{enumerate}
\item[a.i)] a large clique-minor that is highly connected to $W$ or
\item[a.ii)] a structural decomposition of $G$ relative to $W,$ typically of a topological nature, such as a flat wall or, more generally, a near embedding of $G$ in some surface with additional properties tailored to the theorem at hand.
\end{enumerate}

For our results, handling outcome (a.ii) relies on sophisticated topological arguments stemming from the Graph Minor Structure Theorem \cite{RobertsonS2003GraphMinorsXVI} and its subsequent improvements \cite{KawarabayashiTW2020Quickly, GorskySW2025Polynomial} as well as more recent refinements \cite{ThilikosW2024Killing, PaulPTW2024Obstructions, PaulPTS2025LocalIndex}.
While technically demanding, this part largely builds on and adapts existing machinery to deal with annotation in near embeddings \cite{PaulPTS2025LocalIndex}, except for \zcref{thm_restrictiveTreewidthIntro} which uses different techniques which we comment on in the next paragraph.

\paragraph{Rainbow clique-minors and a multicolored Menger-type theorem.}

The genuinely new difficulties and contributions arise in outcome (a.i).
In the colored settings, the mere presence of a large clique-minor provides little information on the structure of the colorful graph with respect to the colored vertices.
For our purposes, it is therefore essential to understand how well the colored vertices are connected to this clique-minor.
Concretely, given a $q$-colorful graph $(G, \chi)$ containing a large clique-minor, we aim to establish the following dichotomy:
\begin{enumerate}
\item[b.i)] there exists a small set of vertices $S$ and a color $i \in [q]$ such that the component of $(G - S, \chi)$ containing most of the clique-minor is free of vertices of color $i,$ or
\item[b.ii)] $(G, \chi)$ contains a large rainbow clique-minor that is highly connected to the original clique-minor.
\end{enumerate}
This results, formalized in \zcref{thm_rainbowclique}, is proved in \zcref{subsec_rainbowclique} and constitutes a central new ingredient underpinning all our structural theorems.

In the context of the Disjoint Paths problem, a similar result was established by Robertson and Seymour \cite{RobertsonS1995Graph} for $1$-colorful graphs with a bounded number of colored vertices.
Our setting is substantially more general: we allow $q$-colorful graphs with an arbitrary number of colored vertices.
Our proof involves a \textsl{one-to-many} strengthening of Menger's Theorem, stated as \zcref{lemma_multicolorlinakge}, which may be of independent interest.
For all positive integers $k, \ell$, given vertex sets $X_{1}, \ldots, X_{\ell}$ and $Y$ in a graph $G,$ either
\begin{enumerate}
\item[c.i)] there exist $k \cdot \ell$ disjoint paths in $G,$ with exactly $k$ paths linking $X_{i}$ to $Y,$ or
\item[c.ii)] there exists an index $i \in [l]$ and a vertex set $S$ of size less than $k \cdot \ell$ that separates $X_{i}$ from $Y.$
\end{enumerate}

It is noteworthy that in the local version, \zcref{thm_localSegregation}, of \zcref{thm_restrictiveTreewidthIntro}, outcome (a.ii) also relies on the one-to-many strengthening of Menger above to control how well-connected the colored vertices are to a given flat wall in the graph.

\subsection{Related work and organization of the paper}

\paragraph{Related work.}

Various combinatorial problems concerning rooted minors with specific root configurations have been studied by Kawarabayashi \cite{Kawarabayashi2004Rooted}, Wollan 
\cite{Wollan2008Extremal}, Fabila-Monroy and Wood \cite{Fabila-MonroyW2013Rooted}, and Moore \cite{Moore2017Rooted}. For $1$-colorful graphs, structural results have been obtained by Marx, Seymour, and Wollan \cite{MarxSW2017Rooted}; Hodor, La, Micek, and Rambaud \cite{Hodor2024Quickly}; and Fiorini, Kober, Seweryn, Shantanam, and Yuditsky \cite{Fiorini2025Face}.

The most extensively studied problem on $1$-colorful graphs is the \textsc{Steiner Tree} problem, which was shown to be solvable in polynomial time when the torso treewidth is bounded \cite{JansenS2024SteinerTree}, and when a rainbow $K_4$ is excluded \cite{GroendlandNK2024Polynomial}.
To the best of our knowledge, no combinatorial or algorithmic study has been undertaken on $q$-colorful graphs as a general and unifying combinatorial framework.

\paragraph{Organization of the paper.}\label{subsec_Organization}

\zcref{sec_preliminaries} introduces all preliminary concepts we require throughout the paper.
\zcref{subsec_graphMinorStructure} particularly introduces key concepts from graph minor theory for handling the structure of $H$-minor-free graphs.
Moreover, \zcref{subsec_WQO} establishes our well-quasi-ordering result and \zcref{subsec_folio} the reduction for colorful minor checking.
Afterwards, \zcref{sec_rainbowclique}, \zcref{sec_RainbowGrid}, and \zcref{sec_SegregatedGrids} concern our structural results.
\zcref{sec_rainbowclique} deals with the proof of the full version of \zcref{thm_intro_ExcludeRainbowClique}, namely \zcref{thm_excludeRainbowClique}, on the structure of colorful graphs excluding a rainbow clique as a colorful minor.
\zcref{sec_RainbowGrid} proceeds with the full version of \zcref{thm_ExcludingRainbowGridIntro}, namely \zcref{thm_ExcludingRainbowGrid}, on the structure of colorful graphs excluding a rainbow grid as a colorful minor.
The last section on structure results, \zcref{sec_SegregatedGrids}, deals with the structure of graphs excluding segregated grids, in the form of \zcref{thm_segregatedGridGlobal}, which is the full version of \zcref{thm_restrictiveTreewidthIntro}.
Going into \zcref{sec_Applications}, we discuss algorithmic applications of our colorful minor framework and we present our two AMTs, namely \zcref{the_MSOL} and \zcref{the_MSOL_tw}, the full versions of \zcref{thm_rwtMeta_Intro} and \zcref{thm_rcMeta_Intro} respectively. 
In \zcref{sec_concl} we discuss possible directions for future research.
Last but not least, in Appendix \ref{sec_AppendixA} we present a (non-exhaustive) list of problems in which our AMTs apply.

\section{Preliminaries}\label{sec_preliminaries}

By $\mathbb{Z}$ we denote the set of integers and by $\mathbb{R}$ the set of reals.
Given any two integers $a,b\in\mathbb{Z},$ we write $[a,b]$ for the set $\{z\in\mathbb{Z} ~\!\colon\!~ a\leq z\leq b\}.$
Notice that the set $[a,b]$ is empty whenever $a>b.$
For any positive integer $c$ we set $[c]\coloneqq [1,c].$

 \paragraph{Separations and linkages.}
 
Let $G$ be a graph.
A \emph{separation} of $G$ is a tuple $(A,B)$ such that $V(G)=A\cup B$ and there is no edge with one end in $A\setminus B$ and the other in $B\setminus A.$
A set $S\subseteq V(G)$ is called a \emph{separator} of $G$ if $G-S$ has strictly more components than $G.$
Given sets $X,Y\subseteq V(G),$ a set $S$ is an \emph{$X$-$Y$-separator} if no component of $G-S$ contains both a vertex of $X$ and a vertex of $Y.$
Note that, somewhat counterintuitively, this means that an $X$-$Y$-separator in $G$ is not necessarily a separator of $G.$

Let $G$ be a graph and $x,y\in V(G).$
A path $P$ is said to be an \emph{$x$-$y$-path} if the endpoints of $P$ are $x$ and $y.$
Let $X,Y\subseteq V(G).$
A path $P$ is said to be an \emph{$X$-$Y$-path} if there exist $x\in X$ and $y\in Y$ such that $P$ is an $x$-$y$-path and no internal vertex of $P$ belongs to $X\cup Y.$
An \emph{$X$-$Y$-linkage} is a collection $\Pcal$ of pairwise vertex-disjoint $X$-$Y$-paths.
A \emph{linkage} is a collection $\Pcal$ of paths such that there exist $X,Y\subseteq V(G)$ for which $\Pcal$ is an $X$-$Y$-linkage. 
The \emph{order} of a linkage $\Pcal$ is the number $|\Pcal|$ of paths in $\Pcal.$ 
In a slight abuse of notation, we identify a linkage $\Pcal=\{ P_1,P_2,\dots,P_n\}$ and the graph $\bigcup_{i\in[n]}P_i.$

\begin{proposition}[Menger's Theorem]\label{prop_menger}
For every positive integer $k,$ every graph $G,$ and every choice of sets $X,Y\subseteq V(G)$ there either exists an $X$-$Y$-linkage of order $k$ or there exists an $X$-$Y$-separator of size less than $k$ in $G.$

Moreover, there exists an algorithm that finds, given $k,$ $G,$ $X,$ and $Y$ either a linkage of order $k$ or an $X$-$Y$-separator of size less than $k$ in time $\Ocal(k \cdot |\!|G|\!|).$
\end{proposition}

\subsection{Minors}

Let $G$ and $H$ be graphs.
We say that $H$ is a \emph{minor} of $G$ if $H$ can be obtained from $G$ by a sequence of \textsl{(i)} vertex deletions, \textsl{(ii)} edge deletions, and \textsl{(iii)} edge contractions.

It is easy to see that $H$ is a minor of $G$ if and only if there exists a collection $\{ G_v\}_{v\in V(H)}$ of connected and pairwise vertex-disjoint subgraphs of $G,$ called the \emph{branch sets}, such that for every $uv\in E(H)$ there exists an edge $x_ux_v\in E(G)$ with $x_u\in V(G_u)$ and $x_v\in V(G_v).$ 
We call $\{ G_v\}_{v\in V(H)}$ a \emph{minor model} of $H$ in $G.$

Now let $(G,\chi)$ and $(H,\psi)$ be colorful graphs.
We say that a minor model $\{ G_v \}_{v\in V(H)}$ of $H$ in $G$ is a \emph{colorful minor model} of $(H,\psi)$ in $(G,\chi)$ if for every $v\in V(H)$ and every $i\in \psi(v)$ there is some vertex $u\in V(G_v)$ such that $i\in\chi(u).$
We often drop the specifier \say{colorful} and just say that $\{ G_v\}_{v\in V(H)}$ is a minor model of $(H,\psi)$ in $(G,\chi).$

\subsection{Highly linked sets and tangles}

\paragraph{Tangles.}

Let $G$ be a graph and $k$ be an integer.
We denote by $\Scal_k$ the set of all separations of order less than $k$ in $G.$
An \emph{orientation} of $\Scal_k$ is a set $\Ocal\subseteq\Scal_k$ such that for every $(A,B)\in\Scal_k$ \textsl{exactly} one of $(A,B)$ and $(B,A)$ belongs to $\Ocal.$

A \emph{tangle} of order $k$ in a graph $G$ is an orientation $\Tcal$ of $\Scal_k$ such that for all $(A_1,B_1),(A_2,B_2),(A_3,B_3)\in\Tcal$ it holds that $G[A_1]\cup G[A_2]\cup G[A_3]\neq G.$

Let $\Tcal$ and $\Tcal'$ be tangles in $G.$
We say that $\Tcal'$ is a \emph{truncation} of $\Tcal$ if $\Tcal'\subseteq\Tcal.$

Let $\Tcal$ be a tangle of order $k$ in a graph $G$ and let $S\subseteq V(G)$ be a set of size at most $k-2.$
Notice that every separation $(A,B)$ of order at most $k-|S|-1$ in $G-S$ corresponds to a separation $(A\cup S,B\cup S)$ of order at most $k-1$ in $G.$
It follows that $\Tcal_S\coloneqq \{ (A,B) \colon (A\cup S,B\cup S)\in\Tcal \}$ is a tangle of order $k-|S|$ in $G-S.$
Now suppose that $G-S$ has components $G_1,\dots,G_p.$
Then there exists $i\in[p]$ such that $V(G_i)\subseteq B$ for all separations of order $0$ in $\Tcal_S.$
We call $G_i$ the \emph{$\Tcal$-big component} of $G-S.$
\smallskip

 \paragraph{Walls, cliques, and tangles.}
Let $t$ be a positive integer and $\Xcal = \{ G_v\}_{v\in V(K_t)}$ be a minor model of $K_t$ in $G.$
Notice that for every separation $(A,B)\in\Scal_t$ there exists a unique $X\in\{ A,B\},$ $Y\in\{ A,B\}\setminus X,$ such that $V(G_v)\subseteq X\setminus Y$ for some $v\in V(K_t).$
We call $X$ the \emph{$\Xcal$-big side} of $(A,B).$
Moreover, if we let $\Tcal$ be the orientation of $\Scal_t$ obtained by taking all separations $(A,B)\in\Scal_k$ such that $B$ is the $\Xcal$-big side of $(A,B),$ then $\Tcal$ is a tangle.
We call $\Tcal$ the \emph{tangle induced by $\Xcal$} and denote it by $\Tcal_{\Xcal}.$
\smallskip

Let $n,m \in \mathbb{N}$ be two positive integers.
The \emph{$(n \times m)$-grid} is the graph $G$ with the vertex set $V(G) = [n] \times [m]$ and the edges
\begin{align*}
E(G) = \big\{ \{ (i, j) , (\ell , k) \} ~\!\colon\!~ & i, \ell \in [n], \ j,k \in [m], \text{ and } \ \\
 & ( |i - \ell| = 1 \text{ and } j = k ) \text{ or } ( |j - k| = 1 \text{ and } i = \ell ) \big\} .
\end{align*}

The \emph{elementary $(n \times m)$-wall} is in turn derived from the $(n \times 2m)$-grid by deleting all edges in the following set
\begin{align*}
\big\{ \{ (i, j) , (i+1 , j) \} ~\!\colon\!~ i \in [n - 1], \ j \in [m], \text{ and } {i \not\equiv j \mod{2} } \big\}
\end{align*}
and removing all vertices of degree at most 1 in the resulting graph.
An \emph{$(n \times m)$-wall} is a subdivision of the elementary $(n \times m)$-wall.
An \emph{$n$-wall} is an $(n \times n)$-wall and a \emph{wall} is any graph that is an $(n\times m)$-wall for some $n,m.$
Notice that every wall with at least two columns and rows is $2$-connected.
\smallskip

Notice that any $(n\times m)$-wall $W$ consists of $n$ pairwise-disjoint paths $P_1,\dots,P_n$ and $m$ pairwise-disjoint paths $Q_1,\dots,Q_m$ such that $W=\bigcup_{i\in[n]}P_i \cup \bigcup_{i\in[m]}Q_i,$ $P_1\cap Q_j$ is a non-empty path for all $i\in[n],$ $j\in[m],$ each $P_i$ meets the $Q_j$s in order, and each $Q_i$ meets the $P_j$s in order.
We call the $P_i$ the \emph{horizontal paths} of $W$ and the $Q_i$ the \emph{vertical paths} of $W$ (sometimes we also refer to them as the \emph{rows} and \emph{columns} of $W$ respectively).
The \emph{perimeter} of $W$ is the unique cycle contained in $P_1\cup Q_1\cup P_n\cup Q_m$ and every vertex of $W$ not contained in $P_1\cup Q_1\cup P_n\cup Q_m$ is said to be an \emph{inner vertex} of $W.$

A \emph{brick} of a wall $W$ is any facial cycle of $W$ that consists of exactly six degree-$3$ vertices of $W.$

Let $r \in \mathbb{N}$ with $r\geq 3,$ let $G$ be a graph, and $W$ be an $r$-wall in $G.$
Let $\Tcal_{W}$ be the orientation of $\Scal_r$ such that for every $(A,B) \in \Tcal_W,$ the set $B \setminus A$ contains the vertex set of both a horizontal and a vertical path of $W,$ we call $B$ the \emph{$W$-majority side} of $(A,B).$
Then $\Tcal_W$ is the tangle \emph{induced} by $W.$
If $\Tcal$ is a tangle in $G,$ we say that $\Tcal$ \emph{controls} the wall $W$ if $\Tcal_W$ is a truncation of $\Tcal.$

Let $G$ and $H$ be (colorful) graphs as well as $\Tcal$ be a tangle in $G.$
We say that a minor model $\{ G_{v} \}_{v \in V(H)}$ of $H$ in $G$ is \emph{controlled} by $\Tcal$ if there does not exist a separation $(A,B) \in \Tcal$ of order less than $|H|$ and an $x \in V(H)$ such that $V(G_{x}) \subseteq A \setminus B.$
\smallskip

 \paragraph{$k$-linked sets.}
Let $\alpha \in [2/3, 1)_{\mathbb{R}}.$
Moreover, let $G$ be a graph and $X \subseteq V(G)$ be a vertex set. 
A set $S \subseteq V(G)$ is said to be an \emph{$\alpha$-balanced separator} for $X$ if for every component $C$ of $G - S$ it holds that $|V(C) \cap X| \leq \alpha|X|.$ 
Let $k$ be a non-negative integer.
We say that $X$ is a \emph{$(k, \alpha)$-linked set} of $G$ if there is no $\alpha$-balanced separator of size at most $k$ for $X$ in $G.$
\smallskip

Given a $(3k, \alpha)$-linked set $X$ of $G$ we define 
\begin{align*}
\Tcal_{X} \coloneqq \{ (A, B) \in \Scal_{k+1}(G) \colon |X \cap B| > \alpha|X| \}.
\end{align*}
Then $\Tcal_{S}$ is a tangle of order $k+1$ in $G.$

Finally, we need an algorithmic way to identify highly linked sets, and find, given a highly linked set, a large wall whose tangle is a truncation of the tangle induced by the highly linked set.
The first is a classic result due to Reed \cite{Reed1992Finding}.
The second is done in \cite{ThilikosW2024Excluding} by turning to an algorithm a proof of Kawarabayashi, Thomas, and Wollan, from \cite{KawarabayashiTW2020Quickly}.

\begin{proposition}[Reed \cite{Reed1992Finding}]\label{prop_balancedseps}
There exists an algorithm that takes as input an integer $k,$ a graph $G,$ and a set $X\subseteq V(G)$ of size at most $3k+1$ and finds, in time $2^{\Ocal(k)} \cdot |\!|G|\!|$, either a $\nicefrac{2}{3}$-balanced separator of size at most $k$ for $X$ or correctly determines that $X$ is $(k,\nicefrac{2}{3})$-linked in $G.$
\end{proposition}

\begin{proposition}[Thilikos and Wiederrecht \cite{ThilikosW2024Excluding} (see Theorem 4.2.)]\label{thm_algogrid}
Let $k\geq 3$ be an integer and $\alpha\in [2/3,1).$
There exist universal constants $c_1, c_2\in\mathbb{N}\setminus\{ 0\},$ and an algorithm that, given a graph $G$ and a $(c_1k^{20},\alpha)$-linked set $X\subseteq V(G)$ computes in time $2^{\Ocal(k^{c_2})} \cdot |G|^2 |\!|G|\!| \log|G|$ a $k$-wall $W\subseteq G$ such that $\Tcal_W$ is a truncation of $\Tcal_X.$
\end{proposition}

\subsection{Graph minor structure}\label{subsec_graphMinorStructure}

We continue by briefly introducing some key concepts for handling the general structure of $H$-minor-free graphs.
The definitions we introduce here act as a strict simplification of the notions and terminology used in most proofs of the Graph Minor Structure Theorem by Robertson and Seymour.
For more in-depth discussions on this topic, the reader is referred to \cite{RobertsonS2003GraphMinorsXVI,KawarabayashiTW2020Quickly,GorskySW2025Polynomial}.
For the purposes of this paper, the simplified definitions are fully sufficient as we are mostly using the corresponding theorems as black boxes.
Most of the notions below are adapted from \cite{KwarabayashiTW2018NewProof}.

By a \emph{surface} we mean a compact $2$-dimensional manifold with or without boundary.

Given a surface $\Sigma$ and a set $X\subseteq \Sigma$ we denote by $\overline{X}$ the \emph{closure} of $X$ and by $\widetilde{X}$ we denote the set $\overline{X}\setminus X.$

A \emph{painting} in a surface $\Sigma$ is a pair $\Gamma=(U,N),$ where $N\subseteq U\subseteq \Sigma,$ $N$ is finite, $U\setminus N$ has finitely many arcwise-connected components, those are called the \emph{cells} of $\Gamma,$ and for every cell $c,$ the closure $\overline{c}$ is a closed disk where $\widetilde{c}=\overline{c}\cap N\subseteq \mathsf{bd}(\overline{c}).$
If $|\tilde{c}|\geq 4$ we call $c$ a \emph{vortex}.
we define $N(\Gamma)\coloneqq N,$ $U(\Gamma)\coloneqq U,$ and denote the set of cells of $\Gamma$ by $C(\Gamma).$

Notice that the cells of a painting define a hypergraph whose vertices are precisely the elements of $N,$ called the \emph{nodes}, and a node $x\in N(\Gamma)$ belongs to some cell $c\in C(\Gamma)$ if and only if $x\in \widetilde{c}.$

 \paragraph{$\Sigma$-renditions.}
Let $G$ be a graph and $\Sigma$ be a surface.
A \emph{$\Sigma$-rendition} of $G$ is a triple $\rho=(\Gamma,\sigma,\pi),$ where
\begin{itemize}
 \item $\Gamma$ is a painting in $\Sigma,$
 \item $\sigma$ assigns to each cell $c\in C(\Gamma)$ a subgraph $\sigma(c)$ of $G,$ and
 \item $\pi\colon N(\Gamma) \to V(G)$ is an injection
\end{itemize}
such that
\begin{description}
 \item[(R1)] $G=\bigcup_{c\in C(\Gamma)}\sigma(c),$
 \item[(R2)] for all distinct $c,c'\in C(\Gamma),$ $\sigma(c)$ and $\sigma(c')$ are edge-disjoint,
 \item[(R3)] $\pi(\widetilde{c})\subseteq V(\sigma(c))$ for every cell $c\in C(\Gamma),$ and
 \item[(R4)] for every cell $c\in C(\Gamma),$ $V(\sigma(c)\cap \bigcup_{c'\in C(\Gamma)\setminus\{ c\}} (\sigma(c')))\subseteq \pi(\widetilde{c}).$
\end{description}
We write $N(\rho)$ for the set $N(\Gamma)$ and similarly, we lift the set of cells from $C(\Gamma)$ to $C(\rho).$
If the $\Sigma$-rendition $\rho$ of $G$ is understood from the context, we usually identify the sets $\pi(N(\rho))$ and $N(\rho)$ along $\pi$ for ease of notation.
\smallskip

 \paragraph{Vortices, linear decompositions, and depth.}
Let $G$ be a graph, $\Sigma$ be a surface and $\rho=(\Gamma,\sigma,\pi)$ be a $\Sigma$-rendition of $G.$
Moreover, let $c\in C(\rho)$ be a vortex.

Let $X=\{ x_1,x_2,\dots,x_{n-1},x_n\}$ be the vertices of $G$ in $\widetilde{c}$ numbered corresponding to the order of encountering them when traversing along $\mathsf{bd}(\overline{c})$ in anti-clockwise (or clockwise) direction starting from $x_1.$
A \emph{linear decomposition} of $\sigma(c)$ is a sequence $\langle X_1,X_2,\dots,X_{n-1},X_n\rangle$ of sets such that
\begin{description}
 \item[(V1)] $X_i\subseteq V(\sigma(c))$ and $x_i\in X_i$ for all $i\in[n],$
 \item[(V2)] $\bigcup_{i\in[n]}X_i = V(\sigma(c))$ and for all $uv\in E(\sigma(c))$ there exists $i\in[n]$ such that $u,v\in X_i,$ and
 \item[(V3)] for every $v\in V(\sigma(c))$ the set $\{ i\in[n] \colon v\in X_i\}$ forms an interval in $[n].$
\end{description}
The \emph{adhesion} of a linear decomposition is $\max\{|X_i\cap X_{i+1}| \colon i\in[n-1] \}.$
The \emph{width} of a linear decomposition is $\max\{ |X_i| \colon i\in[n] \}.$
The \emph{depth} of a vortex $c\in C(\rho)$ is the smallest integer $k$ such that $\sigma(c)$ has a linear decomposition of adhesion at most $k$ and the \emph{width} of $c$ is the smallest integer $k$ such that $\sigma(c)$ has a linear decomposition of width at most $k.$
\smallskip

The \emph{breadth} of a $\Sigma$-rendition $\rho$ of a graph $G$ is the number of vortices of $\rho,$ the \emph{depth} of $\rho$ is the maximum depth over all vortices of $\rho,$ and finally, the \emph{width} of $\rho$ is the maximum width over all of its vortices.

 \paragraph{Grounded subgraphs, traces, and flat walls.}
Let $\Sigma$ be a surface and $G$ be a graph with a $\Sigma$-rendition $\rho=(\Gamma,\sigma,\pi).$
The \emph{simple torso of $G$ under $\rho$} is the multi graph $G_{\rho}$ obtained from the subgraph $G[N(\rho)]$ by adding an edge between any pair of vertices $u,v\in N(\rho)$ for every non-vortex cell $c\in C(\rho)$ such that $u,v\in\widetilde{c}.$

Notice that $G_{\rho}$ is a graph that can be embedded in $\Sigma$ such that for every $c\in C(\rho)$ with $|\widetilde{c}|=3,$ the vertices of $\widetilde{c}$ form a triangle bounding a face as follows:
For any pair of adjacent vertices $u,v\in N(\rho)$ there exists at least one cell $c\in C(\rho)$ such that $u,v\in \widetilde{c}.$
From here on let $c\in C(\rho)$ be any such cell.
Moreover, there exists a curve $\gamma(c, uv)$ between $u,v \in \Sigma$ that is completely contained in $\overline{c}$ and does not intersect any other element of $N(\rho).$
In case $|\widetilde{c}|=3,$ $\gamma(c,uv)$ is unique.
we apply the same construction for cells $c$ with $|\widetilde{c}|=2$ here, there are two possible choices.
Here, we will always implicitly assume that $\rho$ comes with a unique selection for each such curve $\gamma(c,uv).$

Let $X\subseteq G$ be a cycle or path.
The \emph{image} of $X$ in $G_{\rho}$ is the graph $X_{\rho}\subseteq G_{\rho}$ such that
\begin{itemize}
 \item $V(X_{\rho})=V(X) \cap N(\rho),$ and
 \item $uv\in E(X_{\rho})$ if and only if there exists a subpath of $X$ between $u$ and $v$ that does not contain any other vertex of $X_{\rho}.$
\end{itemize}

Notice that any subpath between $u$ and $v$ in $X$ such as the one from the second item of the above definition must be fully contained in a single cell $c$ of $\rho.$
Hence, for any subpath $P$ of $X$ between vertices of $X_{\rho}$ that contains no other element of $V(X_{\rho})$ there exists a unique curve $\gamma(c,uv)$ \emph{mirroring} $P.$
We may identify the edges of $X_{\rho}$ and the curves mirroring the subpaths of $X$ corresponding to those edges.
For every edge $uv\in E(X_{\rho})$ we say that the curve $\gamma(c,uv)$ mirroring $P$ is the \emph{drawing} of $uv.$

If $X$ is a path, we say that $X$ is grounded if $X_{\rho}$ has at least one edge.
In case $X$ is a cycle we say that $X$ is grounded if $X_{\rho}$ is a simple\footnote{That is, if $X_{\rho}$ has at least three vertices and no parallel edges.} cycle.
If $X$ is grounded, we call the subset of $\Sigma$ that is the union of $V(X_{\rho})$ together with the union of the drawings of all edges of $X_{\rho}$ the \emph{trace} of $X.$
Notice that the trace of $X$ is a simple curve in case $X$ is a path and a closed curve in case $X$ is a cycle.

A $2$-connected subgraph $H$ of $G$ is said to be \emph{grounded} if every cycle of $H$ is grounded.
\smallskip

As before let $\Sigma$ be a surface and $G$ be a graph with a $\Sigma$-rendition $\rho=(\Gamma,\sigma,\pi).$
Now let $W$ be an $(n\times m)$-wall, $n,m\geq 3,$ in $G$ such that $W$ is grounded in $\rho.$
We say that $W$ is \emph{flat in $\rho$} if there exists a disk $\Delta\subseteq \Sigma$ such that $\mathsf{bd}(\Delta)$ intersects $\Gamma$ only in $N(\Gamma),$ no vortex of $\rho$ is contained in $\Delta,$ all vertices of $W_{\rho}$ are contained in $\Delta,$ and the trace of the perimeter of $W$ separates all inner vertices of $W_{\rho}$ from $\mathsf{bd}(\Delta).$

We say that a wall $W$ is \emph{flat} in $G$ if there exists a $\Sigma$-rendition $\rho$ of $G$ for some non-empty surface $\Sigma$ such that $W$ is flat in $\rho.$

The following theorem is known as the \emph{Flat Wall Theorem}.
It first appeared in a weaker form and without explicit bounds in Graph Minors XIII by Robertson and Seymour \cite{RobertsonS1995Graph}.
The best known bounds are due to Chuzhoy \cite{Chuzhoy2015Improved}.
The version below is due to Gorsky, Seweryn, and Wiederrecht \cite{GorskySW2025Polynomial}.

\begin{proposition}[Flat Wall Theorem \cite{RobertsonS1995Graph,GorskySW2025Polynomial}]\label{thm_flatwall}
 Let $t, r$ be integers with $t \geq 5,$ $r \geq 3,$ let $n \geq 100t^3(r+2t+2),$ and let $G$ be a graph with an $n$-wall $W.$
 Then either
 \begin{itemize}
 \item there exists a model of $K_t$ in $G$ which is controlled by $\mathcal{T}_{W},$ or
 \item there exist a set $Z$ with $|Z| < 16t^3$ and an $r$-subwall $W'$ of $W$ which is disjoint from $Z$ and flat in $G - Z.$
 \end{itemize}
 Furthermore, there exists a $\poly(t+r) \cdot |\!|G|\!|$-time algorithm which finds either the model of $K_t,$ or the set $Z,$ the subwall $W',$ and a $\Sigma$-rendition witnessing that $W'$ is flat in $G - Z.$
\end{proposition}

A slightly weaker but relevant version of grounded walls is the notion of centrality.
Let $\Sigma$ be a surface and let $\rho$ be a $\Sigma$-rendition of a graph $G$ containing an $r$-wall $W.$
We say that $\delta$ is \emph{$W$-central} if there is no cell $c\in C(\delta)$ such that $V(\sigma(c))$ contains the $W$-majority side of a separation from $\Tcal_W.$
Similarly, let $A\subseteq V(G),$ $|A|\leq r-1,$ let $\Sigma'$ be a surface and $\rho'$ be a $\Sigma'$-rendition of $G-A.$
Then we say that $\rho'$ is \emph{$(W-A)$-central} for $G$ if no cell of $\delta'$ contains the majority side of a separation from $\Tcal_W\cap\Scal_{r-|A|}.$

With this, we are ready to state the so called Local Structure Theorem.
This theorem was originally proven by Robertson and Seymour \cite{RobertsonS2003GraphMinorsXVI}, the version we present here is the recent result of Gorsky, Seweryn, and Wiederrecht \cite{GorskySW2025Polynomial} that also allows for polynomial bounds.

\begin{proposition}[Local Structure Theorem \cite{RobertsonS2003GraphMinorsXVI,GorskySW2025Polynomial}]\label{thm_localstructure}
There exist functions $\apex_{\ref{thm_localstructure}},\depth_{\ref{thm_localstructure}}\colon\mathbb{N}\to\mathbb{N}$ and $\wall_{\ref{thm_localstructure}}\colon\mathbb{N}^2\to\mathbb{N}$ such that for all integers, $t\geq 5,$ and $r\geq 3,$ every graph $H$ on $t$ vertices, every graph $G$ and every $\wall_{\ref{thm_localstructure}}(t,r)$-wall $W\subseteq G$ one of the following holds.
\begin{enumerate}
 \item $G$ has an $H$-minor controlled by $\mathcal{T}_{W},$ or
 \item there exists a set $A\subseteq V(G)$ of size at most $\apex_{\ref{thm_localstructure}}(t)$ and a subwall $W'\subseteq W-A$ such that $G$ has a $\Sigma$-rendition $\rho$ of breadth at most $\nicefrac{1}{2}(t-3)(t-4)$ and depth at most $\depth_{\ref{thm_localstructure}}(t)$ such that $\Sigma$ is a surface where $H$ does not embed, $W'$ is flat in $\rho,$ and $\rho$ is $(W-A)$-central.
\end{enumerate}
Moreover, it holds that $\apex_{\ref{thm_localstructure}}(t),~ \depth_{\ref{thm_localstructure}}(t) \in \Ocal\big(t^{112}\big), \text{ and } \wall_{\ref{thm_localstructure}}(t,r) \in \Ocal\big(t^{115}r \big).$
There also exists an algorithm that, given $t,r,$ a graph $H,$ a graph $G$ and a wall $W$ as above as input finds one of these outcomes in time $\poly(t+r) \cdot |G|^2 |\!|G|\!|.$
\end{proposition}

\subsection{Well-Quasi Ordering}\label{subsec_WQO}

The main goal of this subsection is to prove that for fixed $q,$ the colorful minor relation on colorful graphs on at most $q$ colors is a well-quasi-order.
This result follows directly from a powerful general theorem of Robertson and Seymour \cite{GraphMinorsXXIII}.

In this subsection, and this subsection only, we permit our graphs to have loops.

A directed graph $D$ is an \emph{oriented graph} if there exists a graph $G$ such that $V(D)=V(G)$ and $D$ can be obtained from $G$ by replacing each undirected edge $\{ u,v\}\in E(G)$ by one of the two possible directed edges $(u,v)$ or $(v,u).$

\begin{proposition}[Robertson and Seymour \cite{GraphMinorsXXIII}]\label{prop_superWQO}
Let $\Lambda=(E_{\Lambda},\leq_{\Lambda})$ be a well-quasi-order and $\langle G_i\rangle _{i\in\mathbb{N}}$ be a sequence of oriented graphs.
Moreover, for each $i\in\mathbb{N}$ let $\phi_i\colon V(G_i)\cup E(G_i)\to E_{\Lambda}$ be some function.
Then there exist $j > i\geq 1$ and a map $\eta$ with domain $V(G_i)\cup E(G_i),$ satisfying:
\begin{enumerate}
 \item for each $v\in V(G_i),$ $\eta(v)$ is a connected subgraph of $G_j,$ and there exists $x\in V(\eta(v))$ with $\phi_i(v) \leq_{\Lambda} \phi_j(x)$; and $\eta(v)\cap\eta(w)=\emptyset$ for all pairs of distinct $v,w\in V(G_i),$
 \item for each $e\in E(G_i),$ $\eta(e)$ is an edge of $G_j,$ a loop if and only if $e$ is a loop, with $\phi_i(e)\leq_{\Lambda} \phi_j(\eta(e))$; and $\eta(e)\neq\eta(e')$ for all distinct pairs $e,e'\in E(G_i),$ and
 \item for each $e\in E(G_i)$ with head $u$ and tail $v,$ $\eta(e)$ has head in $V(\eta(u))$ and tail in $V(\eta(v)),$ and $\eta(e)$ is not an edge of $\eta(u)$ or $\eta(v)$ (the last is trivial unless $u=v$).
\end{enumerate}
\end{proposition}

This is all we need to deduce \zcref{thm_WQO}.

\begin{proof}[Proof of \zcref{thm_WQO}.]
We consider the partial order $\Lambda = (2^{[q]},\subseteq).$
Since $q$ is a fixed integer, $2^{[q]}$ is finite and thus, $\Lambda$ is a well-quasi-order.

Suppose there exists an infinite but countable anti-chain $\langle (G_i,\chi_i) \rangle_{i\in\mathbb{N}}$ of colorful graphs on at most $q$ colors under the colorful minor relation.
For each $i\in\mathbb{N}$ we define $\phi_i(v)\coloneqq \chi_i(v)$ for all $v\in V(G_i)$ and $\phi_i(e)\coloneqq \emptyset$ for all $e\in E(G_i).$
Moreover, for each $i\in V(G_i)$ we may assign an arbitrary orientation to the edges of $G_i$ to obtain an oriented graph.

Then there must exist $j> i\geq 1$ such that $G_i$ and $G_j$ satisfy the outcome of \zcref{prop_superWQO} witnessed by the map $\eta.$
Now observe that $\eta$ also witnesses that $G_i$ is a minor of $G_j$ and, in particular, $\chi_i(v)\subseteq \chi_j(\eta(v))$ for all $v\in V(G_i).$
Hence, $(G_i,\chi_i)$ is a colorful minor of $(G_j,\chi_j)$ which contradicts the assumption that $\langle (G_i,\chi_i) \rangle_{i\in\mathbb{N}}$ is an anti-chain.
\end{proof}

 \paragraph{Colorful obstructions.}
Let $q$ be a non-negative integer and let $\Ccal$ be 
a class of $q$-colorful graphs that is closed under taking colorful minors.
We define the (colorful minor) \emph{obstruction set} of $\Ccal,$ denoted by $\obs(\Ccal),$
as the set of all colorful minor minimal $q$-colorful graphs that are not members 
of $\Ccal.$ Notice that the obstruction set of $\Ccal$ is an anti-chain for the 
colorful minor relation, therefore, from \zcref{thm_WQO} we have the following.

\begin{corollary}
\label{cor_finite_set}
Let $q$ be a non-negative integer and let $\Ccal$ be 
a class of $q$-colorful graphs that is closed under taking colorful minors.
Then $\obs(\Ccal)$ is a finite set.
\end{corollary}

\subsection{Colorful minor checking}\label{subsec_folio}

We now proceed with the proof of \zcref{thm_introColorfulMinors}.
Having established that the colorful minor relation is indeed a well-quasi-order, we next show that testing whether some $q$-colorful graph $(H,\psi)$ is a colorful minor of a $q$-colorful graph $(G,\chi)$ can be done efficiently.
Towards proving \zcref{thm_introColorfulMinors} we require the fast minor testing algorithm of Korhonen, Pilipczuk, and Stamoulis \cite{KorhonenPS2024Minor}.
Indeed, we need the rooted variant of their algorithm.

A tuple $(G,r_1,\dots,r_{\ell})$ where $r_i\in V(G)$ for all $i\in[\ell]$ is called a \emph{rooted graph} and the $r_i$ are called the \emph{roots} of $(G,r_1,\dots,r_{\ell})$.
Notice that one may interpret a rooted graph with $\ell$ roots as an $\ell$-colorful graph $(G,\rho)$ where $|\rho^{-1}(i)|=1$ for each $i\in[\ell]$.
Indeed, to each rooted graph $\mathbf{G}=(G,r_1,\dots,r_{\ell})$ we associate a corresponding colorful graph $(G,\rho_{\mathbf{G}})$ such that $\rho_{\mathbf{G}}^{-1}(i)=r_i$ for all $i\in[\ell]$.
We say that a rooted graph $\mathbf{H}$ is a \emph{rooted minor} of a rooted graph $\mathbf{G}$ if $(H,\rho_{\mathbf{H}})$ is a colorful minor of $(G,\rho_{\mathbf{G}})$.

\begin{proposition}[Korhonen, Pilipczuk, Stamoulis \cite{KorhonenPS2024Minor}]\label{prop_minorChecking}
There exists a function $f_{\ref{prop_minorChecking}}\colon \mathbb{N}^2 \to \mathbb{N}$ and an algorithm that takes as input two rooted graph $\mathbf{H}=(H,t_1,\dots,t_{p})$ and $\mathbf{G}=(G,r_1,\dots,r_q)$ and decides if $\mathbf{H}$ is a rooted minor of $\mathbf{G}$ in time $\mathcal{O}(f_{\ref{prop_minorChecking}}(|H|,p))|\!|G|\!|^{1+o(1)}$.
\end{proposition}

The proof of \zcref{thm_introColorfulMinors} is a simple reduction of the colorful minor checking problem to the rooted minor checking problem.

Let $(G,\chi)$ be a $q$-colorful graph.
We denote by $F(G,\chi)$ the rooted graph $(G,(t_{v,i} \colon v\in V(G),~ i\in\chi(v)))$ obtained from $G$ by introducing for every $v\in V(G)$ and $i\in[q]$ a new vertex $t_{v,i}$ adjacent precisely to $v$.

Moreover, let $(H,\psi)$ be another $q$-colorful graph.
We denote by $R_{(H,\psi)}(G,\chi)$ the rooted graph $(G,(r_{v,i} \colon v\in V(H),~ i\in\psi(v)))$ obtained from $G$ by introducing, for every $v\in V(H)$ and $i\in\psi(v)$ a new vertex $r_{v,i}$ adjacent to all vertices of $\chi^{-1}(i)$.

Notice that for a fixed $q$-colorful graph $(H,\psi)$, the number of roots in $F(H,\psi)$ and $R_{(H,\psi)}(G,\chi)$ are the same.
As a convention, we will assume the vertices of $H$ to come with some linear ordering such that in the sequence $(t_{v,i})$ of roots of $F(H,\psi)$, the vertex $t_{u,j}$ has precisely the same position as the vertex $r_{u,j}$ in the sequence $(r_{v,i})$ of roots of $R_{(H,\psi)}(G,\chi)$.
Under this convention, we now obtain the following lemma.

\begin{lemma}\label{lemma_emulateColours}
Let $(H,\psi)$ and $(G,\chi)$ be $q$-colorful graphs.
Then $(H,\psi)$ is a colorful minor of $(G,\chi)$ if and only if $F(H,\psi)$ is a rooted minor of $R_{(H,\psi)}(G,\chi)$.
\end{lemma}

\begin{proof}
Assume first that there exists a colorful minor model $\{ G_v\}_{v\in V(H)}$ of $(H,\psi)$ in $(G,\chi)$.
In this case, for each $v\in V(H)$ and $i\in \psi(v)$ there exists a vertex $x\in V(G_v)$ with $i \in \chi(v)$.
Hence, in $R_{(H,\psi)}(G,\chi)$, $x$ is adjacent to $r_{v,i}$.
It follows that $\{ G_v\}_{v\in V(H)} \cup \{ r_{v,i}\}_{v\in V(H),~i\in\chi(v)}$ is a minor model witnessing that $F(H,\psi)$ is indeed a rooted minor of $R_{(H,\psi)}(G,\chi)$.

Now assume that there exists a minor model $(G_v)_{v\in V(F(H,\psi))}$ of $F(H,\psi)$ in $R_{(H,\psi)}(G,\chi)$.
In this case, for each $v\in V(H)$, $G_v$ must first of all be disjoint from all roots of $R_{(H,\psi)}(G,\chi)$ as each $r_{v,i}$ must be mapped by the minor model to the root $t_{v,i}$ of $F(H,\psi)$.
Moreover, for each $i\in \psi(v)$, there is a vertex $x$ of $G_v$ adjacent to $r_{v,i}$.
It follows from the definition of $R_{(H,\psi)}$ that $i\in \chi(x)$.
Hence, $\{ G_v \}_{v\in V(H)}$ must be a colorful minor model of $(H,\psi)$ in $(G,\chi)$.
\end{proof}

With \zcref{lemma_emulateColours} it now suffices to observe that both $F(H,\psi)$ and $R_{(H,\psi)}(G,\chi)$ can be computed from $(H,\psi)$ and $(G,\chi)$ in time linear in $q \cdot (|H| + |G|)$.
Moreover, the number of roots of $R_{(H,\psi)}(G,\chi)$ is at most $q|H|$ and this is also an upper bound on the number of vertices $F(H,\psi)$ has in excess of $H$ itself.
Finally, the number of edges in $R_{(H,\psi)}(G,\chi)$ is at most $|\!|G|\!| + |G|\cdot |H| \cdot q$.
So we may compute $R_{(H,\psi)}(G,\chi)$ and run the algorithm from \zcref{prop_minorChecking} to check for $F(H,\psi)$ as a rooted minor without any additional impact on its running time in terms of $H$ and $G$.

\section{Excluding a rainbow clique}\label{sec_rainbowclique}

In this section we discuss the structure of colorful graphs excluding a fixed rainbow clique as a colorful minor.
For this, first we establish a general tool that allows us to connect several sets of vertices $X_i$ to a fixed target set $Y$ with many disjoint paths, or to find a small separator that separates at least one of the $X_i$ from $Y.$
From here, we make use of several useful results of Robertson and Seymour from \cite{RobertsonS1995Graph} to establish the following key theorem of this section.

\begin{theorem}\label{thm_rainbowclique}
Let $q,$ $t,$ and $k$ be positive integers with $k\geq \lfloor \nicefrac{3}{2}\cdot qt\rfloor+t.$
Let $(G,\chi)$ be a $q$-colorful graph such that $G$ contains a minor model $\Xcal$ of $K_k.$
Then one of the following is true.
\begin{enumerate}
\item There exists a colorful minor model $\Fcal$ of a rainbow $K_t$ in $(G,\chi)$ such that $\Tcal_{\Fcal}$ is a truncation of $\Tcal_{\Xcal},$ or
\item There exists a set $S\subseteq V(G)$ of size at most $qt-1$ such that the $\Tcal_{\Xcal}$-big component of $G-S$ is restricted.
\end{enumerate}
Furthermore, there exists an algorithm that takes as input $t,$ $(G,\chi),$ and $\Xcal$ as above and finds one of the two outcomes in time $\poly(q+t) \cdot |G|.$
\end{theorem}

Establishing \zcref{thm_rainbowclique} is the key towards our first structure theorem for colorful graphs excluding a rainbow clique minor, i.e.\@ \zcref{thm_intro_ExcludeRainbowClique}.
Robertson and Seymour in \cite{RobertsonS1995Graph}, prove a theorem akin to \zcref{thm_rainbowclique} that however works in the setting of $1$-colorful graphs with a bounded number of colored vertices.
In the following section, we adapt this result to work for $q$ colors and an arbitrary number of colored vertices to suit our needs.

\subsection{Finding a rainbow clique}\label{subsec_rainbowclique}

Within our proofs, in this section and beyond, we will repeatedly make use of a variant of \zcref{prop_menger} that allows us to link a family of a bounded number of vertex sets $X_1,\dots,X_{\ell}$ equally into a single set $Y,$ or separate at least one of the $X_i$ from $Y$ by deleting a small set of vertices.

\begin{lemma}\label{lemma_multicolorlinakge}
	For all positive integers $k,\ell,$ every graph $G$ and every collection $X_1,\dots,X_{\ell},Y\subseteq V(G)$ of (not necessarily disjoint) vertex sets one of the following holds:
	\begin{enumerate}
		\item there exists a linkage $\Pcal$ of order $k\cdot \ell$ in $G$ such that for every $i\in[\ell]$ there exists an $X_i$-$Y$-linkage $\Pcal_i\subseteq\Pcal$ of order $k$ and $\Pcal=\bigcup_{i\in[\ell]}\Pcal_i,$ or
		\item there exists a non-empty set $I\subseteq [\ell]$ and a set $S\subseteq V(G)$ of order less than $k\cdot\ell$ such that $S$ is a $\bigcup_{i\in I}X_i$-$Y$-separator in $G.$ 
	\end{enumerate}
	Moreover, there exists an algorithm that, given $k,\ell,$ $G,$ $X_1,\dots,X_{\ell},$ and $Y$ as above as input finds one of the two outcomes in time $\poly(k+\ell) \cdot |G|.$
	\end{lemma}
	
	\begin{proof}
	Let $k,\ell,$ $G,$ $X_1,\dots,X_{\ell},$ and $Y$ be given.
	We construct an auxiliary graph $G'$ as follows.
	For every $i\in[\ell]$ we introduce a set $Z_i$ of $k$ new vertices and for each $z\in Z_i$ we make $z$ adjacent to all vertices of $X_i.$
	We then apply \zcref{prop_menger} to the graph $G$ with sets $Z=\bigcup_{i\in[\ell]}Z_i$ and $Y$ to either find a $Z$-$Y$-linkage $\Pcal'$ of order $k\cdot\ell$ or a $Z$-$Y$-separator $S$ of order at most $k\cdot\ell-1$ in $G.$
	
	Let us first assume the case where \zcref{prop_menger} returns the linkage $\Pcal'.$
	Then every vertex of $Z$ must be an endpoint of some path in $\Pcal'.$
	For each $i\in[\ell]$ let $\Pcal_i$ the collection of paths $P$ obtained from a path $P'\in\Pcal'$ with one endpoint in $Z_i$ by deleting this endpoints from $P'.$
	It follows that $|\Pcal_i|=k$ and $\Pcal_i$ is an $X_i$-$Y$-linkage.
	Hence, this case produces a linkage as desired.
	
	So now assume \zcref{prop_menger} returns a $Z$-$Y$-separator $S$ of order at most $k\cdot\ell-1$ in $G'.$
	Let $I\subseteq [\ell]$ be the set of all $i\in[\ell]$ such that there is no $X_i$-$Y$-path in $G-(S\setminus Z).$
	All we have to do is to show that $I$ is non-empty.
	Towards this goal let $J\subseteq[\ell]$ be the collection of all $i\in[\ell]$ such that $Z_i\setminus S$ is non-empty.
	Since $|Z|=k\cdot\ell$ and $|S|\leq k\cdot\ell-1$ we know that $J\neq \emptyset.$
	Now observe that there cannot be a $j\in J$ such that there is an $X_j$-$Y$-path in $G-(S\setminus Z).$
	To see this assume that such a path, let us call it $R,$ exists for some $j\in J.$
	Then this path must also exist in $G'-S.$
	However, since there is some $z\in Z_j\setminus S$ and $R$ has one endpoint in $X_j,$ the path $zR$ is a $Z_j$-$Y$-path in $G'.$
	This contradicts the assumption that $S$ is a $Z$-$Y$-separator in $G'.$
	Thus, the set $S\setminus Z$ is the set required by the second outcome of our assertion and the proof is complete. 
	\end{proof}
	
	We are going to need the following theorem of Robertson and Seymour \cite{RobertsonS1995Graph}.

	\begin{proposition}[Robertson and Seymour \cite{RobertsonS1995Graph}]\label{prop_linkaclique}
		Let $k$ and $p$ be positive integers with $p\geq \lfloor \nicefrac{3}{2}\cdot k\rfloor.$
		Further, let $G$ be a graph, $Z\subseteq V(G)$ of size $k,$ and $G_1,\dots,G_p$ be pairwise vertex-disjoint subgraphs of $G,$ such that
		\begin{enumerate}
			\item for every $i\in[p],$ either $G_i$ is connected, or every component of $G_i$ has a vertex in $Z,$
			\item for all $i\neq j\in[p],$ either both $G_i$ and $G_j$ contain a vertex of $Z,$ or there exists an edge $x_ix_j\in E(G)$ with $x_i\in V(G_i)$ and $x_j\in V(G_j),$ and
			\item there is no $i\in[p]$ such that there is a separation $(A,B)$ of $G$ of order less than $k$ with $Z\subseteq A$ and $A\cap V(G_i)= \emptyset.$
		\end{enumerate}
		Then there exists $h\in[0,k]$ together with $p'=p-\lfloor\nicefrac{1}{2}(k-h)\rfloor\geq k$ connected and pairwise vertex-disjoint subgraphs $H_1,\dots,H_{p'}$ of $G,$ such that
		\begin{itemize}
			\item $|V(H_i)\cap Z|=1$ for each $i\in[k]$ and $V(H_i)\cap Z=\emptyset$ for all $i\in[k+1,p'],$ and
			\item for all $i\neq j\in[p'],$ if there does not exist an edge $x_ix_j\in E(G)$ with $x_i\in V(H_i)$ and $x_j\in V(H_j)$ then $i,j\in[k],$ i.e. both $H_i$ and $H_j$ have a non-empty intersection with $Z.$
		\end{itemize}
	\end{proposition}

Notice that \zcref{prop_linkaclique} can be realized in polynomial time, allowing for an algorithmic version.
From here on, all of our proofs are constructive and thus, all of them admit polynomial-time algorithms that yield the respective outcomes.

The following does not find a direct application in this paper, but it makes for a nice observation on its own.

	\begin{observation}\label{obs_homogeniseclique}
		Let $t$ and $q$ be positive integers.
		Let $(G,\chi)$ be a $q$-colorful graph such that there exist $p\coloneqq 2^qt$ pairwise vertex-disjoint connected subgraphs $G_1,\dots,G_p$ such that for all $i\neq j\in[p]$ there exists $x_ix_j\in E(G)$ with $x_i\in V(G_i)$ and $x_j\in V(G_j).$
		Then there exist $I\subseteq [p]$ and $J\subseteq[q]$ such that
		\begin{itemize}
			\item $|I|\geq t,$ and
			\item for all $i\in I$ it holds that $J=\chi(G_i).$
		\end{itemize}
	\end{observation}

	\begin{proof}
		The lemma follows immediately from an iterative application of the pigeon hole principle as follows.
		Let $I_0\coloneqq [p].$
		For each $i\in[q],$ assuming inductively that $I_{i-1}$ with $|I_{i-1}|\geq 2^{q-i+1}t$ has been found, let $I^1_i\coloneqq \{ j\in I_{i-1} \colon i\in\chi(G_j) \}$ and $I^2_i\coloneqq \{ j\in I_{i-1} \colon i\notin\chi(G_j) \}.$
		Then there exists $r\in\{ 1,2\}$ such that $|I^r_i|\geq 2^{q-1}t,$ if this is true for both set $r\coloneqq 1.$
		Let $I_i\coloneqq I^r_i.$
		Then, by construction, the set $I_q$ has the desired properties.
	\end{proof}

\begin{lemma}\label{lemma_linkchromaticclique}
Let $t$ and $q$ be positive integers.
Let $(G,\chi)$ be a $q$-colorful graph such that there exist $p\geq \lfloor \nicefrac{3}{2}\cdot tq\rfloor+t$ connected and vertex-disjoint subgraphs $G_1,\dots,G_p$ of $G,$ such that for every $i\neq j\in[p]$ there exists an edge $x_ix_j\in E(G)$ such that $x_i\in V(G_i)$ and $x_j\in V(G_j).$
Then one of the following is true.
\begin{enumerate}
\item there exists a separation $(A,B)$ of $G$ of order less than $tq,$ a non-empty set $I\subseteq [q]$ such that $\chi^{-1}(I)\subseteq A,$ and there exists $j\in[p]$ such that $V(G_j)\cap A=\emptyset,$ or
\item for $p'=(q+1)t,$ there exist pairwise vertex-disjoint subgraphs $H_1,\dots,H_{p'}$ such that 
\begin{itemize}
\item for every $i\in[q]$ and $j\in[(i-1)t+1,it]$ it holds that $i\in\chi(G_j),$
\item for all $i\in[p']$ either $H_i$ is connected or there exists $j\in[q]$ such that $i\in[(j-1)t+1,jt]$ and every component of $H_i$ contains a vertex $v$ with $j\in\chi(v),$ and
\item for every $i\neq j\in[p']$ there either exists an edge $x_ix_j\in E(G)$ such that $x_i\in V(G_i)$ and $x_j\in V(G_j)$ or $i,j\in[tq].$ 
\end{itemize}
\end{enumerate}
\end{lemma}
\begin{proof}
Let us begin by creating an auxiliary graph $G'$ as follows.
For every $i\in[q]$ introduce a set $Z_i$ of $t$ new vertices such that each of them is adjacent to every vertex $v\in V(G)$ with $i\in\chi(v).$
Let $Z\coloneqq\bigcup_{i\in[q]}Z_i.$
Then $|Z|=qt.$

Now assume that there exists a separation $(A',B')$ of $G'$ of order at most $tq-1$ such that $Z\subseteq A$ and there is some $i\in[p]$ such that $V(G_i)\cap A=\emptyset.$
Let $A\coloneqq A'\cap V(G)$ and let $B\coloneqq B'\cap V(G).$
We claim that there exists some $i\in[q]$ such that $\chi^{-1}(i)\subseteq A$ and $(A,B)$ is a separation of order at most $tq-1$ of $G.$
First notice that, if there is an edge $ab\in E(G)$ with $a\in A\setminus B$ and $B\setminus A$ then this edge also exists in $G'.$
Hence, $(A,B)$ is indeed a separation of $G$ of the desired order.
Next observe that, since $Z\subseteq A',$ $|Z|=tq,$ and $|A'\cap B'|\leq tq-1$ there must exist some $i\in[q]$ such that $Z_i\setminus B'\neq\emptyset.$
It follows from the definition of separations that there cannot be a vertex $v\in B\setminus A$ such that $i\in\chi(v)$ as otherwise $v$ would be adjacent to a vertex in $A'\setminus B'.$
Thus, in this case the first outcome of our assertion holds.
		
From here on we may assume that there does not exist a separation $(A',B')$ in $G'$ such that $Z\subseteq A'$ and there is some $i\in[p]$ such that $V(G_i)\cap A'\neq\emptyset.$
In this situation we may apply \zcref{prop_linkaclique} to $G',$ $Z,$ and the graphs $G_1,\dots,G_{p}.$
It follows that there exists some $p'\in[(q+1)t,p]$ together with connected and pairwise vertex-disjoint subgraphs $H'_1,\dots,H'_{p'}$ of $G'$ such that
\begin{itemize}
\item $|V(H'_i)\cap Z|=1$ for each $i\in[tq]$ and $V(H_i') \cap Z = \emptyset$ for all $i \in [tq + 1, p'],$ and
\item for all $i \neq j\in[p'],$ if there does not exist an edge $x_ix_j\in E(G')$ with $x_i\in V(H'_i)$ and $x_j\in V(H'_j)$ then $i,j\in[tq].$
\end{itemize}
By possibly reordering we may assume that for every $i\in [q]$ and every $j\in[(i-1)t,it]$ it holds that $Z_i\cap V(H'_i)\neq\emptyset.$
It then follows from the definition of $G'$ that $i\in\chi(V(H_j'))$ for each $j\in[(i-1)t,it].$
For every $i\in[p']$ let $H_i\coloneqq H'_i-Z.$
Now, for every $i\in[q]$ and $j\in[(i-1)t+1,it]$ it holds that either $H_i$ is connected, or every component of $H_i$ contains a vertex that is adjacent to a vertex of $Z_i$ in $G'.$
It follows that this vertex carries color $i.$
Finally, the last point of outcome two of our assertion follows directly from the definition of the $H_i$ and our proof is complete.
\end{proof}

\begin{proof}[Proof of \zcref{thm_rainbowclique}.]
Let $\Xcal=\{ G_i\}_{i\in[k]}.$
We may now apply \zcref{lemma_linkchromaticclique}.
		
Suppose \zcref{lemma_linkchromaticclique} returns a separation $(A,B)$ of $G$ of order at most $qt-1$ together with a non-empty set $I\subseteq[q]$ such that $\chi^{-1}(I)\subseteq A$ and there is some $i\in[k]$ for which $V(G_i)\cap A=\emptyset.$
It follows that $G[B\setminus A]$ is restricted.
Moreover, the unique $\Tcal_{\Xcal}$-big component of $G-(A\cap B)$ must be fully contained in $B\setminus A.$
Hence, with $S\coloneqq A\cap B$ we have found the desired set.

Hence, we may assume that \zcref{lemma_linkchromaticclique} returns $p'$ pairwise vertex-disjoint subgraphs $H_1,\dots,H_{p'}$ where $p'=(q+1)t$ and
\begin{itemize}
\item for every $i\in[q]$ and $j\in[(i-1)t+1,it]$ it holds that $i\in\chi(G_j),$
\item for all $i\in[p']$ either $H_i$ is connected or there exists $j\in[q]$ such that $i\in[(j-1)t+1,jt]$ and every component of $H_i$ contains a vertex $v$ with $j\in\chi(v),$ and
\item for every $i\neq j\in[p']$ there either exists an edge $x_ix_j\in E(G)$ such that $x_i\in V(G_i)$ and $x_j\in V(G_j)$ or $i,j\in[tq].$ 
\end{itemize}
Now, for each $i\in[t]$ let us define a connected graph $F_i$ as follows.
For each $j\in [q]$ there exists a component $D_{i,j}$ of $H_{(j-1)t+i}$ such that $j\in\chi(D_{i,j})$ and there exists an edge $uv\in E(G)$ with $u\in V(D_{i,j})$ and $v\in H_{qt+i}.$
We set $F_i\coloneqq H_{qt+i}\cup \bigcup_{j\in[q]}D_{i,j}.$
It is now easy to see that $\Fcal=\{ F_i\}_{i\in[t]}$ is indeed a colorful minor model of a rainbow $K_t.$
\end{proof}

\subsection{Colorful graphs without large rainbow clique minors}\label{subsec_cliqueexclusionthm}

We now have everything in place to prove \zcref{thm_intro_ExcludeRainbowClique}.
In fact, we will prove a slightly stronger statement which will then imply \zcref{thm_intro_ExcludeRainbowClique}.

\begin{theorem}\label{thm_excludeRainbowClique}
There exist functions $\adhesion_{\ref{thm_excludeRainbowClique}},\genus_{\ref{thm_excludeRainbowClique}},\apex_{\ref{thm_excludeRainbowClique}},\breadth_{\ref{thm_excludeRainbowClique}},\vortex_{\ref{thm_excludeRainbowClique}}\colon\mathbb{N}^2\to\mathbb{N}$ such that for all integers $q,t\geq 1$ and every $q$-colorful graph $(G,\chi)$ one of the following holds:
\begin{enumerate}
 \item $(G,\chi)$ contains a rainbow $K_t$ as a colorful minor, or
 \item $(G, \chi)$ has a tree-decomposition $(T, \beta)$ of adhesion at most $\adhesion_{\ref{thm_excludeRainbowClique}}(q, t)$ and a (possibly empty) subset $L \subseteq V(T)$ of the leaves of $T$ such that for every $x \in V(T)$ one of the following holds:
 \begin{itemize}
 \item $x \in L$ is a leaf with parent $d\in V(T)$ and $(G[\beta(x) \setminus \beta(d)], \chi)$ is restricted, or
 \item if $(G_{x}, \chi_{x})$ is the colorful torso of $G$ at $x,$ then there exists a set $A \subseteq V(G_{x})$ of size at most $\apex_{\ref{thm_excludeRainbowClique}}(q,t)$ and a surface of Euler-genus at most $\genus_{\ref{thm_excludeRainbowClique}}(q,t)$ such that $G_{x} - A$ has a $\Sigma$-rendition of breadth at most $\breadth_{\ref{thm_excludeRainbowClique}}(q,t)$ and width at most $\vortex_{\ref{thm_excludeRainbowClique}}(q,t).$
 \end{itemize}
\end{enumerate}
Moreover, it holds that $\adhesion_{\ref{thm_excludeRainbowClique}}(q,t),~\apex_{\ref{thm_excludeRainbowClique}}(q,t),~ \vortex_{\ref{thm_excludeRainbowClique}}(q,t) \in \Ocal\big((qt)^{2300}\big)$ and $\genus_{\ref{thm_excludeRainbowClique}}(q,t),~\breadth_{\ref{thm_excludeRainbowClique}}(q,t)\in\Ocal((qt)^2).$ 

There also exists an algorithm that, given $t$ and $(G,\chi)$ as input, finds either a colorful minor model of a rainbow $K_{t}$ in $(G,\chi)$ or a tree-decomposition $(T,\beta)$ as above in time $2^{\poly(qt)} \cdot |G|^{3}|\!|G|\!|\cdot \log|G|$
\end{theorem}

Our proof of \zcref{thm_excludeRainbowClique} proves a slightly stronger variant that allows for a straightforward induction. 

Let $G$ be a graph.
A tuple $(T,r,\beta)$ is a \emph{rooted} tree-decomposition for $G$ if $(T,\beta)$ is a tree-decomposition for $G$ and $r\in V(T).$
We also define notions like \emph{width} and \emph{adhesion} for a rooted tree-decomposition $(T,r,\beta)$ by forwarding them to the tree-decomposition $(T,\beta).$

\begin{theorem}\label{thm_RainbowClique_induction}
There exist functions $\link_{\ref{thm_RainbowClique_induction}},\apex_{\ref{thm_RainbowClique_induction}},\breadth_{\ref{thm_RainbowClique_induction}},\vortex_{\ref{thm_RainbowClique_induction}}\colon\mathbb{N}^2\to\mathbb{N}$ such that for every $t\geq 1,$ every $q\geq 1,$ every $q$-colorful graph $(G,\chi),$ and every vertex set $X\subseteq V(G)$ with $|X|\leq 3\link_{\ref{thm_RainbowClique_induction}}(q, t) + 1$ one of the following holds:
\begin{enumerate}
 \item $(G,\chi)$ contains a rainbow $K_t$ as a colorful minor, or
 \item $(G, \chi)$ has a rooted tree-decomposition $(T, r, \beta)$ of adhesion at most $3\link_{\ref{thm_RainbowClique_induction}}(q,t)+\apex_{\ref{thm_RainbowClique_induction}}(q,t)+\vortex_{\ref{thm_RainbowClique_induction}}(q,t)+3$ and a (possibly empty) subset $L \subseteq V(T)$ of the leaves of $T$ such that for every $x \in V(T)$ one the following holds:
 \begin{itemize}
 \item $x \in L$ is a leaf with parent $d\in V(T)$ and $(G[\beta(x) \setminus \beta(d)], \chi)$ is restricted, or
 \item if $(G_{x}, \chi_{x})$ is the colorful torso of $G$ at $x,$ then there exists a set $A \subseteq V(G_{x})$ of size at most $3\link_{\ref{thm_RainbowClique_induction}}(q, t) + \apex_{\ref{thm_RainbowClique_induction}}(q,t) + 1$ and a surface of Euler-genus at most $((2q + 1)t)^{2}$ such that $G_{x} - A$ has a $\Sigma$-rendition of breadth at most $\breadth_{\ref{thm_RainbowClique_induction}}(q,t)$ and width at most $2\link_{\ref{thm_RainbowClique_induction}}(q,t) + \vortex_{\ref{thm_RainbowClique_induction}}(q,t) + 1.$ In addition, if $x = r$ then $X \subseteq A.$
 \end{itemize}
\end{enumerate}
Moreover, it holds that $\apex_{\ref{thm_RainbowClique_induction}}(q,t),~\breadth_{\ref{thm_RainbowClique_induction}}(q,t),~ \vortex_{\ref{thm_RainbowClique_induction}}(q,t) \in \Ocal\big((qt)^{115}\big),~\link_{\ref{thm_RainbowClique_induction}}(q,t)\in\Ocal\big( (qt)^{2300}\big).$

There also exists an algorithm that, given $t,$ and $(G,\chi)$ as input, finds either a colorful minor model of a rainbow $K_t$ in $(G,\chi)$ or a tree-decomposition $(T,\beta)$ as above in time $2^{\poly(qt)} \cdot |G|^3|\!|G|\!|\log|G|.$
\end{theorem}
\begin{proof}

We begin by declaring the functions involved.
Let $c_1$ be the constant from \zcref{thm_algogrid}.
\begin{align*}
 \apex_{\ref{thm_RainbowClique_induction}}(q,t) & \coloneqq \apex_{\ref{thm_localstructure}}((2q+1)t)\\
 \vortex_{\ref{thm_RainbowClique_induction}}(q,t) & \coloneqq 2\depth_{\ref{thm_localstructure}}((2q+1)t)\\
 \breadth_{\ref{thm_RainbowClique_induction}}(q,t) & \coloneqq \nicefrac{1}{2}((2q+1)t-3)((2q+1)t-4)\\
 \link_{\ref{thm_RainbowClique_induction}}(q,t) & \coloneqq c_1\wall_{\ref{thm_localstructure}}((2q+1)t, 3)^{20} + \apex_{\ref{thm_RainbowClique_induction}}(q,t) + \vortex_{\ref{thm_RainbowClique_induction}}(q,t) + 1
\end{align*}
Notice that the required bounds on all four functions follow directly from \zcref{thm_localstructure} and the above definitions.

We now proceed by induction on $|V(G)\setminus X|$ and start by discussing two fundamental cases.
In the following we assume that $(G,\chi)$ does not contain a rainbow $K_t$ as a colorful minor.

\smallskip
\textbf{Fundamental case 1:}
In case $|G|\leq 3\link_{\ref{thm_RainbowClique_induction}}(q, t) + 1,$ we may select $T$ to consist solely of the vertex $r$ and set $\beta(r)\coloneqq V(G)$ which satisfies the second outcome of our theorem.
Thus, we may assume
 $|G| > 3\link_{\ref{thm_RainbowClique_induction}}(q, t) + 1.$

\smallskip
\textbf{Fundamental case 2:}
Similarly, if $|X|<3\link_{\ref{thm_RainbowClique_induction}}(q, t) + 1,$ we may select any vertex $v\in V(G)$ and set $X'\coloneqq X\cup \{ v\},$ thereby achieving $|V(G)\setminus X'|<|V(G)\setminus X|$ and so we are done by induction.
Thus, we may assume $|X| = 3\link_{\ref{thm_RainbowClique_induction}}(q, t) + 1.$
\smallskip

We now call \zcref{prop_balancedseps} for $X$ and $k = \link_{\ref{thm_RainbowClique_induction}}(q, t)$ which has two possible outcomes.
Either, $X$ is $(k,\nicefrac{2}{3})$-linked, or there is a $\nicefrac{2}{3}$-balanced separator of size at most $k$ for $X.$
We treat these two cases separately.

\smallskip
\textbf{Case 1:}
There exists a $\nicefrac{2}{3}$-balanced separator $S$ of size at most $k$ for $X.$

In this case, for each component $H$ of $G-S$ we consider the graph $H' \coloneqq G[V(H) \cup S]$ together with the set $X_H \coloneqq S\cup (V(H)\cap X).$
Notice that, since $S$ is a $\nicefrac{2}{3}$-balanced separator for $X$ of size at most $k,$ we have that
\begin{align*}
 |X_H| \leq k + \left\lfloor\frac{2}{3}(3k+1)\right\rfloor \leq 3k < 3k+1.
\end{align*}
Hence, $H'$ and $X_H$ satisfy the properties of one of the two fundamental cases.
In either case, we obtain a rooted tree-decomposition $(T_H,r_H,\beta_H)$ together with a set of leaves $L_H \subseteq V(T_H)$ for $H'$ with $X_H \subseteq \beta_H(r_H).$
From here we define $(T,r,\beta)$ and $L$ as follows.
The tree $T$ is the tree obtained from the disjoint union of all $T_H$ by introducing a new vertex $r$ adjacent to all $r_H.$
The bags of $T$ are set to be $\beta(t) \coloneqq \beta_H(t)$ in case $t\in V(T_H)$ for some component $H$ of $G - S$ and otherwise we have $t=r$ and set $\beta(r)\coloneqq X\cup S.$
Moreover, we set $L$ to be the union of all $L_H.$
Notice that $(T,r,\beta)$ is indeed a rooted tree-decomposition as required and, in particular, $|\beta(r)|\leq 4k+1$ and so $G_r$ satisfies the second bullet point of the second outcome of our claim.

\smallskip
\textbf{Case 2:}
The set $X$ is $(k,\nicefrac{2}{3})$-linked.

Since $k \geq c_1\big(\wall_{\ref{thm_localstructure}}((2q+1)t)\big)^{20},$ we may now apply \zcref{thm_algogrid} to $X$ and obtain a $\wall_{\ref{thm_localstructure}}((2q+1)t)$-wall $W_0$ such that $\Tcal_{W_0} \subseteq \Tcal_X.$
This allows us to apply \zcref{thm_localstructure} to $W_0,$ asking for a clique minor of size $(2q+1)t.$
There are two possible outcomes:
\begin{enumerate}
 \item $G$ contains a $K_{(2q+1)t}$-minor controlled by $\mathcal{T}_{W_0},$ or
 \item there exists a set $A\subseteq V(G)$ of size at most $\apex_{\ref{thm_localstructure}}((2q+1)t)$ and a subwall $W_1\subseteq W_0-A$ such that $G - A$ has a $\Sigma$-rendition $\rho$ of breadth at most $\nicefrac{1}{2}((2q+1)t-3)((2q+1)t-4)$ and depth at most $\depth_{\ref{thm_localstructure}}((2q+1)t)$ such that $\Sigma$ is a surface where $K_{(2q+1)t}$ does not embed, $W_1$ is flat in $\rho,$ and $\rho$ is $(W_0-A)$-central.
\end{enumerate}
We further distinguish two subcases based on these outcomes.

\smallskip
\textbf{Case 2.1:}
The graph $G$ contains a model $\Xcal$ of $K_{(2q+1)t}$ controlled by $\mathcal{T}_{W_0}.$

In this case, we apply \zcref{thm_rainbowclique} to $\Xcal$ in $(G,\chi).$
Since, by our assumption, $(G,\chi)$ does not contain a rainbow $K_t$ as a colorful minor, the second outcome of \zcref{thm_rainbowclique} must hold.
That is, we obtain a set $S\subseteq V(G)$ of size at most $qt-1 < k$ such that the $\Tcal_{\Xcal}$-big component $H_0$ of $G-S$ is restricted.

Now let $H$ be any other component of $G-S$ apart from $H_0.$
Let $X_H\coloneqq S\cup (V(H)\cap X).$
Since $H\neq H_0,$ $H_0$ is the $\Tcal_{\Xcal}$-big component of $G-S,$ and $\Tcal_{\Xcal}\subseteq \Tcal_{W_0}$ as $\Xcal$ is controlled by $\mathcal{T}_{W_0},$ we obtain that
\begin{align*}
 |X_H| < k + \left\lfloor \frac{2}{3}|X| \right\rfloor < 3k.
\end{align*}
Hence, for each component $H$ of $G-S$ other than $H_0,$ we have that $H$ together with $X_H$ falls into one of the two fundamental cases.
Therefore, for each such $H$ we obtain a rooted tree-decomposition $(T_H,r_H,\beta_H)$ and a set $L_H$ of leaves of $T_H$ with $X_H\subseteq\beta_H(r_H)$ that satisfies the second outcome of our claim.

We define the desired rooted tree-decomposition $(T,r,\beta)$ with the set $L$ of leaves as follows.
Let $T$ be the tree obtained from the disjoint union of the $T_H$ for all components $H\neq H_0$ of $G-S$ by introducing a new isolated vertex $d$ and then a new vertex $r$ adjacent to $d$ and all $r_H.$
We set $\beta(d) \coloneqq V(H_{0}) \cup X \cup S,$ $\beta(r)\coloneqq X\cup S,$ and $\beta(t)\coloneqq \beta_H(t)$ for all $t\in V(T)$ such that there is a component $H\neq H_0$ of $G-S$ with $t\in V(T_H).$
Finally, we set $L$ to be the union of all $L_H$ together with the set $\{ d\}.$
As before, it is straightforward that $(T,r,\beta)$ together with $L$ is indeed the required rooted tree-decomposition as demanded by the second outcome of our claim.
This settles the first subcase.

\smallskip
\textbf{Case 2.2:}
There exists a set $A\subseteq V(G)$ of size at most $\apex_{\ref{thm_localstructure}}((2q+1)t)$ and a subwall $W_1\subseteq W_0-A$ such that $G - A$ has a $\Sigma$-rendition $\rho$ of breadth at most $\nicefrac{1}{2}((2q+1)t-3)((2q+1)t-4)$ and depth at most $\depth_{\ref{thm_localstructure}}((2q+1)t)$ such that $\Sigma$ is a surface where $K_{(2q+1)t}$ does not embed, $W_1$ is flat in $\rho,$ and $\rho$ is $(W_0-A)$-central.

Since $\rho$ is $(W_0-A)$-central, no cell of $\rho$ contains the majority side of any separation in $\Tcal_{W_0}-A.$
Since $\Tcal_{W_0}\subseteq \Tcal_X$ it follows that no cell of $\rho$ contains the big side of a separation from $\Tcal_X-A.$
Therefore, we have $|V(\sigma(c))\cap X|\leq \lfloor \nicefrac{2}{3} |X| \rfloor$ for all non-vortex cells $c\in C(\rho).$
Moreover, for every vortex $c\in C(\rho)$ and every bag $Y_c^i$ from the linear decomposition of adhesion at most $\depth_{\ref{thm_localstructure}}((2q+1)t)$ -- which is guaranteed to exist due to our assumptions on $\rho$ -- we have that also $|Y_c^i\cap X|\leq \lfloor \nicefrac{2}{3} |X| \rfloor.$

For every non-vortex cell $c\in C(\rho)$ let $H_c\coloneqq G[A\cup V(\sigma(c))]$ and $X_c\coloneqq A \cup \widetilde{c} \cup (X\cap V(\sigma(c))).$
Similarly, for every vortex cell $c\in C(\rho)$ and every $i\in[\ell_c],$ where $\ell_c$ is the number of bags in the linear decomposition of $c,$ let $H^i_c \coloneqq G[Y^i_c \cup A]$ and $X^i_c$ be the union of $A,$ $X\cap V(\sigma(c)),$ the -- at most two -- adhesion sets of $Y^i_c,$ and the vertex of $Y^{i}_{c}$ in $\tilde{c}.$
It follows that
\begin{align*}
 |X_c| & \leq \apex_{\ref{thm_localstructure}}((2q+1)t) + \left\lfloor \frac{2}{3}|X|\right\rfloor + 3 < 3k\text{, and}\\
 |X_v^i| & \leq 2\depth_{\ref{thm_localstructure}}((2q+1)t) + \apex_{\ref{thm_localstructure}}((2q+1)t) + 1 + \left\lfloor \frac{2}{3}|X| \right\rfloor < 3k
\end{align*}
for all non-vortex cells $c\in C(\rho),$ all vortex cells $v\in C(\rho),$ and all $i\in[\ell_v].$

Now for each $H_c$ and each $H^i_v,$ as before, we may observe that, together with their respective $X_c$ and $X^i_v,$ they fall into one of the two fundamental cases.
From here, we obtain for each of them, a rooted tree-decomposition $(T_c,r_c,\beta_c)$ and a set $L_c$ ($(T_v^i,r_v^i,\beta_v^i)$ and $L_v^i$ for the vortices $v\in C(\rho)$ and $i\in[\ell_v]$ respectively) as desired.
Finally, all that is left to do is to combine them into the rooted tree-decomposition $(T,r,\beta)$ with set $L$ as follows.

Let $T$ be the tree obtained from the disjoint union of all $T_c$ and $T_v^i$ by adding a new vertex $r$ adjacent to the $r_c$ and $r_c^i.$
We set $\beta(r)\coloneqq N(\rho)\cup A\cup X\cup B$ where $B$ denotes the union of all the adhesion sets of the linear decompositions of all of the vortices $v\in C(\rho).$
From our inequalities above it follows that 
\begin{align*}
|\beta(r)\cap Y^i_v| & \leq 2\depth_{\ref{thm_localstructure}}((2q+1)t) + 1 + \left\lfloor \frac{2}{3}|X| \right\rfloor\\
& \leq 2\link_{\ref{thm_RainbowClique_induction}}(q,t) + 2\depth_{\ref{thm_localstructure}}((2q+1)t) + 1
\end{align*}
as anticipated.
For all $t\in V(T)$ such that there is either a non-vortex cell $c\in C(\rho)$ with $t\in V(T_c)$ or a vortex cell $v\in C(\rho)$ and some $i\in[\ell_i]$ such that $t\in V(T^i_v)$ we set $\beta(t)\coloneqq \beta_c(t)$ and $\beta(t)\coloneqq \beta_v^i(t)$ respectively.
Finally, the set $L$ is the union of all sets $L_c$ and $L_v^i.$
With this, our proof is complete.
\end{proof}

\section{Excluding a rainbow grid}\label{sec_RainbowGrid}
We now move on from rainbow cliques to rainbow grids:
In this section, we explain how \zcref{thm_ExcludingRainbowGridIntro} can be deduced from a powerful variant of \zcref{thm_localstructure} by Paul, Protopapas, Thilikos and Wiederrecht \cite{PaulPTS2025LocalIndex,PaulPTW2024Obstructions}.

\subsection{Colors in a surface}\label{subsec_epall1}
The main purpose of this subsection is to import a set of tools from the work of Paul, Protopapas, Thilikos and Wiederrecht in \cite{PaulPTS2025LocalIndex} (see also \cite{PaulPTW2024Obstructions}).
These tools were developed to homogenize $\Sigma$-renditions with respect to the rooted minors one can find within their cells.
In \cite{PaulPTS2025LocalIndex} rooted minors are represented using colors, a feature that comes in handy for our purposes.

 \paragraph{Colorful Boundaried Indices.}

Consider a fixed non-negative integer $q.$
A \emph{boundaried $q$-colorful graph} is a quadruple $(G,\chi,B,\rho),$ where $(G,\chi)$ is a $q$-colorful graph, $B \subseteq V(G),$ and $\rho \colon B \to [|B|]$ is a bijection.
In other words $(G,\chi,B,\rho)$ is an alternative compact notation for a rooted $q$-colorful graph where all roots are pairwise distinct, i.e., for the rooted $q$-colorful graph $(G,\chi,\rho^{-1}(1),\ldots,\rho^{-1}(|B|)).$

Two boundaried $q$-colorful graphs $(G_1,\chi_1,B_1,\rho_1)$ and $(G_2,\chi_2,B_2,\rho_2)$ are \emph{compatible} if $\rho_{2}^{-1}\circ\rho_{1}$ is an isomorphism from $G_1[B_{1}]$ to $G_2[B_{2}]$ (where $B_{2}=\rho_{2}^{-1}\circ\rho_{1}(B_1)$), and for each $i\in[|B|],$ $\chi_{1}(\rho_1(i))=\chi_{2}(\rho_2(i)).$
Moreover, they are \emph{isomorphic} if there exists an isomorphism $\sigma \colon V(G_{1}) \to V(G_{2})$ between $G_{1}$ and $G_{2}$ such that $\rho_{2}^{-1}\circ\rho_{1} \subseteq \sigma$ and for every $u \in V(G) \setminus B_{1},$ $\chi_{1}(u) = \chi_{2}(\sigma(u)).$ 
We denote by $\Bcal^{q}$ the set of all (pairwise non-isomorphic) boundaried $q$-colorful graphs whose boundary contains at least $1$ and at most $3$ vertices.
A \emph{$q$-colorful boundaried index} is a function $\iota \colon \Bcal^{q} \rightarrow \Nbbb_{\geq 1}$ such that, if for two $q$-colorful boundaried graphs $\mathbf{B}_{1} = (G_{1}, \chi_{1},B_{1}, \rho_{1}) \in \Bcal$ and $\mathbf{B}_{2} = (G_{2}, \chi_{2},B_{2}, \rho_{2}) \in \Bcal^{q},$ $\iota(\mathbf{B}_{1}) = \iota(\mathbf{B}_{2}),$ then $\mathbf{B}_{1}$ and $\mathbf{B}_{2}$ are compatible.
Moreover, we demand that if $|B_{1}| < |B_{2}|,$ then $\iota(\mathbf{B}_{1}) < \iota(\mathbf{B}_{2}).$
We define the capacity of $\iota$ as $\mathsf{cap}(\iota) = \max_{\mathbf{B} \in \Bcal} \iota(\mathbf{B}).$
We moreover demand that $\mathsf{cap}(\iota) = [\ell],$ for some positive integer $\ell.$ 

\paragraph{Index representation.}

For the purposes of this section, and this section only, we assume that every $q$-colorful graph $(G, \chi)$ comes equipped with a \emph{universal ranking}, i.e. an injection $\zeta_{G} \colon V(G) \to |V(G)|.$
This is useful as it allows us to uniquely define boundaried graphs on its subgraphs.
For instance, if $(H, \psi)$ is a subgraph of $(G, \chi)$ and $\emptyset \neq B \subseteq V(H),$ then we may define the boundaried graph $(H, \psi, B, \rho),$ where for every $u \in B,$ $\rho(u) = i$ if and only if $\zeta_{G}(u)$ is the $i$-th smallest integer in $\zeta_{G}(B).$
This allows us to utilize the simpler notation $(H, \psi, B)$ instead of $(H, \psi, B, \rho).$

Next, we would like to express what it means for a collection of indices to be bidimensionally represented within a $\Sigma$-rendition.
This is exactly what will allow us to obtain the rainbow grid later on.

Let $\iota$ be a $q$-colorful boundary index, $(G, \chi)$ be a $q$-colorful graph, and $\rho$ be a $\Sigma$-rendition of $G.$
Given $t \in \mathbb{N}_{\geq 5},$ we say that a minor model $\{ G_{u} \}_{u \in V(\Gamma)}$ of a $(t \times t)$-grid $\Gamma$ in $G$ \emph{represents} $\iota$ in $\delta$ if
\begin{itemize}
\item $\bigcup_{u \in V(\Gamma)} G_{u}$ is disjoint from $V(\sigma(c)) \setminus N(\rho)$ for all vortices $c \in C(\rho)$ and
\item for every vertex $u \in V(\Gamma),$ for every $\alpha \in \mathsf{cap}(i),$ if there exists a non-vortex cell $c \in C(\rho)$ with $\iota(c) = \alpha,$ then there is a cell $c' \in C(\rho)$ with $\iota_{\rho}(c') = \alpha$ and $\tilde{c} \subseteq V(G_{u}).$
Here $\iota_{\rho}(c)$ denotes the value $\iota(\sigma(c), \chi, \tilde{c}).$
\end{itemize}

 \paragraph{An indexed local structure theorem.}
We are now ready to state a key result from \cite{PaulPTS2025LocalIndex}.
A preliminary version was already used in \cite{PaulPTW2024Obstructions}, but \cite{PaulPTS2025LocalIndex} contains the following, much more streamlined version.

\begin{proposition}\label{prop_bidim_structure}
There exists functions $\wall_{\ref{prop_bidim_structure}}\colon\mathbb{N}^4\to\mathbb{N},$ and $\apex_{\ref{prop_bidim_structure}}, \breadth_{\ref{prop_bidim_structure}}, \depth_{\ref{prop_bidim_structure}} \colon \Nbbb^{3} \to \Nbbb$ such that, for every non-negative integer $q,$ for every $q$-colorful boundary index $\iota$ with $\mathsf{cap}(\iota) = p \in \Nbbb_{\geq 1},$ every $r \in \mathbb{N},$ every $k, t \in \Nbbb_{\geq 5},$ every graph $H$ on $t$ vertices, every $q$-colorful graph $(G, \chi),$ and every $\wall_{\ref{prop_bidim_structure}}(r, t, k, p)$-wall $W \subseteq G,$ one of the following holds:
\begin{enumerate}
\item $G$ contains $H$ as a minor controlled by $\mathcal{T}_{W},$ or
\item there is a set $A \subseteq V(G),$ a $\Sigma$-rendition $\rho$ of $G - A,$ and a $(r + k)$-wall $W'$ such that
\begin{itemize}
\item $|A| \leq \apex_{\ref{prop_bidim_structure}}(t, k, p),$
\item $\Sigma$ is a surface where $H$ does not embed,
\item $\rho$ has breadth at most $\breadth_{\ref{prop_bidim_structure}}(t, k, p)$ and depth at most $\depth_{\ref{prop_bidim_structure}}(t, k, p),$
\item $W'$ is grounded in $\rho$ and $\Tcal_{W'} \subseteq \Tcal_W,$ and
\item there is a minor model of a $(k \times k)$-grid in $G - A$ that represents $\iota$ in $\rho.$
\end{itemize}
\end{enumerate}
Moreover,
\begin{align*}
\wall_{\ref{prop_bidim_structure}}(r, t, k, p) \ &\in \ r^{2^{\Ocal(p)}} \cdot t^{2^{\Ocal(p)}} + 2^{\Ocal((t^{2} + k^{2}) \cdot \log k \cdot {2^{\Ocal(p)})}}\\
\apex_{\ref{prop_bidim_structure}}(t, k, p), \depth_{\ref{prop_bidim_structure}}(t, k, p) \ &\in \ 2^{\Ocal((t^{2} + k^{2}) \cdot \log k \cdot {2^{\Ocal(p)})}}\text{, and}\\
\breadth_{\ref{prop_bidim_structure}}(t, k, p) \ &\in \ \Ocal(t^{2} + k^{2} \cdot 2^{\Ocal(p)}).
\end{align*}
Also, given an algorithm $\mathbf{A}$ for $\iota,$ there exists an algorithm that, given $k,$ $t,$ a graph $H,$ a graph $G,$ and a wall $W$ as above as input, finds one of these outcomes in time $2^{2^{(t+k)^{\Ocal(1)}2^{\Ocal(p)}}} \cdot |G|^{2} |\!|G|\!|.$
\end{proposition}

We should stress that the original formulation of \zcref{prop_bidim_structure} is not stated in terms of $q$-colorful graphs and $q$-colorful boundary indices.
However, a variant like \zcref{prop_bidim_structure} is straightforward to obtain as in the original proof the boundary index is only used to assign a coloring to the non-vortex cells of a $\Sigma$-rendition of $G,$ which in this case would depend on $\chi.$
Besides that, $\chi$ would not play any other role in its proof and can safely be disregarded thereafter.

\subsection{Locally excluding a rainbow grid}\label{subsec_localRainbowGrid}
\zcref{prop_bidim_structure} allows us to immediately deduce a strong local structure theorem for $q$-colorful graphs excluding rainbow grids.

\smallskip
For our purposes, we make use of a highly simplified case of a colorful boundaried index.
We wish to express $q$-colorful graphs as $q$-colorful boundaried graphs in a way that preserves the information of colors the vertices on the boundary have access to.
We have to be slightly careful so as to not map two non-equivalent $q$-colorful boundaried graphs to the same index even if their boundaries have access to the same set of colors.

To achieve this, let $\Pi^{q} = \{ \mathcal{B}^{q}_{1}, \ldots, \mathcal{B}^{q}_{\ell} \}$ be a partition of $\Bcal^{q}$ to maximal sets of equivalent $q$-colorful boundaried graphs such that for every $i \leq j \in [\ell],$ if $(G_{i}, \chi_{i}, B_{i}, \rho_{i}) \in \mathcal{B}^{q}_{i}$ and $(G_{j}, \chi_{j}, B_{j}, \rho_{j}) \in \mathcal{B}^{q}_{j},$ $|B_{i}| \leq |B_{j}|.$
Notice that since we consider boundaries of size at most $3,$ $|\Pi^{q}|$ is bounded by a constant depending only on $q.$
Consider a bijection $\phi_{q} \colon 2^{[q]} \times \Pi^{q} \to [2^{q} \cdot |\Pi^{q}|]$ such that for every $i < j \in [\ell]$ and every $Q, Q' \in 2^{[q]},$ $\phi_{q}(Q, \mathcal{B}^{q}_{i}) < \phi_{q}(Q', \mathcal{B}^{q}_{j}).$

Now, for every non-negative integer $q,$ we define a $q$-colorful boundaried index $\iota^{q}$ so that for every $q$-colorful boundaried graph $(G, \chi, B, \rho) \in \mathcal{B}^{q}_{i}$ for some $i \in [\ell],$ $\iota^{q}(G, \chi, B, \rho) = \phi_{q}(\bigcup_{v \in B} \chi(G_{v}), \mathcal{B}^{q}_{i})$ where $G_{v}$ is the component of $G$ that contains $v \in B.$
In other words, $\iota^{q}$ expresses the set of all colors that can be represented by one of the boundary vertices.

\begin{lemma}\label{lemma_localRainbowGrid}
There exist functions $\wall_{\ref{lemma_localRainbowGrid}}\colon\mathbb{N}^3\to\mathbb{N},$ $\apex_{\ref{lemma_localRainbowGrid}}, \breadth_{\ref{lemma_localRainbowGrid}}, \depth_{\ref{lemma_localRainbowGrid}} \colon \Nbbb^{2} \to \Nbbb$ such that, 
for all non-negative integers $q$ and $r,$ every $k \in \Nbbb_{\geq 3},$ every $q$-colorful graph $(G,\chi),$ and every $\wall_{\ref{lemma_localRainbowGrid}}(r, q, k)$-wall $W \subseteq G,$ one of the following holds:
\begin{enumerate}
\item $(G,\chi)$ contains the rainbow $(q,k)$-grid as a colorful minor, 
\item there exists a set $A\subseteq V(G)$ of size at most $\apex_{\ref{lemma_localRainbowGrid}}(q, k)$ such that the $\Tcal_{W}$-big component of $(G-A,\chi)$ is restricted, or
\item there exists a set $A \subseteq V(G),$ a $\Sigma$-rendition $\rho$ of $G - A,$ and a $(r+k)$-wall $W'$ such that
\begin{itemize}
\item $|A| \leq \apex_{\ref{lemma_localRainbowGrid}}(q, k),$
\item $\Sigma$ is a surface of Euler-genus in $\Ocal(q^2k^4),$
\item $\rho$ has breadth at most $\breadth_{\ref{lemma_localRainbowGrid}}(q, k)$ and depth at most $\depth_{\ref{lemma_localRainbowGrid}}(q, k),$ 
\item $W'$ is grounded in $\rho$ and $\Tcal_{W'} \subseteq \Tcal_W,$ and
\item if $q\neq 0,$ there exists $\emptyset \neq I\subseteq [q]$ such that $I\cap \chi(\sigma(c))=\emptyset$ for all non-vortex cells $c\in C(\rho).$
\end{itemize}
\end{enumerate}
Moreover,
\begin{align*}
\wall_{\ref{lemma_localRainbowGrid}}(r, q, k) \ &\in \ r^{2^{2^{\Ocal(q)}}} \cdot k^{2^{2^{\Ocal(q)}}} + 2^{\Ocal(k^{4} \cdot \log k) \cdot {2^{2^{\Ocal(q)}}}}\\
\apex_{\ref{lemma_localRainbowGrid}}(q, k), \depth_{\ref{lemma_localRainbowGrid}}(q, k) \ &\in \ 2^{\Ocal(k^{4} \cdot \log k) \cdot {2^{2^{\Ocal(q)}}}}\text{, and}\\
\breadth_{\ref{lemma_localRainbowGrid}}(q, k) \ &\in \ \Ocal(k^{4} \cdot 2^{2^{\Ocal(q)}}).
\end{align*}
Also, there exists an algorithm that, given $k,$ a $q$-colorful graph $(G,\chi),$ and a wall $W$ as above as input, finds one of these outcomes in time $2^{2^{k^{\Ocal(1)}2^{2^{\Ocal(q)}}}} \cdot |G|^{2} |\!|G|\!|.$ \end{lemma}

\begin{proof}
First, we set up our functions as follows:
\begin{align*}
\wall_{\ref{lemma_localRainbowGrid}}(r, q, k) \coloneqq~ & \wall_{\ref{prop_bidim_structure}}(r, (2q + 1)k^{2}, k, 2^{c_{q} q})\\
\apex_{\ref{lemma_localRainbowGrid}}(q, k) \coloneqq~ & \apex_{\ref{prop_bidim_structure}}((2q+1)k^2, k, 2^{c_{q} q})+qt-1\\ 
\breadth_{\ref{lemma_localRainbowGrid}}(q, k) \coloneqq~ & \breadth_{\ref{prop_bidim_structure}}((2q+1)k^2, k, 2^{c_{q} q})\\
\depth_{\ref{lemma_localRainbowGrid}}(q, k) \coloneqq~ & \depth_{\ref{prop_bidim_structure}}((2q+1)k^2, k, 2^{c_{q} q}).
\end{align*}
for a non-negative integer $c_{q}$ which depends only on $q$ and such that $2^{c_{q} q} \geq 2^{q} \cdot |\Pi^{q}|.$

We proceed by applying \zcref{prop_bidim_structure} with $H=K_{(2q+1)k^2},$ $t=(2q+1)k^2,$ and $\iota = \iota^{q}.$
Notice that this means that $\mathsf{cap}(\iota_{\chi}) = 2^{c_{q} q}.$
This yields one of two outcomes as follows.

 \paragraph{$G$ contains a $K_{(2q+1)k^2}$-minor controlled by $\mathcal{T}_{W}.$}
Let $\Xcal$ be the minor model of $K_{(2q+1)k^2}$ and notice that $\Tcal_{\Xcal}$ is a truncation of $W.$
We now apply \zcref{thm_rainbowclique}.
If this results in a colorful minor model of a rainbow $K_{k^2},$ we have found every $q$-colorful rainbow graph on $k^2$ vertices as a colorful minor in $(G,\chi)$ and thus, in particular, we have found the rainbow $(q,k)$-grid.
Hence, we may assume that we enter the second outcome of \zcref{thm_rainbowclique} which yields a set $A\subseteq V(G)$ of size at most
\begin{align*}
qt-1 \leq \apex_{\ref{lemma_localRainbowGrid}}(q, k)
\end{align*}
such that the $\Tcal_{\Xcal}$-big component of $(G-S,\chi)$ is restricted.
Since $\Tcal_{\Xcal}\subseteq \Tcal_W$ this means we have reached the second outcome of our assertion and thus we are done.

 \paragraph{$G$ does not contain a $K_{(2q+1)k^2}$-minor controlled by $\mathcal{T}_{W}.$}
In this case we now find a set $A\subseteq V(G)$ of size at most $\apex_{\ref{lemma_localRainbowGrid}}(q, k),$ a $(r+k)$-wall $W'$ whose tangle is a truncation of the tangle of $W,$ a surface $\Sigma$ where $K_{(2q+1)k^2}$ does not embed, and a $\Sigma$-rendition of breadth at most $\breadth_{\ref{lemma_localRainbowGrid}}(q, k)$ and depth at most $\depth_{\ref{lemma_localRainbowGrid}}(q, k)$ such that $W'$ is grounded in $\rho$ and $\rho$ bidimensionally $k$-represents $\iota_{\chi}.$

Notice that the Euler-genus of $\Sigma$ is in $\Ocal(q^2k^4)$ as desired.

It follows that there exists a minor model $\Xcal$ of the $(k\times k)$-grid in $G$ which is disjoint from the interiors of the vortices of $\rho$ such that for each $X\in\Xcal$ and every $\alpha\in[2^q]$ for which there is a non-vortex $c\in C(\rho)$ with $\iota^{q}_{\rho}(c)=\alpha$ there exists a non-vortex cell $c'\in C(\rho)$ with $\iota^{q}_{\rho}(c')=\alpha$ such that $\widetilde{c}\subseteq X.$

Notice that the boundary of any non-vortex cell of $\rho$ acts as a separator of size at most $3$ in $G-A.$
Hence, for any $X\in\Xcal$ and each non-vortex cell $c\in C(\rho)$ with $\widetilde{c}\subseteq X$ it follows that each component $J$ of $\sigma(c)$ that contains a vertex from $\widetilde{c}$ must be vertex-disjoint from all $Y\in\Xcal\setminus \{ X\}.$
Hence, for each such $J$ we may assume $V(J)\subseteq X$ and, moreover, since $X$ is connected, this means that $\chi(J)\subseteq \chi(X).$

Suppose for every $i\in[q]$ there exists $\alpha\in[2^q]$ and $\Bcal^{q} \in \Pi^{q}$ such that $(i, \Bcal^{q}) \in \phi^{-1}_{q}(\alpha)$ and $\iota^{q}_{\rho}(c) = \alpha$ for some $c \in C(\rho).$
Let us denote by $\Vcal \subseteq C(\rho)$ the set of all vortex cells of $\rho.$
Then we have that
\begin{align*}
 [q] = \{ i \in [q] \mid c\in C(\rho)\setminus\Vcal \text{ and } (i, \Bcal^{q}) \in \phi^{-1}_{q}(\iota^{q}_{\rho}(c))\}.
\end{align*}
Since $\Xcal$ represents $\iota^{q}$ in $\rho,$ this implies that $\chi(X)=[q]$ for all $X\in\Xcal$ which means that $\Xcal$ is indeed a colorful minor model of the rainbow $(q,k)$-grid.

Let
\begin{align*}
I \coloneqq \big\{ i\in[q] \mid \text{there is no }c\in C(\rho)\setminus\Vcal\text{ with } (i, \Bcal^{q}) \in \phi^{-1}_{q}(\iota^{q}_{\rho}(c)) \big\}.
\end{align*}
By our discussion it follows that, if $q\neq 0$ and $I=\emptyset,$ $(G,\chi)$ contains the rainbow $(q,k)$-grid as a colorful minor which is the first outcome of our assertion.
Hence, we may assume $I\neq \emptyset.$
Moreover, for every $c\in C(\rho)\setminus\Vcal(\rho)$ we may assume that every component of $\sigma(c)$ contains a vertex from $\widetilde{c}$ as we may remove any component of $\sigma(c)$ that does not contain a vertex of $\widetilde{c}$ and add it to $\sigma(v)$ for some $v\in \Vcal$ without altering the essential properties of $\rho.$
In the case where $\Vcal=\emptyset$ we may artificially introduce a new cell that acts as the storage for all such components.

Therefore, if $q\neq 0,$ we have that $I\neq\emptyset$ and $I\cap \chi(\sigma(c))=\emptyset$ for all non-vortex cells $c\in C(\rho)$ as desired.
This completes the proof.
\end{proof}

\subsection{The proof of \zcref{thm_ExcludingRainbowGridIntro}}\label{subsec_GlobalRainbowGrid}
We are now ready for the proof of \zcref{thm_ExcludingRainbowGridIntro}.
The proof will largely follow a similar strategy as the one for \zcref{thm_intro_ExcludeRainbowClique}.
Indeed, we will prove a stronger version similar to \zcref{thm_excludeRainbowClique} and its strengthened form \zcref{thm_RainbowClique_induction}.
\zcref{thm_ExcludingRainbowGridIntro} will then follow as a corollary.
This slightly stronger version reads as follows:

\begin{theorem}\label{thm_ExcludingRainbowGrid}
There exist functions $\link_{\ref{thm_ExcludingRainbowGrid}},$ $\apex_{\ref{thm_ExcludingRainbowGrid}},$ $\breadth_{\ref{thm_ExcludingRainbowGrid}},$ $\vortex_{\ref{thm_ExcludingRainbowGrid}}\colon\mathbb{N}^2\to\mathbb{N}$ such that for all $q,k\in\mathbb{N}_{\geq 1},$ and all $q$-colorful graphs $(G,\chi)$ one of the following holds:
\begin{enumerate}
 \item $(G,\chi)$ contains the rainbow $(q,k)$-grid as a colorful minor, or
 \item $(G,\chi)$ has a tree-decomposition $(T,\beta)$ of adhesion at most $\link_{\ref{thm_ExcludingRainbowGrid}}(q,k)$ such that for all $t\in V(T)$ one the following holds:
 \begin{itemize}
 \item $t$ is a leaf with unique neighbor $d$ such that $(G[\beta(t) \setminus \beta(d)], \chi)$ is restricted, or
 \item if $(G_t,\chi_t)$ is the colorful torso of $(G,\chi)$ at $t,$ then there exists a set $A\subseteq V(G_t)$ of size at most $\apex_{\ref{thm_ExcludingRainbowGrid}}(q,k),$ a surface of Euler-genus at most $\Ocal(q^2k^4),$ and a non-empty set $I\subseteq[q],$ such that $G_t-A$ has a $\Sigma$-rendition of breadth at most $\breadth_{\ref{thm_ExcludingRainbowGrid}}(q,k)$ and width at most $\vortex_{\ref{thm_ExcludingRainbowGrid}}(q,k)$ such that $\widetilde{c} = V(\sigma(c))$ for all non-vortex cells $c\in C(\rho),$ and $I\cap \chi_t(x)\neq\emptyset$ if and only if $x\in V(\sigma(v))\setminus\widetilde{v}$ for some vortex $v\in C(\rho).$
 \end{itemize}
\end{enumerate}
Moreover,
\begin{align*}
\link_{\ref{thm_ExcludingRainbowGrid}}(q,k), \apex_{\ref{thm_ExcludingRainbowGrid}}(q,k), \vortex_{\ref{thm_ExcludingRainbowGrid}}(q,k) \ &\in \ 2^{\Ocal(k^{4} \cdot \log k) \cdot {2^{2^{\Ocal(q)}}}}\text{, and}\\
\breadth_{\ref{thm_ExcludingRainbowGrid}}(q,k) \ &\in \ \Ocal(k^{4} \cdot 2^{2^{\Ocal(q)}}).
\end{align*}
Also, there exists an algorithm that, given $k,$ and a $q$-colorful graph $(G,\chi)$ as input, finds one of these outcomes in time $2^{2^{k^{\Ocal(1)}2^{2^{\Ocal(q)}}}} \cdot |G|^{3} |\!|G|\!| \log|G|.$ \end{theorem}

As announced, we strengthen \zcref{thm_ExcludingRainbowGrid} further to allow for a smooth induction.

\begin{theorem}\label{thm_ExcludingRainbowGrid_Induction}
There exist functions $\link_{\ref{thm_ExcludingRainbowGrid_Induction}},$ $\apex_{\ref{thm_ExcludingRainbowGrid_Induction}},$ $\breadth_{\ref{thm_ExcludingRainbowGrid_Induction}},$ $\vortex_{\ref{thm_ExcludingRainbowGrid_Induction}}\colon\mathbb{N}^2\to\mathbb{N}$ such that for all $q,k\in\mathbb{N}_{\geq 1},$ $k\geq 3,$ all $q$-colorful graphs $(G,\chi),$ and every set $X\subseteq V(G)$ with $|X|\leq 3\link_{\ref{thm_ExcludingRainbowGrid_Induction}}(q,k)+1$ one of the following holds:
\begin{enumerate}
 \item $(G,\chi)$ contains the rainbow $(q,k)$-grid as a colorful minor, or
 \item $(G,\chi)$ has a rooted tree-decomposition $(T,r,\beta)$ of adhesion at most $3\link_{\ref{thm_ExcludingRainbowGrid_Induction}}(q,k)+\apex_{\ref{thm_ExcludingRainbowGrid_Induction}}(q,k)+\vortex_{\ref{thm_ExcludingRainbowGrid_Induction}}(q, k) + 3$ and a (possibly empty) subset $L \subseteq V(T)$ of leaves of $T$ such that $X \subseteq \beta(r)$ and for all $t \in V(T)$ one the following holds:
 \begin{itemize}
 \item $t \in L$ is a leaf with a unique neighbor $d$ such that $(G[\beta(t) \setminus \beta(d)],\chi)$ is restricted, or
 \item if $(G_t, \chi_t)$ is the colorful torso of $(G, \chi)$ at $t,$ then there exists a set $A \subseteq V(G_t)$ of size at most $3\link_{\ref{thm_ExcludingRainbowGrid_Induction}}(q,k)+\apex_{\ref{thm_ExcludingRainbowGrid_Induction}}(q,k) + 1,$ a surface of Euler-genus at most $\Ocal(q^2k^4),$ and a non-empty set $I\subseteq[q],$ such that
 \begin{enumerate}
 \item $G_t-A$ has a $\Sigma$-rendition of breadth at most $\breadth_{\ref{thm_ExcludingRainbowGrid_Induction}}(q,k)$ and width at most $2\link_{\ref{thm_ExcludingRainbowGrid_Induction}}(q,k) + \vortex_{\ref{thm_ExcludingRainbowGrid_Induction}}(q,k) + 1$ such that $\widetilde{c} = V(\sigma(c))$ for all non-vortex cells $c\in C(\rho),$ and $I \cap \chi_t(x)\neq\emptyset$ if and only if $x \in V(\sigma(v)) \setminus \widetilde{v}$ for some vortex $v \in C(\rho),$ and
 \item if $d$ is a child of $t$ with $\beta(d)\cap\beta(t)\subseteq A\cup \widetilde{c}$ for some non-vortex cell $c\in C(\rho),$ then $d\in L.$
 \end{enumerate}
 \end{itemize}
\end{enumerate}
Moreover,
\begin{align*}
\link_{\ref{thm_ExcludingRainbowGrid_Induction}}(q,k), \apex_{\ref{thm_ExcludingRainbowGrid_Induction}}(q,k), \vortex_{\ref{thm_ExcludingRainbowGrid_Induction}}(q,k) \ &\in \ 2^{\Ocal(k^{4} \log k) \cdot {2^{2^{\Ocal(q)}}}}\text{, and}\\
\breadth_{\ref{thm_ExcludingRainbowGrid_Induction}}(q,k) \ &\in \ \Ocal(k^{4} \cdot 2^{2^{\Ocal(q)}}).
\end{align*}
Also, there exists an algorithm that, given $k,$ and a $q$-colorful graph $(G,\chi)$ as input, finds one of these outcomes in time $2^{2^{k^{\Ocal(1)}2^{2^{\Ocal(q)}}}} \cdot |G|^{3} |\!|G|\!| \log|G|.$ 
\end{theorem}

\begin{proof}
First, we define the functions involved as follows where $c_1$ is the constant from \zcref{thm_algogrid}.
\begin{align*}
 \apex_{\ref{thm_ExcludingRainbowGrid_Induction}}(q,k) & \coloneqq \apex_{\ref{lemma_localRainbowGrid}}(q, k)\\
 \vortex_{\ref{thm_ExcludingRainbowGrid_Induction}}(q,k) & \coloneqq 2\depth_{\ref{lemma_localRainbowGrid}}(q, k)\\
 \breadth_{\ref{thm_ExcludingRainbowGrid_Induction}}(q,k) & \coloneqq \breadth_{\ref{lemma_localRainbowGrid}}(q, k)\\
 \link_{\ref{thm_ExcludingRainbowGrid_Induction}}(q,k) & \coloneqq c_1\wall_{\ref{lemma_localRainbowGrid}}(\apex_{\ref{thm_ExcludingRainbowGrid_Induction}}(q,k) + \vortex_{\ref{thm_ExcludingRainbowGrid_Induction}}(q,k)+4,q, k)^{20}.
\end{align*}
Moreover, in the following let $k' \coloneqq \apex_{\ref{thm_ExcludingRainbowGrid_Induction}}(q,k)+ \vortex_{\ref{thm_ExcludingRainbowGrid_Induction}}(q,k) + k + 4.$

We now proceed by induction on $|V(G)\setminus X|$ and start by discussing two fundamental cases.
Exactly as in our proof of \zcref{thm_RainbowClique_induction}, we first discuss the two \textbf{fundamental cases} where either $G$ has only few vertices, or $|X|<3\link_{\ref{thm_ExcludingRainbowGrid_Induction}}(q,k)+1.$
These two cases are fully identical to the ones in \zcref{thm_RainbowClique_induction} and we therefore omit them here.
Then two further cases follow:

\textbf{Case 1}, where $X$ has a small balanced separator, and \textbf{Case 2}, where $X$ does not have a small balanced separator.
As before, for \textbf{Case 1} there is no difference to the proof of \zcref{thm_RainbowClique_induction} and thus, we may omit it here.

This leaves us with the last case.
\smallskip

\textbf{Case 2:}
The set $X$ is $(\link_{\ref{thm_ExcludingRainbowGrid_Induction}}(q,k),\nicefrac{2}{3})$-linked.
\smallskip

Since $\link_{\ref{thm_ExcludingRainbowGrid_Induction}}(q,k) = c_1\wall_{\ref{lemma_localRainbowGrid}}(k' - k, q, k)^{20},$ we may now apply \zcref{thm_algogrid} to $X$ and obtain a $\wall_{\ref{lemma_localRainbowGrid}}(k' - k, q, k)$-wall $W_0$ such that $\Tcal_{W_0} \subseteq \Tcal_X.$
This further allows for the application of \zcref{lemma_localRainbowGrid} to the wall $W_0$ in $(G,\chi).$
In the following we distinguish three cases based on the three possible outcomes of \zcref{lemma_localRainbowGrid}.
\smallskip

\textbf{Case 2.1:}
$(G,\chi)$ contains the rainbow $(q,k)$-grid as a colorful minor.
\smallskip

This is one of the possible outcomes of our assertion and thus, there is nothing more to show.
\medskip

\textbf{Case 2.2:}
There exists a set $A\subseteq V(G)$ of size at most $\apex_{\ref{lemma_localRainbowGrid}}(q, k)=\apex_{\ref{thm_ExcludingRainbowGrid_Induction}}(q,k)$ such that the $\Tcal_{W_0}$-big component of $(G-A,\chi)$ is restricted.
\smallskip

Let $H_0$ be the $\Tcal_{W_0}$-big component of $G-A$ and let $H_1,\dots,H_{\ell}$ denote all other components of $G-A.$
Observe that, since $\Tcal_{W_0}\subseteq\Tcal_X$ we have that $|V(H_0)\cap X|\geq 2\link_{\ref{thm_ExcludingRainbowGrid_Induction}}(q,k) + 1$
And thus, $|V(H_i)\cap X|\leq \link_{\ref{thm_ExcludingRainbowGrid_Induction}}(q,k)$ for all $i\in[\ell].$

We now describe how to create the rooted tree-decomposition $(T,r,\beta)$ as required by the assertion.

First, for each $i\in[\ell]$ select a vertex $v_i\in V(H_i)$ and let $X_i\coloneqq (X\cap V(H_i))\cup A\cup \{v_i\}.$
Notice that, by definition $\link_{\ref{thm_ExcludingRainbowGrid_Induction}}(q, k) > \apex_{\ref{lemma_localRainbowGrid}}(q, k) = \apex_{\ref{thm_ExcludingRainbowGrid_Induction}}(q,k)$ and therefore we have that $|X_i|\leq 3\link_{\ref{thm_ExcludingRainbowGrid_Induction}}(q,k)+1$ and, moreover, $|V(H_i)\setminus X_i| < |V(G)\setminus X|$ for all $i\in[\ell].$
Thus, we may apply the induction hypothesis to each $(H_i,\chi),$ thereby obtaining a rooted tree-decomposition $(T_i,r_i,\beta_i)$ with leaf sets $L_i$ as described in the assertion for each $i\in[\ell].$
In particular this means $X_i\subseteq \beta_i(r_i).$

Let now $T$ be the tree obtained from the disjoint union of the $T_i$ by introducing a vertex $d$ and a vertex $r$ where $r$ is adjacent to $d$ and all $r_i,$ $i\in[\ell].$
Moreover, we set $L\coloneqq \{ d\}\cup \bigcup_{i\in[\ell]}L_i.$
We set $\beta(r)\coloneqq X\cup A$ and $\beta(d) \coloneqq V(H_0)\cup X\cup A$ as well as $\beta(t)\coloneqq \beta_i(t)$ for each $t\in V(T_i)$ and every $i\in[\ell].$
One can now easily observe that $(T,r,\beta)$ is the desired rooted tree-decomposition of $(G,\chi).$
\medskip

\textbf{Case 2.3:}
There exists a set $A \subseteq V(G),$ a $\Sigma$-rendition $\rho$ of $G - A,$ and a $k'$-wall $W_1$ such that
\begin{itemize}
\item $\Tcal_{W_1}\subseteq \Tcal_{W_0}$
\item $|A| \leq \apex_{\ref{thm_ExcludingRainbowGrid_Induction}}(q,k),$
\item $\Sigma$ is a surface of Euler-genus in $\Ocal(q^2k^4),$
\item $\rho$ has breadth at most $\breadth_{\ref{thm_ExcludingRainbowGrid_Induction}}(q,k)$ and depth at most $\depth_{\ref{lemma_localRainbowGrid}}(q, k),$ 
\item $W_1$ is grounded in $\rho,$ and
\item there exists $\emptyset \neq I\subseteq [q]$ such that $I\cap \chi(\sigma(c))=\emptyset$ for all non-vortex cells $c\in C(\rho).$
\end{itemize}
\smallskip

Let us begin with a number of minor observations.
First of all, $W_1$ is grounded in $\rho.$
This means that for every non-vortex cell $c\in C(\rho)$ no entire row or column of $W_1$ can be contained in $\sigma(c).$
Since $\widetilde{c}$ has size at most $3$ and separates $\sigma(c)$ from the rest of $G-A$ in $G-A,$ and $\Tcal_{W_1}\subseteq \Tcal_{W_0}\subseteq \Tcal_X$ it follows that $|X\cap V(\sigma(c))|\leq \link_{\ref{thm_ExcludingRainbowGrid_Induction}}(q,k)+3$ where the additional $3$ vertices may only be contained in $\widetilde{c}.$

Similarly, let $v\in V(\rho)$ be any vortex with linear decomposition $\langle Y^v_i\rangle_{i\in[\ell_{v}]}$ of adhesion at most $\depth_{\ref{lemma_localRainbowGrid}}(q, k).$
It follows that there exists a set $S^v_i\subseteq Y^v_i$ of size at most $\vortex_{\ref{thm_ExcludingRainbowGrid_Induction}}(q,k)+1< k'$ that separates $Y^v_i\setminus S^v_i$ from the rest of $G-A.$
Moreover, $Y^v_i$ is disjoint from $W_1.$
Indeed, one can observe that no column or row of $W_1$ may be fully contained in $Y^v_i.$
Hence, we deduce that $|(Y^v_i\setminus S^v_i)\cap X|\leq \link_{\ref{thm_ExcludingRainbowGrid_Induction}}(q,k).$

Now, for every non-vortex cell $c\in C(\rho)$ let $H_c\coloneqq G[V(\sigma(c))\cup A],$ select $x_c\in V(H_c)\setminus (A\cup X)$ if it exists, and let $X_c\coloneqq A\cup (X\cap V(\sigma(c)))\cup \{x_c\}.$
Notice that $|X_c| \leq 3\link_{\ref{thm_ExcludingRainbowGrid_Induction}}(q,k)+1$ and $|V(G)\setminus X|>|V(H_c)\setminus X_c|$ or $V(H_c)= X_c.$
Moreover, from the fact that $I\neq \emptyset$ we may deduce that $\chi(V(\sigma(c))) \neq [q].$
Hence, each of the $(\sigma(c),\chi)$ is restricted.

Similarly, for every vortex $v\in C(\rho)$ and $i\in[\ell_v]$ let $H^v_i\coloneqq G[Y^v_i\cup A],$ select $x^v_i\in Y^v_i\setminus(A\cup X)$ if it exists and let $$X^i_v\coloneqq A\cup S^v_i\cup (X\cap Y^v_i)\cup \{x^v_i\}.$$
As before, one may observe that $|X^i_v|\leq 3\link_{\ref{thm_ExcludingRainbowGrid_Induction}}(q,k)+1$ and $|V(G)\setminus X|>|V(H^i_v)\setminus X^i_v|$ or $V(H^v_i)= X^v_i.$

These observations give us the right to call the induction hypothesis for every pair of a vortex $v\in C(\rho)$ and $i\in[\ell_v]$ on $(H^v_i,\chi)$ and $X^v_i.$
This yields, in the absence of a rainbow $(q,k)$-grid as a colorful minor, rooted tree-decompositions $(T^v_i,r^v_i,\beta^v_i)$ as described by the assertion for each such pair $v,i.$

We are now ready to construct the rooted tree-decomposition $(T,r,\beta)$ as follows.
Let $T$ be the tree obtained from the disjoint union of the $T^v_i$ for all $i\in[\ell_i]$ and all vortices $v\in C(\rho)$ by introducing, for each non-vortex cell $c\in C(\rho)$ a new vertex $d_c,$ then introducing a new vertex $r,$ and finally joining $r$ to all $d_c$ and all $r^v_i.$
We set 
\begin{align*}
\beta(r)\coloneqq X\cup A\cup N(\rho)\cup \bigcup_{\substack{v\in C(\rho)\\\text{vortex}}}\bigcup_{i\in[\ell_v]}S^v_i.
\end{align*}
For each non-vortex cell $c\in C(\rho)$ we set $\beta(c) = X_c\cup V(\sigma(c)),$ and for every vortex $v\in C(\rho)$ and each $i\in[\ell_v]$ we set $\beta(t)\coloneqq \beta^v_i(t)$ for all $t\in V(T^v_i).$
It is now easy to check that $(T,\beta)$ is indeed a tree-decomposition and from our discussion above it follows that it has exactly the properties required by our assertion.
\end{proof}

\section{Excluding segregated grids}\label{sec_SegregatedGrids}

In this section, we are concerned with the proof of \zcref{thm_restrictiveTreewidthIntro}, which is an analogue of the Grid Theorem of Robertson and Seymour for \textsl{torso treewidth}.
Similar to the organization of \zcref{sec_RainbowGrid}, we begin with some preliminary results and a \say{local} variant of the structure theorem before proceeding to the actual proof of \zcref{thm_restrictiveTreewidthIntro}.

\subsection{A homogeneous flat wall}\label{sec_HomogeneousFlatWall}

In this section we prove an auxiliary result on flat walls which has been proven in several variations before.
We include our own proof specifically for the setting of colorful graphs for the sake of completeness.

Let $(G,\chi)$ be a $q$-colorful graph with a flat $k$-wall $W\subseteq G$ where $k\geq 3$ and a $\Sigma$-rendition $\rho$ witnessing that $W$ is flat in $G$ and where $\Sigma$ is the sphere.
Let $B$ be a brick of $W.$
Then let $\Delta_B\subseteq\Sigma$ be the minimal disk such that $H_B\coloneqq \bigcup_{\substack{c\in C(\rho)\\c\subseteq \Delta_B}}\sigma(c)$ contains $B$ but not the perimeter of $W.$
Notice that by the minimality of $\Delta_B$ it follows that $B$ is the only brick of $W$ contained in $H_B.$

We say that $W$ is \emph{homogeneous} if for all bricks $B$ and $B'$ of $W,$ we have $\chi(V(H_B)) = \chi(V(H_{B'})).$

\begin{lemma}\label{lemma_HomogeneousFlatWall}
There exists a function $\hfw\colon\mathbb{N}^3\to\mathbb{N}$ such that for all integers $q\geq 0,$ $t\geq 1,$ and $r\geq 3$ and every $q$-colorful graph $(G,\chi),$ and every $\hfw(q,t,r)$-wall $W$ in $G$ one of the following holds:
\begin{enumerate}
 \item there exists a model of $K_t$ in $G$ which is controlled by $\mathcal{T}_{W},$ or
 \item there exists a set $Z$ with $|Z|<16t^3$ and an $r$-subwall $W'$ of $W$ which is disjoint from $Z,$ and flat and homogeneous in $G-Z.$
\end{enumerate}
Moreover, $\hfw(q,t,r)\in t^{\Ocal(1)}r^{2^{\Qcal(q)}}$ and there exists a $(t+r)^{2^{\Ocal(q)}} \cdot |\!|G|\!|$-time algorithm which finds either the model of $K_t$ or the set $Z,$ the subwall $W',$ and a $\Sigma$-rendition witnessing that $W'$ is flat in $G-Z.$
\end{lemma}

\begin{proof}
We begin by defining an auxiliary recursive function as follows.
\begin{align*}
 f(0,y) & \coloneqq y\\
 f(x,y) & \coloneqq f(x-1,y)^2 + 2f(x-1,y) +2\text{ for }x\in\mathbb{N}_{\geq 1}.
\end{align*}
From this recursion it follows that $f(x,y) \in y^{2^{\Ocal(x)}}.$
We now set
\begin{align*}
\hfw(q,t,r) \coloneqq 100t^3(f(q,r) + 2t + 2).
\end{align*}
This allows us to call upon \zcref{thm_flatwall} to find either a model of $K_t$ controlled by $\mathcal{T}_{W},$ in which case we would be done, or a set $Z\subseteq V(G)$ of size less than $16t^3$ and an $f(q,r)$-subwall $W_0\subseteq W$ that is flat in $G-Z.$

From here on, we will be further refining $W_0.$
For now let $J_0\subseteq [q]$ be the -- possibly empty -- collection of all colors $j$ such that $j\in\chi(V(H_B))$ for all bricks $B$ of $W_0$ and let $I_0\subseteq [q]$ the -- also possibly empty -- collection of all colors $i$ such that there does not exist a single brick $B$ of $W_0$ with $i\in\chi(V(H_B)).$

Our goal is to find sets $J_i\supseteq J_0$ and $I_i\supseteq I_0$ together with a subwall $W_i$ of $W_0$ of order at least $r$ such that $W_i$ is homogeneous, $\chi(V(H_B))=J_i$ for all bricks $B$ of $W_i,$ and $I_i=[q]\setminus J_i.$
This is done by induction over the number $q$ of colors and in each step we add precisely one new color to either the set $J_i$ or the set $I_i$ until all colors have been distributed.
Clearly, such a process must always stop after at most $q$ steps.
Moreover, the base case $q=0$ is trivially true.

Suppose for $i\in[0,q-1]$ we have already constructed the sets $J_i$ and $I_i$ as well as an $f(q-i,r)$-wall $W_i\subseteq W_0$ such that $J_i\subseteq \chi(V(H_B))$ and $I_i\cap \chi(V(H_B))=\emptyset$ for all bricks $B$ of $W_i,$ and $|J_i\cup I_i|\geq i.$
Then fix $j\in[q]\setminus (J_i\cup I_i)$ to be as small as possible.
In case no such $j$ exists we are already done, so we may assume that there is a possible choice.

Next, let us denote the vertical paths of $W_i$ by $P_1,\dots,P_{f(q-i,r)}$ and the horizontal paths of $W_i$ by $Q_1,\dots, Q_{f(q-i,r)}.$
Next, for each $x\in[f(q-i-1,r)+2]$ let $h_x\coloneqq x+(x-1)f(q-i-1,r).$
Let $W_a$ be the $(f(q-i-1,r)+2)$-wall consisting of the vertical paths $P_{h_x}$ and horizontal paths $Q_{h_x}$ for all $x\in[f(q-i-1,r)+2].$
Observe that for every brick $B$ of $W_a$ the graph $H_B$ contains an $(f(q-i-1,r)+2)$-subwall $W_B$ of $W_i$ whose perimeter is precisely the cycle $B.$

Now we are in one of two cases:
Either $j\in \chi(V(H_B))$ for all bricks $B$ of $W_a,$ or there exists a brick $B'$ of $W_a$ such that $j\notin\chi(V(H_{B'})).$

In the first case we let $W_{i+1}$ be the $f(q-i-1,r)$-subwall of $W_a$ that avoids its perimeter, and set $J_{i+1}\coloneqq J_i\cup \{ j\},$ and $I_{i+1}\coloneqq I_i.$
In this case we are done by induction.

In the second case, we let $W_{i+1}$ be the $f(q-i-1,r)$-subwall of $W_{B'}$ that avoids $B',$ and we set $J_{i+1}\coloneqq J_i$ and $I_{i+1}\coloneqq I_i\cup\{ j\}.$
Hence, also here we are done by induction and our proof is complete.
\end{proof}

\subsection{Locally excluding segregated grids}
We are now ready to provide a full local structure theorem for $q$-colorful graphs that exclude all $(q,k)$-segregated grids as colorful minors.

\begin{theorem}\label{thm_localSegregation}
There exist functions $\link_{\ref{thm_localSegregation}},\apex_{\ref{thm_localSegregation}}\colon\mathbb{N}^2\to\mathbb{N}$ such that for all integers $q\geq 1$ and $k,$ every $q$-colorful graph $(G,\chi),$ and every $(\link_{\ref{thm_localSegregation}}(q,k),\nicefrac{2}{3})$-linked set $X\subseteq V(G)$ one of the following holds:
\begin{enumerate}
 \item $G$ contains a $(q,k)$-segregated grid as a colorful minor, or
 \item there exists a set $A\subseteq V(G)$ with $|A|\leq \apex_{\ref{thm_localSegregation}}(q,k)$ such that the $\Tcal_X$-big component of $(G-A,\chi)$ is restricted. 
\end{enumerate}
Moreover, $\apex_{\ref{thm_localSegregation}}(q,k)\in \Ocal(q^9k^6),$ $\link_{\ref{thm_localSegregation}}(q,k)\in k^{2^{\Ocal(q)}},$ and there exists an algorithm that, given $q,$ $k,$ $(G,\chi),$ and $X$ as above as input finds one of these outcomes in time $2^{k^{2^{\Ocal(q)}}} \cdot |G|^2|\!|G|\!| \log|G|.$
\end{theorem}

\begin{proof}
As it is custom, we begin by setting up our numbers.
Let $c_1$ be the constant from \zcref{thm_algogrid} and $\hfw$ be the function from \zcref{lemma_HomogeneousFlatWall}.
\begin{align*}
 \apex_{\ref{thm_localSegregation}}(q,k) \coloneqq~ & 16(2q+1)^3(qk)^6 + (q+1)^2k \\
 \link_{\ref{thm_localSegregation}}(q,k) \coloneqq~ & c_1(\hfw(q,(2q+1)(qk)^2, 3(q+1)^2k+2))^{20}. 
\end{align*}

As a first step, we apply \zcref{thm_algogrid} to $X$ and obtain a $\hfw(q,(2q+1)(qk)^2, 3qk+1)$-wall $W_0.$
Next, we apply \zcref{lemma_HomogeneousFlatWall} to $W_0.$
This either results in a $K_{(2q+1)(qk)^2}$-minor model $\Xcal$ with $\Tcal_{\Xcal}\subseteq \Tcal_{W_0}\subseteq \Tcal_X,$ or in a set $Z\subseteq V(G)$ of size at most $16(2q+1)^3k^6$ together with a homogeneous flat wall $W_1\subseteq W_0$ in $G-Z$ of order $3(q+1)^2k+2.$

We treat these two cases separately.
\smallskip

\textbf{Case 1:}
We find a $K_{(2q+1)(qk)^2}$-minor model $\Xcal$ with $\Tcal_{\Xcal} \subseteq \Tcal_X.$
\smallskip

We now apply \zcref{thm_rainbowclique}.
This either yields a rainbow $K_{(qk)^2}$ as a colorful minor in $(G,\chi)$ or a set $A\subseteq V(G)$ of size at most $q(qk)^2 \leq 16(2q+1)^3(qk)^6$ such that the $\Tcal_{\Xcal}$-big component of $G-A$ is restricted.
With $\Tcal_{\Xcal}$ being a truncation of $\Tcal_X,$ the second outcome matches the second outcome of our assertion.
Hence we may assume to find a rainbow $K_{(qk)^2}$ as a colorful minor.
This, however, means that $(G,\chi)$ contains every $q$-colorful graph on $(qk)^2$ vertices as a colorful minor and this includes all $(q,k)$-segregated grids.
Hence, we are done with this case.
\medskip

\textbf{Case 2:}
We obtain the set $Z\subseteq V(G)$ of size at most $16(2q+1)^3k^6$ together with a homogeneous flat wall $W_1\subseteq W_0$ in $G-Z$ of order $3(q+1)^2k+2.$
\smallskip

We now layer our wall $W_1$ as follows.
Notice that for any $(k+2)$-wall $M$ with perimeter $P,$ $M-P$ contains a $k$-wall $M'$ that is a subwall of $M.$
If $k=k'+2$ for some positive $k'$ we may repeat this process.

Let now $C_1\dots,C_{(q+1)^2k}$ be the first $(q+1)^2k$ cycles obtained from $W_1$ by iteratively removing perimeters as described in the process above and let $W_2$ be the remaining $((q + 1)^{2}k + 2)$-wall.
We assume the cycles $C_i$ as above are numbered in order of their removal, that is, $C_1$ is the perimeter of $W_1.$
If we number the horizontal paths of $W_1$ starting from the top one as $Q_1,\dots,Q_{3(q+1)^2k+2},$ then the horizontal paths of $W_2$ are subpaths of $Q_{(q+1)^2k+1},\dots,Q_{2(q+1)^2k+2}.$
Let us select $S$ to be the set of the first vertex on $W_2$ from each $Q_i$ with $i\in[(q+1)^2k+1,(2(q+1)^2k)+2]$ when traversing along $Q_i$ starting on its left-most endpoint.

Now let $J\coloneqq \chi(V(H_B))$ for some brick $B$ of $W_2.$
Since $W_1$ is homogeneous, we know that also $W_2$ is homogeneous, so the choice of $B$ here is arbitrary.
Let $I\coloneqq [q]\setminus J$ and for each $i\in I$ let us denote by $Y_i$ the set $\chi^{-1}(i).$
We now apply \zcref{lemma_multicolorlinakge} to $G-Z-(V(W_2-S)),$ the sets $Y_i,$ $i\in I,$ and the set $S$ asking for $(q+1)k$ many paths for each member of $I.$
This yields one of two outcomes.
Either we find a linkage $\Pcal$ of order $q(q+1)k$ such that $\Pcal$ contains a $Y_i$-$S$-linkage $\Pcal_i$ of order $(q+1)k$ for each $i\in I,$ $\Pcal=\bigcup_{i\in I}\Pcal_i,$ and the paths in $\mathcal{P}$ are all internally vertex-disjoint from $W_2,$ or there is a set $Z'$ of size less than $(q+1)^2k$ such that there exists some $j\in I$ and there is no $Y_j$-$S$-path in $G-Z-Z'-(V(W_2-S).$

We discuss these two outcomes independently.
\smallskip

\textbf{Case 2.1:}
There exists a linkage $\Pcal$ of order $q(q+1)k$ such that $\Pcal$ contains a $Y_i$-$S$-linkage $\Pcal_i$ of order $(q+1)k$ for each $i\in I$ and $\Pcal=\bigcup_{i\in I}\Pcal_i.$
\smallskip

In this case we consider the left most vertical path $L$ of $W_2$ which contains all vertices of $S.$
We then partition $L$ into $q$ disjoint paths as follows.
Let $L_1$ be the shortest subpath of $L$ containing the top-most endpoint of $L$ such that there exists $\pi(1)\coloneqq i\in I$ and there are $k$ paths in $\Pcal_i$ with an endpoint in $L_1.$
Notice that for all $j\in I\setminus\{ i\},$ at most $k-1$ paths from $\Pcal_j$ have an endpoint on $L_1.$

Now suppose, for some $h\in[|I|]$ we have already found disjoint and consecutive subpaths $L_1,\dots,L_h$ of $L$ such that there is an injection $\pi\colon[h]\to I$ and for every $i\in[h],$ $k$ paths from $\Pcal_{\pi(i)}$ have an endpoint on $L_i.$
Moreover, we may assume that for each $j\in I$ with $\pi(i)\neq j$ for all $i\in[h],$ at least $(q+1)k-h(k-1)\geq (q-h+1)k$ paths of $\Pcal_j$ have an endpoint on $L-\bigcup_{i\in[h]}L_i.$

Next, let $L_{h+1}$ be the shortest subpath of $L-\bigcup_{i\in[h]}L_i$ containing a neighbor of $L_h$ such that there exists $j\in I,$ $j\neq\pi(i)$ for all $i\in[h],$ and $k$ paths from $\Pcal_j$ have their endpoints on $L_{h+1}.$
Clearly, by our choices of the $L_i$ so far, such a path $L_{h+1}$ and color $j$ exist.
Moreover, if we set $\pi(h+1)\coloneqq j$ we now have $h+1$ consecutive subpaths $L_i$ of $L$ such that for each $i\in[h+1],$ $k$ paths from $\Pcal_{\pi(i)}$ have their endpoint in $L_i,$ $\pi$ now is an injection from $[h+1]$ to $I,$ and for all remaining colors $x$ in $I\setminus \{ \pi(i) \mid i\in[h+1] \},$ there are at least $(q+1)k-(h+1)(k-1) \geq (q-h)k$ endpoints of paths from $\Pcal_x$ on $L-\bigcup_{i\in[h+1]}L_i.$

Hence, we find a bijection $\pi$ from $[|I|]$ to $I$ and linkages $\Pcal'_i$ of order $k$ for each $i\in I$ such that for every $j\in[|I|],$ all paths in $\Pcal'_{\pi(j)}$ have their endpoint on $L_j.$
Let now $W_3$ be the $qk$-subwall of $W_2$ obtained by removing, similar to before, all but $\nicefrac{1}{2} \cdot qk+1$ layers from $W_2.$
Notice that these layers give rise to a family $\Ncal$ of cycles where
\begin{align*}
 |\Ncal| = & \lfloor\frac{1}{2}((q+1)^2k - qk)\rfloor = \lfloor\frac{1}{2}((q^2+q+1)k)\rfloor \geq qk.
\end{align*}
Moreover, there are $qk$ pairwise vertex-disjoint paths -- subpaths of horizontal paths of $W_2$ -- that connect the left most vertical path of $W_3$ to the left most vertical path of $W_2$ and which have all their internal vertices in $W_2-W_3.$
Let now $S_1\subseteq V(L)$ be the set of all endpoints of the paths in $\bigcup_{i\in I}\Pcal'_{\pi(i)}.$
Moreover, let $T_1$ be the set of endpoints of the top most $|I|k$ horizontal paths of $W_3$ and the left most vertical path of $W_3.$

We claim that there exists an $S_1$-$T_1$-linkage $\Lcal$ of size $k|I|$ such that all internal vertices of the paths of $\Lcal$ belong to $W_2-W_3.$
To see this, notice that in the absence of such a linkage \zcref{prop_menger} would imply the existence of a separation of order at most $qk-1$ between $S_1$ and $T_1.$
However, each vertex of $S_1\cup T_1$ belongs to a horizontal path of $W_2$ that can be extended in at least $qk$ different ways, intersecting only in the horizontal path, to a path that intersects all cycles at least $qk$ of $\Ncal.$
Hence, such a separator must leave at least one cycle $N$ of $\Ncal,$ one way to get from some vertex $s$ of $S_1$ to $N,$ and one way to get from $N$ to $T_1$ intact.
Hence, deleting any set of size at most $qk-1$ still leaves at least one $S_1$-$T_1$-path intact and therefore, $S_1$ and $T_1$ cannot be separated by fewer than $qk$ vertices in $W_2-W_3.$

\begin{figure}[ht]
 \centering
 \scalebox{1}{
 \begin{tikzpicture}

 \pgfdeclarelayer{background}
		\pgfdeclarelayer{foreground}
			
		\pgfsetlayers{background,main,foreground}
			
 \begin{pgfonlayer}{main}
 \node (M) [v:ghost] {};

 \end{pgfonlayer}{main}

 \begin{pgfonlayer}{background}
 \pgftext{\includegraphics[width=10cm]{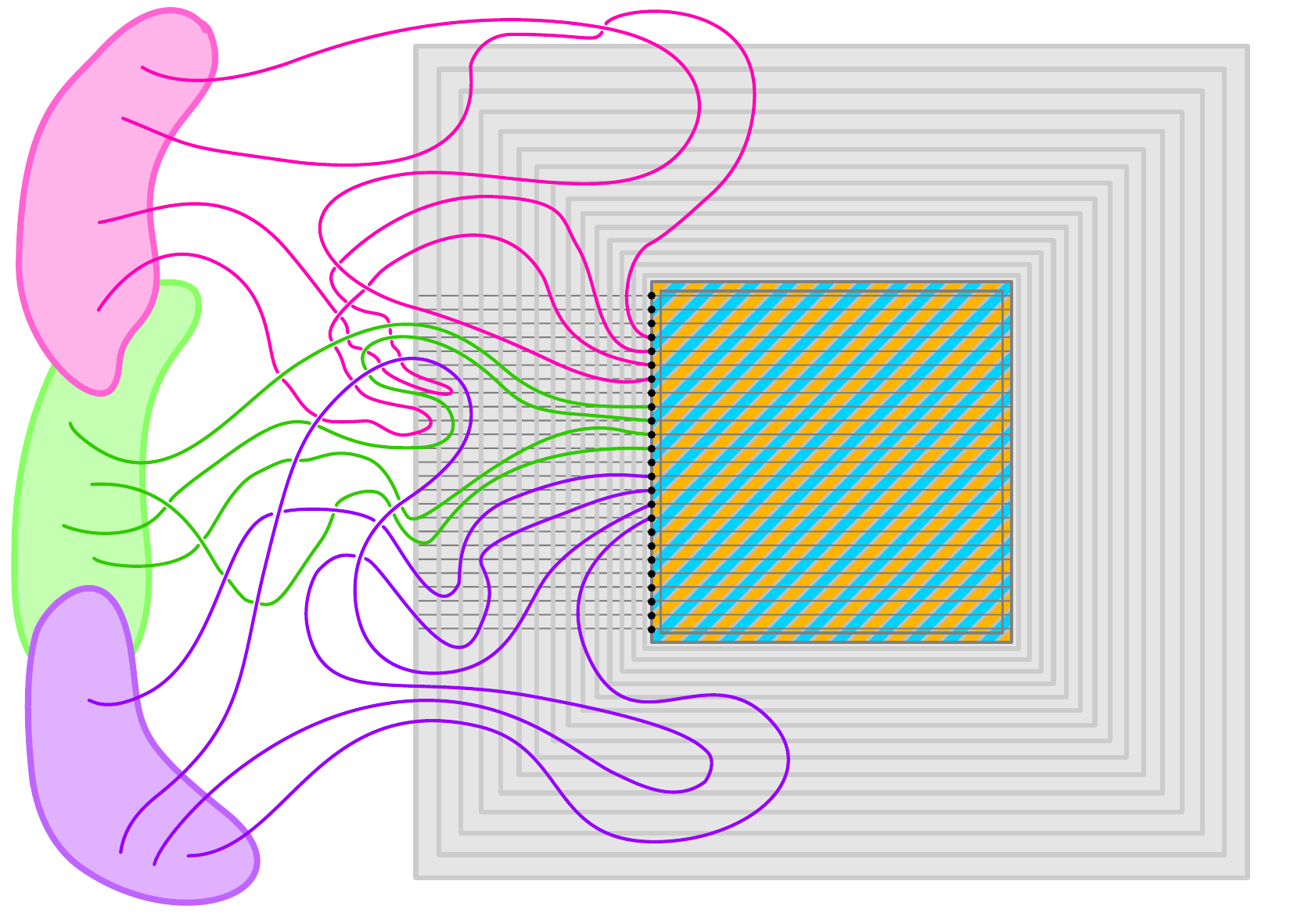}} at (M.center);
 \end{pgfonlayer}{background}
 
 \begin{pgfonlayer}{foreground}
 \end{pgfonlayer}{foreground}

 \end{tikzpicture}}
 \caption{A schematic representation of the construction of a $(5,4)$-segregated grid in the proof of \zcref{thm_localSegregation}.
 The striped box in the middle represents the wall $W_3$ whose bricks all represent colors $2$ and $3,$ while colors $1,$ $4,$ and $5$ are connected from the outside of $W_2$ to the perimeter of $W_3$ via disjoint paths.}
 \label{fig_segregatedGridConstruction}
\end{figure}

Let $D$ be the perimeter of $W_3$ and $W_4$ be the $qk$-subwall of $W_3$ that avoids $D.$
It is now easy to see that together with the paths from $\Pcal'_{\pi(i)},$ $i\in [|I|],$ the paths from $\Lcal,$ the horizontal paths of $W_3,$ and the bricks of $W_4$ that intersect $D,$ we may construct a $(q,k)$-segregated grid as a colorful minor of $(G,\chi).$
See \zcref{fig_segregatedGridConstruction} for an illustration.
\medskip

\textbf{Case 2.2:}
There exists a set $Z'$ of size less than $(q+1)^2k$ such that there exists some $j\in I$ and there is no $Y_j$-$S$-path in $G-Z-Z'-V(W_2-S).$
\smallskip

In this case we claim that the $\Tcal_{W_1}$-big component of $G-Z-Z'$ is restricted.
To see this, simply notice that we removed $(q+1)^2k$ cycles from $W_1$ to create $W_2.$
Moreover, there exist $(q+1)^2k$ pairwise vertex-disjoint horizontal paths of $W_1$ that contain a horizontal path of $W_2$ and therefore a vertex of $S.$
Since $|Z'|<(q+1)^2k,$ in $G-Z-Z'$ there must still be a path from the perimeter of $W_1$ to $S.$
Indeed, take any pair of one horizontal and one vertical path of $W_1$ that is disjoint from $Z',$ then there exists a path from this pair to $S$ in the $\Tcal_{W_1}$-big component of $G-Z-Z'.$
Moreover, there does exist such a path which is internally vertex-disjoint from $W_2.$
Thus, if the $\Tcal_{W_1}$-big component of $G-Z-Z'$ would contain a vertex $v$ with $j\in \chi(v),$ then there would still be a $Y_j$-$S$-path in $G-Z-Z'-V(W_2-S)$ which is impossible.
Hence, in this case the second outcome of our assertion holds.
\end{proof}

\subsection{A Grid Theorem for torso treewidth}\label{sec_globalSegregatedGrid}

We are now ready for the proof of \zcref{thm_restrictiveTreewidthIntro}.
Similar to our proofs of \zcref{thm_excludeRainbowClique,thm_ExcludingRainbowGridIntro} we do this by proving a slightly different and stronger statement which, in turn, is implied by a further strengthening that lends itself better to induction.

\begin{theorem}\label{thm_segregatedGridGlobal}
There exists a function $\sg_{\ref{thm_segregatedGridGlobal}}\colon\mathbb{N}^2\to\mathbb{N}$ such that for all positive integers $q$ and $k$ and every $q$-colorful graph $(G,\chi)$ one of the following holds:
\begin{enumerate}
 \item $(G,\chi)$ contains a $(q,k)$-segregated grid as a colorful minor and the torso treewidth of $(G, \chi)$ is at least $k - 1,$ or
 \item there exists a tree-decomposition $(T,\beta)$ for $(G,\chi)$ such that there is a set $L\subseteq V(T)$ of leaves of $T$ and
 \begin{itemize}
 \item $(T,\beta)$ has adhesion at most $\sg_{\ref{thm_segregatedGridGlobal}}(q,k),$
 \item $|\beta(t)|\leq\sg_{\ref{thm_segregatedGridGlobal}}(q,k)$ for all $t\in V(T),$ and
 \item for each $d \in L$ with neighbor $t \in V(T),$ $(G[\beta(d) \setminus \beta(t)], \chi)$ is restricted.
 \end{itemize}
\end{enumerate}
Moreover, $\sg(q,k)\in k^{2^{\Ocal(q)}}$ and there exists an algorithm that, given $q,$ $k,$ and $(G,\chi)$ as input, finds either a colorful minor model of a $(q,k)$-segregated grid in $(G,\chi)$ or the tree-decomposition as above in time $2^{k^{2^{\Ocal(q)}}} \cdot |G|^3|\!|G|\!|\cdot \log|G|.$
\end{theorem}

Notice that the tree-decomposition from the second outcome of \zcref{thm_segregatedGridGlobal} is a witness for the torso treewidth of $(G, \chi)$ being at most $\sg_{\ref{thm_segregatedGridGlobal}}(q,k).$
Hence, \zcref{thm_segregatedGridGlobal} does indeed act as an analogue of the Grid Theorem of Robertson and Seymour for torso treewidth.

We first prove that $(q,k)$-segregated grids indeed act as lower bounds on the torso treewidth of a graph.

\begin{lemma}\label{lemma_segregatedLowerbound}
For all positive integers $q$ and $k,$ and every $q$-colorful graph $(G,\chi),$ if $(G,\chi)$ contains a $(q,k)$-segregated grid as a colorful minor, then the torso treewidth of $(G, \chi)$ is at least $k - 1.$
\end{lemma}
\begin{proof}
Suppose, towards a contradiction, that $(G,\chi)$ contains a minor model $\Xcal$ of some $(q,k)$-segregated grid and, at the same time, that there exists a set $X\subseteq V(G)$ such that $\tw(G_X)<k-1$ where $G_X$ denotes the torso of $X$ in $G.$

Let $(T,\beta)$ be a tree-decomposition of width at most $k-2$ for $G_X$ and let $dt\in E(T).$
For each $x\in \{ d,t\}$ let $T_x$ be the component of $T-dt$ containing $x$ and denote by $\beta(T_x)$ the set $\bigcup_{y\in V(T_x)}\beta(y).$
It follows directly from the definition of tree-decompositions that $(\beta(T_d),\beta(T_t))$ is a separation of order at most $k-2$ in $G_x.$

Recall that for every clique $K$ of a graph $H$ and every tree-decomposition $(F,\gamma)$ of $H$ there must exist $x\in V(F)$ such that $V(K)\subseteq \gamma(x).$
Hence, for each component $J$ of $G-X$ there exists $x\in V(T)$ such that $N(J)\subseteq\beta(x).$
This observation allows us to make a fundamental assumption on the structure of $(T,\beta).$
From now on we will assume that for every component $J$ of $G-X$ there exists a leaf $l_J\in V(T)$ such that $\beta(l_J) = N(J).$

For each $x\in \{ d,t\}$ let $A^{dt}_x$ denote the union of $\beta(T_x)$ and the vertex sets of all components $J$ of $G-X$ such that $l_J\in T_x.$
Then $(A^{dt}_d,A^{dt}_t)$ is a separation of order at most $k-2$ in $G.$

Since $\Xcal$ is a minor model of the $(qk\times qk)$-grid in particular it follows that $\Tcal_{\Xcal}$ is a tangle of order $qk > k-1$ in $G.$
Hence, we may assign to each edge $xy\in E(T)$ with $(A^{xy}_x,A^{xy}_y)\in\Tcal_{\Xcal}$ the orientation $(x,y)$ and thereby obtain an orientation $\vec{T}$ of $T.$

Suppose first that for every leaf $x$ of $T$ with unique neighbor $y$ we have $(x,y)\in E(\vec{T}).$
In this case there must exist a unique sink $s\in V(T)$ that is not a leaf in $\vec{T}.$
Notice that this means, in particular, that for all pairs of components $J,J'$ of $G-X,$ if $l_J$ and $l_{J'}$ belong to different components of $T-s,$ then $J$ and $J'$ are not connected in $G-\beta(s).$
This is, because by the tangle axioms it follows that no vertex of $\vec{T}$ can have two outgoing edges.
Let $t_1,\dots,t_{\ell}$ be the neighbors of $s.$
Notice that for every $i\in[\ell]$ it holds that $(A^{st_i}_{t_i}\cup\beta(s),A^{st_i}_s)\in\Tcal_{\Xcal}.$
Hence, if $\ell \leq 3$ we have that $\bigcup_{i\in[\ell]}G[A^{st_i}_{t_i}\cup \beta(s)]=G$ which contradicts $\Tcal_{\Xcal}$ being a tangle.
So we may assume that $\ell \geq 4.$

We prove by induction on $j\in[\ell-2]$ that for any subset $I\subseteq[\ell]$ of order at most $j$ we have $(\beta(s)\cup\bigcup_{i\in I}A^{st_i}_{t_i},\bigcap_{i\in I}A^{st_i}_{s})\in\Tcal_{\Xcal}.$
For $j=1$ this is trivially true by choice of $s$ and the definition of $\vec{T}.$
So assume the claim holds for all $j' <j$ and let $I\subseteq[\ell]$ be any set of size $j.$
In case we have that $(\bigcap_{i\in I}A^{st_i}_{s},\beta(s)\cup\bigcup_{i\in I}A^{st_i}_{t_i})\in\Tcal_{\Xcal}$ we know, by induction that for any $h\in I$ we have
\begin{align*}
(\beta(s)\cup\bigcup_{i\in I\setminus \{h\}}A^{st_i}_{t_i},\bigcap_{i\in I\setminus\{ h\}}A^{st_i}_{s}) &\in \Tcal_{\Xcal}\text{ and}\\
(A^{t_hs}_{t_h}\cup \beta(s),A^{t_hs}_s) &\in \Tcal_{\Xcal}.
\end{align*}
Moreover, it holds that
\begin{align*}
G[\big(\beta(s)\cup\bigcup_{i\in I}A^{st_i}_{t_i}\big) \setminus \big(\beta(s)\cup\bigcup_{i\in I\setminus \{h\}}A^{st_i}_{t_i}\big)] \subseteq G[A^{t_hs}_{t_h}].
\end{align*}
Hence, we obtain that
\begin{align*}
G[\bigcap_{i\in I}A^{st_i}_{s}] \cup G[\beta(s)\cup\bigcup_{i\in I\setminus \{h\}}A^{st_i}_{t_i}] \cup G[A^{t_hs}_{t_h}] = G
\end{align*}
which is a contradiction to $\Tcal_{\Xcal}$ being a tangle.

Consequently, we may now assume that there exists a leaf $l\in V(T)$ which is the unique sink of $\vec{T}.$
In case there does not exist a component of $G-X$ such that $l=l_J$ we may apply a similar argument as above and obtain a contradiction to $\Tcal_{\Xcal}$ being a tangle.
Hence, there exists a component $J$ of $G-X$ such that $l=l_J.$
Let $t\in V(T)$ be the unique neighbor of $l.$
Then $J$ is the $\Tcal_{\Xcal}$-big component of $G-(\beta(l)\cup\beta(t)).$
Thus, there exists some vertical path $P$ of a $(q,k)$-segregated grid such that there is a path $P'$ in $J$ which meets each $X_v\in\Xcal$ with $v\in V(P).$

Indeed, it follows that for each $i\in [q]$ there exists a set $Y_i$ of $k$ vertices $y$ with $i\in\chi(y)$ such that there is a $V(P')$-$Y_i$-linkage of size $k$ in $G.$
Hence, for every color $i\in[q]$ there exists a $V(J)$-$\chi^{-1}(i)$-linkage of size $k$ in $G.$
However, by our assumption we have that $J$ is restricted and $|\beta(l)\cap \beta(t)|\leq k-1.$
Since this is clearly absurd we have found our final contradiction and our proof is complete.
\end{proof}

All that is left to do is to prove the following strengthening of the upper bound part of \zcref{thm_segregatedGridGlobal}. 

\begin{theorem}\label{thm_segregatedGridGlobal_Induction}
There exists a function $\sg_{\ref{thm_segregatedGridGlobal_Induction}}\colon\mathbb{N}^2\to\mathbb{N}$ such that for all positive integers $q$ and $k,$ every $q$-colorful graph $(G,\chi),$ and every set $X\subseteq V(G)$ of size at most $3\sg_{\ref{thm_segregatedGridGlobal_Induction}}(q,k)+1$ one of the following holds:
\begin{enumerate}
 \item $(G,\chi)$ contains a $(q,k)$-segregated grid as a colorful minor, or
 \item there exists a rooted tree-decomposition $(T,r,\beta)$ for $(G,\chi)$ such that there is a set $L\subseteq V(T)$ of leaves of $T$ and
 \begin{itemize}
 \item $X\subseteq V(G)$
 \item $(T,\beta)$ has adhesion at most $4\sg_{\ref{thm_segregatedGridGlobal_Induction}}(q,k)+1,$
 \item $|\beta(t)|\leq 4\sg_{\ref{thm_segregatedGridGlobal_Induction}}(q,k)+1$ for all $t\in V(T),$ and
 \item for each $d\in L$ there exists a neighbor $t\in V(T)$ of $d$ such that $(G[\beta(d)\setminus\beta(t)],\chi)$ is restricted.
 \end{itemize}
\end{enumerate}
Moreover, $\sg(q,k)\in k^{2^{\Ocal(q)}}$ and there exists an algorithm that, given $q,$ $k,$ $X,$ and $(G,\chi)$ as input, finds one of these outcomes in time $2^{k^{2^{\Ocal(q)}}} \cdot |G|^3|\!|G|\!| \log|G|.$
\end{theorem}

\begin{proof}
We quickly go through the part of the proof that is analogous to previous variants of a \say{local-to-global}-type argument.

First, the function $\sg_{\ref{thm_segregatedGridGlobal_Induction}}$ is defined as follows.
\begin{align*}
 \sg_{\ref{thm_segregatedGridGlobal_Induction}}(q,k) \coloneqq \apex_{\ref{thm_localSegregation}}(q,k) + \link_{\ref{thm_localSegregation}}(q,k).
\end{align*}

We now proceed by induction on $|V(G)\setminus X|$ and start by discussing two fundamental cases.
Exactly as in our proof of \zcref{thm_RainbowClique_induction}, we first resolve the two \textbf{fundamental cases} where either $G$ has only few vertices, or $|X|<3\sg_{\ref{thm_segregatedGridGlobal_Induction}}(q,k)+1$ in precisely the same way as before.
Then two further cases follow:

\textbf{Case 1}, where $X$ has a small balanced separator, and \textbf{Case 2}, where $X$ does not have a small balanced separator.
As before, for \textbf{Case 1} there is no difference to the proof of \zcref{thm_RainbowClique_induction} and thus, we immediately move on to \textbf{Case 2}.

This leaves is with the last case.
\smallskip

\textbf{Case 2:}
The set $X$ is $(\sg_{\ref{thm_segregatedGridGlobal_Induction}}(q,k),\nicefrac{2}{3})$-linked.
\smallskip

In contrast to previous iterations of this proof, here we do not have to distinguish more subcases.
We immediately may apply \zcref{thm_localSegregation} to $(G,\chi)$ and $X$ which either yields a $(q,k)$-segregated grid as a colorful minor in $(G,\chi)$ -- and there would be nothing more to prove -- or it provides us with a set $A\subseteq V(G)$ with $|A| \leq \apex_{\ref{thm_localSegregation}}(q,k) < \sg_{\ref{thm_segregatedGridGlobal_Induction}}(q,k)$ such that the $\Tcal_X$-big component $K$ of $G-A$ is segregated.

In this second case we may immediately observe that $|V(J)\cap X| \leq \sg_{\ref{thm_segregatedGridGlobal_Induction}}(q,k)$ for every component $J$ of $G-A$ other than $K.$
For every such $J$ where there exists a vertex $v_J \in V(J)\setminus X$ we set $X_J\coloneqq A \cup \{ v_J\}\cup (V(J)\cap X).$
For all other components $J$ we set $X_J \coloneqq A \cup (V(J)\cap X).$
Moreover, for each component $J$ of $G-A$ other than $K,$ we set $H_J\coloneqq G[A\cup V(J)].$

Notice that we have $|X_J|\leq 3\sg_{\ref{thm_segregatedGridGlobal_Induction}}(q,k)+1$ and $|V(H_J)\setminus X_J| < |V(G)\setminus X|$ in all cases.
Hence, we may now apply our induction hypothesis to every pair $(H_J,\chi)$ and $X_J.$
This either yields a $(q,k)$-segregated grid as a colorful minor of $(G,\chi),$ or for each $J$ we find a rooted tree-decomposition $(T_J,r_J,\beta_J)$ satisfying the properties required by the assertion and with $X_J\subseteq \beta_J(r_J).$

Let now $T$ be the tree obtained from the disjoint union of all $T_J$ by introducing a new vertex $t_K,$ a new vertex $r,$ and then joining $r$ with an edge to $t_K$ and all $r_J.$
Finally, we set $\beta(t_K) \coloneqq V(K)\cup X\cup A,$ $\beta(r)\coloneqq X\cup A,$ and $\beta(x)\coloneqq \beta_J(x)$ for every $J$ and $x\in V(T_J).$
It follows directly from our construction that $(T,r,\beta)$ is a rooted tree-decomposition precisely as required by the assertion.
\end{proof}

\section{Algorithmic applications}\label{sec_Applications}

To present the algorithmic applications of our results, we first introduce a unifying concept that encompasses all of them.

A graph parameter $\p$ is called \emph{decomposable} if it is subgraph-monotone and, for every separation $(A,B)$ of $G,$ we have $\p(G) \leq \max\{\p(\torso(G,A)), |B|\}.$
Examples of decomposable parameters considered in this paper include treewidth $\tw$ and the Hadwiger number $\hw.$

\subsection{Star Decompositions of Colorful Graphs}\label{subsec_StarDecompositions}
Given a tree $T$ with at least one vertex, we denote by $L(T)$ the set of its leaves.
A \emph{star} is any connected subgraph of $K_{1,t}$ for some non-negative integer $t.$
If $|L(T)| \geq 2,$ then $T$ has exactly one non-leaf vertex, called its \emph{center}.
If $|L(T)| \leq 1,$ any vertex of $T$ may be designated as the center.

A \emph{star decomposition} of a $q$-colorful graph $(G,\chi)$ is a tree decomposition $\Tcal=(T,\beta)$ where $T$ is a star with center $v^{\star}$ and $\chi^{-1}([q])\subseteq \beta(v^{\star}),$ i.e. all vertices with non-empty palettes appear in the center's bag.
Without loss of generality, we assume $T$ is chosen so that $L(T)$ is maximized.
This ensures that for each $\ell \in L(T),$ the set $\beta(\ell) \setminus \beta(v^{\star})$ is connected in $G.$ 

For a graph parameter $\p,$ the \emph{$\p$-width} of $\Tcal$ is defined as
\begin{align*}
\max\left\{\p\left(\torso(G,\beta(v^\star))\right)\right\} \cup \left\{|\beta(\ell)\cap\beta(v^{\star})| \mid \ell\in L(T)\right\}.
\end{align*}
The following observation follows by iteratively applying the definition of parameter decomposability to each leaf in a star decomposition.

\begin{observation}\label{obs_boundPWidth}
Let $\p$ be a decomposable parameter.
If a graph $G$ admits a star decomposition $(T,\beta)$ of $\p$-width at most $k,$ then 
\begin{align*}
\p(G) \leq \max\{k\} \cup \{|\beta(\ell)| \mid \ell\in L(T)\}.
\end{align*}
\end{observation}

The next two results follow from \zcref{thm_segregatedGridGlobal} and \zcref{thm_excludeRainbowClique}, respectively, applied to $1$-colorful graphs.
These structural results can be stated as follows:

\begin{corollary}\label{thm_SRTW}
There exists a function $f_{\ref{thm_SRTW}} \colon \Nbbb \to \Nbbb$ and an algorithm that, given a $1$-colorful graph $(G,\chi)$ and an integer $k\in\Nbbb,$ either outputs a colorful minor model of a $(1,k)$-segregated grid or a star decomposition $(T,\beta)$ of $(G,\chi)$ with $\tw$-width at most $f_{\ref{thm_SRTW}}(k),$ in $\Ocal_{k}(|G|^{3}| \!|G|| \log|G|)$ time.
\end{corollary}

\begin{corollary}\label{thm_SRC}
There exists a function $f_{\ref{thm_SRC}}\colon\Nbbb\to\Nbbb$ and an algorithm that, given a $1$-colorful graph $(G,\chi)$ and an integer $k\in\Nbbb,$ either outputs a colorful minor model of a $1$-colorful rainbow $K_{k}$ or a star decomposition $(T,\beta)$ of $(G,\chi)$ with $\hw$-width at most $f_{\ref{thm_SRC}}(k),$ in $\Ocal_{k}(|G|^{3}| \!|G|| \log|G|)$ time.
\end{corollary}

In the following sections, we present algorithmic applications of \zcref{thm_SRTW} and \zcref{thm_SRC}.

\subsection{An algorithmic reduction}\label{subsec_algoReduction}

We now introduce the main algorithmic machinery of this section.
It consists of an algorithmic reduction that relies on two key algorithmic assumptions and the existence of a star decomposition along with a graph parameter $\p$ associated with it.
The first assumption is that, for the problem under consideration, there exists an efficient method for replacing boundaried subgraphs of a colorful graph with equivalent but smaller ones.
The second assumption is that the problem can be solved efficiently when the value of $\p$ is bounded.

 \paragraph{Equivalences and representations.}

For every non-negative integer $q,$ we denote the class of all $q$-colorful graphs by $\Ccal^{q}_{\mathsf{all}}.$ 
We define a \emph{functional problem}\footnote{We use the term “functional problem on $q$-colorful graphs” instead of “$q$-colorful graph parameter” to emphasize that we are addressing a computational problem.} on $q$-colorful graphs to be a map $\funp\colon\Ccal^{q}_{\mathsf{all}}\to\Nbbb.$ 
Our goal is to identify conditions that enable the design of algorithms for solving functional problems on $q$-colorful graphs.

We say that two boundaried $q$-colorful graphs $\mathbf{G}$ and $\mathbf{G}'$ are \emph{$\funp$-equivalent} if they are compatible and, for every boundaried $q$-colorful graph $\mathbf{C}$ that is compatible with $\mathbf{G}$ (and therefore with $\mathbf{G}'$ as well), it holds that $$\funp(\mathbf{C} \oplus \mathbf{G})=\funp(\mathbf{C}\oplus \mathbf{G}'), $$
where given two compatible boundaried $q$-colorful graphs $(G_{1},\chi_{1}, B_{1}, \rho_{1})$ and $(G_{2}, \chi_{2}, B_{2}, \rho_{2}),$ $(G_{1},\chi_{1}, B_{1}, \rho_{1}) \oplus (G_{2}, \chi_{2}, B_{2}, \rho_{2})$ denotes the boundaried $q$-colorful graph obtained by identifying the vertex $\rho_{1}(i)$ with the vertex $\rho_{2}(i),$ for every $i \in [|B_{1}|].$

A boundaried $q$-colorful graph $(G,\chi,B,\rho)$ is \emph{confined} if 
$\chi(G)\subseteq B.$
A \emph{representation} of $\funp$ is a pair $(\rep, \bnd),$ where $\rep$ maps \textsl{confined} boundaried $q$-colorful graphs to $\funp$-equivalent \textsl{confined} boundaried $q$-colorful graphs, and $\bnd$ maps integers to integers such that, for every confined boundaried $q$-colorful graph $(G, \chi, B, \rho),$ the graph $\rep(G, \chi, B, \rho)$ has at most $\bnd(|B|)$ vertices.
We emphasize that the notion of \emph{confinedness} plays a crucial role in both our definitions and results.

\begin{theorem}\label{thm_reduceToStars}
Let $q\in\Nbbb,$ $\funp\colon \Ccal^{q}_{\mathsf{all}}\to\Nbbb$ be a functional problem on $q$-colorful graphs, and $\p\colon\gall\to\Nbbb$ be a decomposable graph parameter.
Suppose that there exists functions $f^\mathsf{a}, f^\mathsf{r} \colon \Nbbb\to\Nbbb$ such that the following hold:
\begin{enumerate}
\item There exists an algorithm that, given a $q$-colorful graph $(G,\chi),$ computes $\funp(G,\chi)$ in $\Ocal_{q+\p(G)+\mathsf{m}(G,\chi)}(f^\mathsf{a}(|G|))$ time.
\item There exists a representation $(\rep,\bnd)$ of $\funp$ and algorithms that, given a boundaried $q$-colorful graph $(G, \chi, B, \rho),$ computes $\rep(G, \chi, B, \rho)$ and $\bnd(|B|)$ in $\Ocal_{q+|B|}(f^\mathsf{r}(|G|))$ and $\Ocal_{q+|B|}(1)$ time respectively.
\end{enumerate}
Then, there exists an algorithm that, given a star decomposition $\Tcal = (T,\beta)$ of a $q$-colorful graph $(G,\chi)$ with $\p$-width at most $k,$ outputs $\funp(G,\chi)$ in $\Ocal_{q+k+\mathsf{m}(G,\chi)}(f^\mathsf{a}(|G|)+f^\mathsf{r}(|G|))$ time.
\end{theorem}
\begin{proof}

Let $v^{\star}$ be the center of $T.$
For each leaf $\ell$ of $T,$ consider the boundaried $q$-colorful graph $\mathbf{G}^{\ell}=(G[\beta(\ell)],\chi,\beta(\ell)\cap \beta(v^{\star}),\rho^\ell),$ where $\rho^\ell$ is a bijection from $\beta(\ell)\cap \beta(v^\star)$ to $[|\beta(\ell)\cap \beta(v^\star)|].$
Notice that $\mathbf{G}^{\ell}$ is confined by definition of a star decomposition.
Using the algorithm from (ii), compute $\mathbf{G}^{\ell}_{\text{new}} \coloneqq \rep(\mathbf{G}^{\ell}).$
Initialize $(G',\chi')\coloneqq (G,\chi).$
For each leaf $\ell,$ update
\begin{align*}
(G',\chi')\coloneqq (G'-(\beta(\ell)\setminus \beta(v^*)),\chi',\beta(\ell)\cap \beta(v^\star),\rho^\ell)\oplus \mathbf{G}^{\ell}_{\text{new}},
\end{align*}
effectively replacing $\mathbf{G}^{\ell}$ with $\mathbf{G}^{\ell}_{\text{new}}$ in $(G',\chi').$
By $\funp$-equivalence it holds that $\funp(G,\chi)=\funp(G',\chi).$
This construction runs in $\Ocal_{q+k}(f^\mathsf{r}(|G|))$ time.

Now consider the star decomposition $(T, \beta')$ of $(G',\chi'),$ where $\beta'(\ell)$ is the vertex set of $\mathbf{G}_{\text{new}}^{\ell}$ for each $\ell\in L(T).$
By \zcref{obs_boundPWidth} we know that $\p(G') \leq k'\coloneqq\max\{k,\bnd(k)\}.$
Using the algorithm from (i), we compute $\funp(G',\chi')$ (which is equal to $\funp(G,\chi)$) in $\Ocal_{q + k'}(f^\mathsf{a}(|G|)) = \Ocal_{q + k}(f^\mathsf{a}(|G|))$ time, as required.
\end{proof}

All subsequent meta-algorithmic results will rely on \zcref{thm_reduceToStars}.

\subsection{Logics and AMTs on colorful graphs}\label{subsec_ATMs}
We now present the logical frameworks for two Algorithmic Meta-Theorems (AMTs) that we will extend to the colorful setting.

 \paragraph{\CMSO Logic on Colorful Graphs.}
We consider logical formulas on colorful graphs with the following variables:
A vertex set variable $\mathsf{V}$ and an edge set variable $\mathsf{E}.$
Also we consider a collection of set variables $\mathsf{X}_{1}, \mathsf{X}_{2}, \ldots$ representing colors.
A $q$-colorful graph $(G,\chi)$ is represented as the structure $\mathbf{G} = (G=(V,E), X_{1},\ldots,X_{q}),$ where:
$V$ interprets $\mathsf{V},$ $E$ interprets $\mathsf{E},$ and $X_{i} = \chi^{-1}(i)$ for $i\in[q]$ interprets $\mathsf{X}_{i}.$
These variables are always free in $\mathbf{G}.$
We also consider:
Sequences of vertex/edge variables $\mathsf{v}_{1}, \mathsf{v}_{2}, \ldots$ / $\mathsf{e}_{1}, \mathsf{e}_{2}, \ldots$ and vertex/edge set variables $\mathsf{V}_{1}, \mathsf{V}_{2}, \ldots$ / $\mathsf{E}_{1}, \mathsf{E}_{2}, \ldots,$ which are quantified in formulas.
These formulas describe structures $(G, X_{1}, \ldots, X_{q}, \overline{v}, \overline{e}, \overline{V}, \overline{E}),$ where $\overline{v}, \overline{e}, \overline{V}, \overline{E}$ are \emph{indexed sets} interpreting the corresponding variables.
For instance, $\overline{v} = \{v_{i} \mid i\in I\}$ for some $I \subseteq \Nbbb_{\geq 1},$ where $v_{i}$ interprets $\mathsf{v}_{i}.$
We call $I$ the \emph{index} of $\overline{v},$ with analogous definitions for $\overline{e}, \overline{V},$ and $\overline{E}.$

A \emph{monadic second-order logic (\textsf{MSO}) formula} may be:
\begin{itemize}
 \item An atomic formula: 
 \begin{itemize}
 \item \say{$\mathsf{v}_{i} = \mathsf{v}_{j}$} (vertex equality),
 \item \say{$\mathsf{v}_{i} \sim \mathsf{v}_{j}$} (vertex adjacency),
 \item \say{$\mathsf{v}_{i} \in \mathsf{V}_{j}$} or \say{$\mathsf{e}_{i} \in \mathsf{E}_{j}$} (membership).
 \end{itemize}
 \item A compound formula built using $\vee,$ $\neg,$ and $\exists$ (quantifying over vertex/edge/set variables).
\end{itemize}

For a structure $\mathbf{G} = (G, X_{1}, \ldots, X_{q}, \overline{v}, \overline{e}, \overline{V}, \overline{E})$ and an index $i$ not in $\overline{v}$'s index, we say $\mathbf{G} \models \exists \mathsf{v_{i}} \phi$ if there exists a vertex $v_i \in V(G)$ such that $(G, X_{1}, \ldots, X_{q}, \overline{v} \cup \{v_i\}, \overline{e}, \overline{V}, \overline{E})$ models $\phi$ with $\mathsf{v}_i$ interpreted as $v_i.$
Quantifiers $\exists \mathsf{e}_i,$ $\exists \mathsf{V}_i,$ and $\exists \mathsf{E}_i$ are defined similarly.
We write $\forall \mathsf{v}_i \phi$ as shorthand for $\neg (\exists \mathsf{v}_{i} (\neg \phi)),$ with analogous definitions for $\forall \mathsf{e}_i \phi,$ $\forall \mathsf{V}_i \phi,$ and $\forall \mathsf{E}_i \phi.$

A \emph{counting monadic second-order logic (\textsf{CMSO}) formula} extends \textsf{MSO} with predicates $\mathsf{card}_p(\mathsf{V}_{i}),$ which hold when $|\mathsf{V}_i|$ is divisible by $p \in \Nbbb_{\geq 1}.$
For a colorful graph $(G,\chi)$ and \textsf{CMSO}-formula $\phi,$ we write $(G,\chi) \models \phi$ if $(G, X_{1}, \ldots, X_{q})$ models $\phi.$

 \paragraph{The $\mathsf{CMSO/tw}$ Logic.}
The logic $\mathsf{CMSO/tw}$ modifies $\mathsf{CMSO}$ by replacing:
$\exists \mathsf{V}_i,$ $\exists \mathsf{E}_i$ with $\exists_{k} \mathsf{V}_i,$ $\exists_{k} \mathsf{E}_i$ ($k \in \mathbb{N}$) and $\forall \mathsf{V}_i,$ $\forall \mathsf{E}_i$ with $\neg (\exists_k \mathsf{V}_{i} \neg \phi),$ $\neg (\exists_k \mathsf{E}_{i} \neg \phi).$
We next define the special quantifier $\exists_{k}$: 
For a structure $\mathbf{G} = (G, X_{1}, \ldots, X_{q}, \overline{v}, \overline{e}, \overline{V}, \overline{E})$ and index $i$ not in $\overline{V}$'s index, $\mathbf{G} \models \exists_{k} \mathsf{V_{i}} \phi$ if there exists a vertex set $V_i \subseteq V(G)$ such that $(G, X_{1}, \ldots, X_{q}, \overline{v}, \overline{e}, \overline{V} \cup \{V_i\}, \overline{E}) \models \phi$ and, moreover, $\bidim(G,V_{i}) \leq k.$ 
The quantifier $\exists_{k} \mathsf{E_{i}} \phi$ is defined analogously, requiring $\bidim(G,V_{E_{i}}) \leq k$ for the edge set $E_i$'s endpoints $V_{E_{i}}.$\footnote{In the definition of \cite{SauST25Parameterizing}
it is demanded that $\tw(G,V_{i})\leq k$ that is the maximum treewidth of an $X$-rooted minor of $G.$ Clearly, $\tw$ and $\bidim$ can be seen as $1$-colorful graph parameters that are equivalent because 
of the main result in \cite{chuzhoy2021towards}. In this paper we prefer to define 
$\mathsf{CMSO/tw}$ using $\bidim$ instead of $\tw$ as this is more adequate to the colorful setting of the paper.} Also, given a $\mathsf{CMSO/tw}$ formula $\phi,$ we use the notation $\mathsf{k}(\phi)$ in order 
to denote the bidimensionality bound that it imposes on the interpretation of $\exists_{k} \mathsf{V}_i,$ $\exists_{k} \mathsf{E}_i.$ In order to unify notation, in case of a \CMSO-formula $\phi,$ where there is no restriction on quantification, we assume that $\mathsf{k}(\phi)=\infty$

 \paragraph{Disjoint Paths Extensions.}
Given a $k \in \mathbb{N}_{\geq 1},$ we next define the $2k$-ary \emph{disjoint-paths predicate} $\mathsf{dp}_k(\mathsf{x}_1, \mathsf{y}_1, \ldots, \mathsf{x}_k, \mathsf{y}_k).$
The logic $\mathsf{CMSO/tw}\!+\!\mathsf{dp}$ extends $\mathsf{CMSO/tw}$ by allowing atomic formulas of this form.
For a structure $\mathbf{G},$ $\mathsf{dp}_k(\mathsf{x}_1, \mathsf{y}_1, \ldots, \mathsf{x}_k, \mathsf{y}_k)$ holds if there exist vertex-disjoint paths $P_1, \ldots, P_k$ in $G$ connecting each $\mathsf{x}_i$ to $\mathsf{y}_i.$

We define $\CMSO$ logic as the set of all \textsf{CMSO}-formulas and, similarly, we define $\mathsf{CMSO}\!+\!\mathsf{dp}$ logic.
We stress that $\mathsf{CMSO}\!+\!\mathsf{dp}$ has the same expressibility as $\mathsf{CMSO}$ as the disjoint paths predicate can be expressed by $\mathsf{CMSO}$-formulas.
However, $\mathsf{CMSO/tw}$ is a proper subset of $\mathsf{CMSO/tw}\!+\!\mathsf{dp}$ which is more expressive as the disjoint paths predicate is not expressible in $\mathsf{CMSO/tw}.$

 \paragraph{\CMSO-Definable Functional Problems on $q$-Colorful Graphs.}
Let $q$ be a non-negative integer, and let $\phi$ be an \CMSO-formula over $(q+1)$-colorful graphs.
We define the \emph{maximization problem} on $\phi$ as the functional problem $\opt_{\phi}^{\max}(G,\chi) \colon \Ccal_{\mathsf{all}}^{q}\to\Nbbb,$ where
%
\begin{align}
\opt_{\phi}^{\max}(G,\chi)& = \max\{|X|\mid X\subseteq V(G), \bidim(G,X)\leq \mathsf{k}(\phi), \text{~and~} (G,\chi\!+\!X)\models \phi\}\cup\{0\},\label{def_eq_opt}
\end{align}
where $\chi\!+\!X = \chi\cup\{(x,q+1)\mid x\in X\},$ i.e., we add a new color $q+1$ to all vertices of $X.$
We also define the \emph{minimization problem} $\opt_{\phi}^{\min}(G,\chi)$ analogously, by replacing \say{$\max$} with \say{$\min$} and \say{$\{0\}$} by \say{$\{\infty\}$} in the above.

In general, when we do not need to specify $\bullet\in\{\max,\min\},$ we simply use the term \emph{optimization problem} (on $\phi$).
Notice that in \eqref{def_eq_opt} we insist that the vertex set $X$ carrying the \say{new color} is a subset of the set of vertices that already carry some color.

The following result is the optimization extension of Courcelle's Theorem (\cite{Courcelle1990Monadic, Courcelle1992Monadic}) proved by Arnborg and Andrzej Proskurowski \cite{ArnborgP1989Linear}, stated in terms of colorful graphs (see also \cite{BoriePT1992Automatic,ArnborgLS1991Easy}).

\begin{proposition}\label{prop_TwAnnotated}
There exists an algorithm that, given $q\in\Nbbb,$ $\bullet\in\{\min,\max\},$ a formula $\phi\in\mathsf{CMSO}$ over $(q + 1)$-colorful graphs, and a $q$-colorful graph $(G,\chi),$ computes $\opt_{\phi}^{\bullet}(G,\chi)$ in time $\Ocal_{q+|\phi|+\tw(G)}(|G|).$
\end{proposition}

Sau, Stamoulis, and Thilikos \cite{SauST25Parameterizing} proved the following Algorithmic Meta-Theorem (AMT).
While in its full generality the result is stated for general structures, we present it here for colorful graphs.
Actually, we present the following optimization version of it, that can be easily derived from the results of \cite{SauST25Parameterizing} (for this optimization version, it is important to demand $\bidim(G,X)\leq \mathsf{k}(\phi)$ in the definition of $\opt_{\phi}^{\bullet}(G,\chi)$ in \eqref{def_eq_opt}).
 

\begin{proposition}\label{prop_HadAnnotated}
There exists an algorithm that, given $q \in \Nbbb,$ $\bullet\in\{\min,\max\},$ a formula $\phi\in\mathsf{CMSO/tw}\!+\!\mathsf{dp}$ over $(q+1)$-colorful graphs, and a $q$-colorful graph $(G,\chi),$ computes $\opt_{\phi}^{\bullet}(G,\chi)$ in time $\Ocal_{q + |\phi|+\hw(G)}(|G|^2).$ 
\end{proposition}

We say that a $q$-colorful graph parameter $\funp\colon\Ccal^{q}_{\mathsf{all}}\to\Nbbb$ is \emph{\CMSO-definable} by some \MSO-formula $\phi$ on $(q+1)$-colorful graphs and some choice of $\bullet\in\{\max,\min\},$ if $\funp=\opt_{\phi}^{\bullet}.$
If $\phi$ is a $\mathsf{CMSO/tw}\!+\!\mathsf{dp}$-formula, then we say that $\funp$ is $\mathsf{CMSO/tw}\!+\!\mathsf{dp}$-definable (by $\phi$ and $\bullet$).

By \zcref{prop_TwAnnotated} (resp. \zcref{prop_HadAnnotated}), \textsf{CMSO}-definability (resp. $\mathsf{CMSO/tw}\!+\!\mathsf{dp}$-definability) ensures Condition (i) of \zcref{thm_reduceToStars} when $\p = \tw$ (resp. $\p = \hw$).
This naturally raises the question of what assumptions may guarantee Condition (ii).
One answer, explored in this section, involves the notion of \textsl{topological minor folios}.

\subsection{Folio representation}

\zcref{thm_reduceToStars} has two algorithmic conditions.
The first one will be derived from the conditions of \zcref{prop_TwAnnotated,prop_HadAnnotated}.
The second will follow from the ability to replace some part $(G,\chi,B,\rho)$ of the input colorful graph with another in which we can route the same set of rooted topological minors (up to a given size, depending on the problem).
Provided that this replacement can be performed algorithmically, it will yield an additional combinatorial condition for our algorithmic meta-theorems.

 \paragraph{Topological folios.}

Especially for boundaried $0$-colorful graphs, we omit the empty coloring and simply write triples $(H,B,\rho),$ which we call \emph{boundaried graphs}.

Given two boundaried graphs $(G,B,\rho)$ and $(H,B',\rho'),$ we say $(H',B',\rho')$ is a \emph{topological minor} of $(G,B,\rho)$ if there exists a pair $(M,T)$ where:
\begin{itemize}
 \item $M$ is a subgraph of $G$ with $B \subseteq T \subseteq V(M),$ and all vertices in $V(M)\setminus T$ have degree 2,
 \item Dissolving all vertices in $V(M)\setminus T$ yields a boundaried graph isomorphic to $(H,B',\rho').$
\end{itemize}

Notice that if $(H,B,\rho)$ is a {topological minor} of $(G,B,\rho),$ and $\chi \colon B\to 2^{[q]}$ then $(H,\chi,B,\rho)$ and $(G,\chi,B,\rho)$ are compatible confined boundaried $q$-colorful graphs.
 
The \emph{$d$-topological-folio} of a confined boundaried $q$-colorful graph $(G,\chi,B,\rho)$ is the set of all confined boundaried $q$-colorful graphs $(H,\chi,B,\rho)$ where $(H,B,\rho)$ is a topological minor of $(G,B,\rho)$ and $|H| \leq d.$ 
We denote this set by $d\mbox{-}\tfolio(G,\chi,B,\rho).$

\begin{theorem}\label{thm_compression}
There exists a computable function $f_{\ref{thm_compression}}\colon\Nbbb\to\Nbbb$ such that for every $d\in\Nbbb$ and every confined boundaried $q$-colorful graph $(G,\chi,B,\rho),$ there exists a (compatible) confined boundaried $q$-colorful graph $(G',\chi,B,\rho)$ where:
$|G'| \leq f_{\ref{thm_compression}}(d,|B|)$ and $d\mbox{-}\tfolio(G,\chi,B,\rho) = d\mbox{-}\tfolio(G',\chi,B,\rho').$
\end{theorem}

\begin{proof}
This follows from \cite[Lemma 2.2]{GroheKMW11Finding} for boundaried graphs.
The confined $q$-colorful case holds because $d\mbox{-}\tfolio(G,\chi,B,\rho)$ depends only on $(G,B,\rho),$ not on $\chi.$
\end{proof}

\zcref{thm_compression} associates to every confined boundaried $q$-colorful graph $(G,\chi,B,\rho)$ a new one with the same $d$-folio but with a bounded number of vertices.
We denote this new confined boundaried $q$-colorful graph by $d\mbox{-}\mathsf{frep}(G,\chi,B,\rho).$ 

\begin{lemma}\label{u7y87uyidb3}
There exists an algorithm that computes $d\mbox{-}\mathsf{frep}(G,\chi,B,\rho)$ in $\Ocal_{d+|B|}(|G|^3)$ time.
\end{lemma}

\begin{proof}
The algorithm enumerates all $\Ocal_{d+|B|}(1)$ confined boundaried $q$-colorful graphs $(H,\chi,B,\rho)$ with $|H| \leq f_{\ref{thm_compression}}(d,|B|)$ and returns one with the same $d$-topological-folio as $(G,\chi,B,\rho).$
By \zcref{thm_compression}, such a graph exists.
The computation of $d\mbox{-}\tfolio(G,\chi,B,\rho)$ can be done in $\Ocal_{d+|B|}(|G|^3)$ time, according to the 
main result of \cite{GroheKMW11Finding}.
\end{proof}

\subsection{AMTs for colorful graphs}

We now have all the necessary ingredients for our AMT's.
We conclude with the following corollary of \zcref{thm_reduceToStars} when applied for $\tw$ and taking into account \zcref{thm_SRTW}, \zcref{prop_TwAnnotated}, along with \zcref{u7y87uyidb3} and \zcref{thm_compression}.

Given a $q \in \Nbbb,$ a functional problem $\funp\colon\Ccal^{q}_{\mathsf{all}}\to\Nbbb$ is \emph{folio-representable} if there exists a function $g\colon\Nbbb\to\Nbbb$ such that $(g(|B|)\mbox{-}\mathsf{frep}, f_{\ref{thm_compression}})$ is a representation of $\funp.$

\begin{theorem}\label{the_MSOL}
Let $q \in \Nbbb$ and $\funp\colon\Ccal^{q}_{\mathsf{all}}\to\Nbbb$ be a folio-representable, \MSO-definable functional problem (via a formula $\phi$ and a $\bullet \in \{\max,\min\}$). There exists an algorithm that, given a $q$-colorful graph $(G,\chi),$ computes $\funp(G,\chi)$ in $\Ocal_{k}(|G|^{3}|\!|G|\!| \log|G|)+\Ocal_{q+|\phi|+k}(|G|^{3})$ time, where $k$ is the maximum monodimensionality of a color in $(G, \chi)$.
\end{theorem}
\begin{proof}

Let $(G, \chi')$ be the $1$-colorful graph such that for each vertex $u \in V(G),$ $\chi'(u) = \{ i \}$ if and only if $\chi(u) \neq \emptyset$ and $i$ is the color of maximum monodimensionality in $(G, \chi).$
Note that the monodimensionality of $i$ in $(G, \chi')$ is of order $\mathcal{O}(k)$ (provided that $q$ is a fixed number).

First, apply \zcref{thm_SRTW} to $(G, \chi')$ to obtain a star decomposition $\Tcal = (T,\beta)$ of $(G,\chi)$ with $\tw$-width $\leq k' \coloneqq f_{\ref{thm_SRTW}}(k)$ in $\Ocal_{k}(|G|^{3}|\!|G|\!| \log|G|)$ time.
Then apply \zcref{thm_reduceToStars}:
Condition (i) holds by \zcref{prop_TwAnnotated} where $f^\mathsf{a}(n)=n$ and $\p=\tw.$
For Condition ii), we have that $(g(|B|)\mbox{-}\mathsf{frep},f_{\ref{thm_compression}})$ is a representation of $\funp$ and the algorithms computing $g(|B|)\mbox{-}\mathsf{frep}$ and $f_{\ref{thm_compression}}$ respectively are given by \zcref{u7y87uyidb3} and \zcref{thm_compression} for $f^\mathsf{r}(n)=n^3.$
\end{proof}

If now in the above proof we use \zcref{thm_SRC} instead of \zcref{thm_SRTW}, we replace $\tw$ by $\hw,$ and we apply \zcref{prop_HadAnnotated} instead of \zcref{prop_TwAnnotated} we have the following.

\begin{theorem}
\label{the_MSOL_tw}
Let $q \in \Nbbb$ and $\funp\colon\Ccal^{q}_{\mathsf{all}}\to\Nbbb$ be a folio-representable, $\mathsf{CMSO/tw}\!+\!\mathsf{dp}$-definable functional problem (via formula $\phi$ and a $\bullet \in \{\max,\min\}$).
There exists an algorithm that, given a $q$-colorful graph $(G,\chi),$ computes $\funp(G,\chi)$ in $\Ocal_{k}(|G|^{3}|\!|G|\!| \log|G|)+\Ocal_{q+|\phi|+k}(|G|^{3})$ time, where $k$ is the maximum Hadwiger number of a color in $(G, \chi)$.
\end{theorem}

Notice that \zcref{the_MSOL} and \zcref{the_MSOL_tw} can be viewed as adaptations of \zcref{prop_TwAnnotated} and \zcref{prop_HadAnnotated} to optimization problems on colorful graphs, where the combinatorial condition captures how the colored vertices of the input should be distributed in $G,$ rather than imposing a sparsity condition on the input graph that ignores the colors.

 \paragraph{Folio representability for decision problems on colorful graphs.}

A decision problem on $q$-colorful graphs can be seen as a subset $\Ccal\subseteq\Ccal^{q}_{\mathsf{all}}.$ For completeness, we derive the notion of folio representability for decision problems from its functional counterpart as follows:
A decision problem $\Ccal\subseteq\Ccal^{q}_{\mathsf{all}}$ is \emph{folio representable} if there exists a function $g\colon\Nbbb\to\Nbbb$ such that for every $q$-boundaried graph $\mathbf{G}$ with boundary $B$ and every compatible confined $q$-boundaried graph $\mathbf{C},$ it holds that
\[
\mathbf{C}\oplus\mathbf{G}\in\Ccal \iff \mathbf{C}'\oplus\mathbf{G}\in\Ccal,
\]
where $\mathbf{C}' = g(|B|)\mbox{-}\mathsf{frep}(\mathbf{C}).$
Recall that in the introduction we presented the decision problem analogues of \zcref{the_MSOL} and \zcref{the_MSOL_tw}, namely \zcref{thm_rwtMeta_Intro} and \zcref{thm_rcMeta_Intro}.
These can be easily derived as follows:
Let $q\in\Nbbb,$ and let $\Ccal\subseteq \Ccal^{q}_{\mathsf{all}}$ be a decision problem on $q$-colorful graphs that is folio-representable and \CMSO-definable (resp. $\mathsf{CMSO/tw}\!+\!\mathsf{dp}$-definable) by some formula $\phi.$
We define a new formula $\phi'$ on $(q+1)$-colorful graphs so that
$(G,\chi+X)\models\phi'$ if and only if the equivalence $(G,\chi)\models \phi \iff X\neq\emptyset$ holds.
Then, $(G,\chi)\models \phi$ if and only if $\p_{\phi'}^{\min}(G,\chi)=1.$
We then apply the algorithm of \zcref{the_MSOL} (resp. \zcref{the_MSOL_tw}) to compute the value of $\p_{\phi'}^{\min}(G,\chi).$
If this value is $1,$ then $(G,\chi)\in \Ccal,$ while if it is $0,$ then $(G,\chi)\notin \Ccal.$

\section{Conclusion}
\label{sec_concl}

As previously noted, the algorithmic contributions of this paper are limited to showcasing the algorithmic aspects of the colorful minor framework we developed.
Below, we outline several promising directions for further research on problems involving colorful graphs.

\subsection{Packing and covering colorful minors}
\label{subsec_packing_covering}

Let $q$ be a non-negative integer and let $(H,\psi)$ and $(G,\chi)$ be $q$-colorful graphs.
A \emph{packing of $(H,\psi)$ in $(G,\chi)$ of size $k$} is a collection of $k$ pairwise vertex-disjoint $q$-colorful subgraphs of $(G,\chi),$ each of which contains $(H,\psi)$ as a colorful minor.
We denote by $\pack_{H,\psi}(G,\chi)$ the largest $k$ for which such a packing exists, and by $\cover_{H,\psi}(G,\chi)$ the smallest integer $k$ for which there is a set $S\subseteq V(G)$ with $|S|\leq k$ such that $(G-S,\chi)$ does not contain $(H,\psi)$ as a colorful minor.
As a general research project, one might consider the problem of computing these two parameters
for different instantiations of $(H,\psi).$
Next, we give examples of known problems that can be expressed by these parameters, given suitable choices of $(H,\psi).$


 \paragraph{Problems expressible via packing and covering colorful minors.}

First of all, assume that $(H, \psi)$ is the $2$-colorful rainbow $K_{1}$ and observe that by Menger's theorem, for every $2$-colorful graph $(G,\chi),$ it holds that $\cover_{H,\psi}(G,\chi) = \pack_{H,\psi}(G,\chi),$ and the two parameters are obviously computable in polynomial time.
This is extended by \zcref{lemma_multicolorlinakge} implying that, for every positive $r,$ 
$\cover_{H_{r},\psi}$ and $\pack_{H_{r},\psi}$
are (linearly) equivalent when we consider the 
$(r+1)$-colorful graph $(H_{r},\psi_{r})$ defined as follows: 
$H_{r}$ is the union of $r$ disjoint copies $H_{1}, \ldots, H_{r}$ 
of $K_{1}$ and, for $i\in[r],$ $\psi_{r}$ assigns to the vertex of $H_{i}$ the colors $\{1, i+1\}.$

Let now $(H,\psi)$ be the $1$-colorful rainbow $K_{2}.$ 
Then $\pack_{H,\psi}$ is the maximum number of vertex disjoint paths between pairs of distinct vertices from $\chi^{-1}(1).$
This is the \textsc{Maximum Disjoint T-Paths} problem that can be solved in $\Ocal(n^4),$ according to \cite{Hananya1994disjoint}.

Notice that for the same $(H,\psi),$ $\cover_{H,\psi}(G,\chi)\leq k$ asks whether there are $\leq k$ vertices that, when
removed, leave all remaining colored vertices in different connected components.
This is the \textsc{Node Multiway Cut with Deletable Terminals} problem that is \NP-complete even when the input graph is restricted to be planar \cite{JohnsonMPPSL2024multiway}.

We now proceed with an example of a problem that can be reduced to a covering problem on $1$-colorful graphs.
%

This is \textsc{Node Multicut with Deletable Terminals} where the input is a graph $G,$ a sequence 
$(s_{1},t_{1}),\ldots,(s_{r},t_{r})$ of pairs of distinct terminals, and a $k\in\Nbbb.$
The question is whether there is a set $S$ of at most $k$ vertices 
such that, for $i\in[r],$ if $s_{i},t_{i}\in V(G-S)$ 
then there is no $s_{i}$-$t_{i}$ path in $G-S.$
We now define a $1$-colorful graph $(G', \chi)$ as follows.
We first subdivide once all edges of $G$ and then 
add a set $X = \{ v_{1}^{1}, \ldots, v_{1}^{k+1}, \ldots, v^{1}_{r}, \ldots, v_{r}^{k+1} \}$ of $(k+1)r$ new vertices and for $i \in [r]$ and $j \in [k+1]$ make the vertex $v_{i}^{j}$ adjacent to the vertices $s_{i}$ and $t_{i}.$ 
We also define $\chi$ so that $\chi^{-1}(1) = X.$
Let $H$ be a cycle on 4 vertices and let $\psi$ be a coloring assigning the color 1 to some of the vertices of $H.$
It is easy to observe that $G,(s_{1},t_{1}),\ldots,(s_{r},t_{r}), k$ is a \yes instance of \textsc{Node Multicut with Deletable Terminals} if and only if $\cover_{H,\psi}(G',\chi) \leq k.$

 \paragraph{Constructing algorithms for packing and covering.}

It is easy to observe that both parameters $\cover_{H,\psi}$ and $\pack_{H,\psi}$ are colorful minor-monotone, for every instantiation of $(H,\psi).$
Therefore, checking whether $\cover_{H,\psi}(G,\chi)\leq k$ (resp. $\pack_{H,\psi}(G,\chi) \geq k$) can be done by an algorithm running in $\Ocal_{q+|H|+k}(|G| |\!|G|\!|^{1+o(1)})$ time because 
of \zcref{thm_NonConstructiveColorfulParameters}.
However, as already commented after \zcref{thm_NonConstructiveColorfulParameters}, this argument is not constructive and only yields the existence of the claimed algorithms. 

One may construct an algorithm for checking whether $\pack_{H,\psi}(G,\chi) \geq k$ as follows: We first check whether $G$ contains a $q$-colorful rainbow $K_{(k+1)|H|}$ using \zcref{thm_introColorfulMinors} in $\Ocal_{q+|H|+k}(|G|^{\Ocal(1)})$ time.
If this is the case, then we report a {positive} answer. If not, 
then, as $\Ccal_{k}=\{(G,\chi)\mid \pack_{H,\psi}(G,\chi) < k\}$ is folio representable and $\mathsf{CMSO/tw}+\mathsf{dp}$-definable,
we can apply \zcref{thm_rcMeta_Intro} and construct an algorithm running deciding whether $(G,\chi)\in\Ccal_{k}$ in time 
 $\Ocal_{q +|H|+k}(|G|^{\Ocal(1)}).$


On the other hand, a general constructive algorithm for checking whether $\cover_{H,\psi}(G,\chi) \leq k$ requires more work.
We may certainly begin by checking whether the $r$-colorful rainbow $K_{|H| + k}$ is a colorful minor of $(G,\chi),$ and if so, we report negatively.
Again, may attempt to apply \zcref{thm_rcMeta_Intro} for
$\Ccal_{k} = {(G,\chi) \mid \cover_{H,\psi}(G,\chi) \leq k}.$
While $\Ccal_{k}$ is $\mathsf{CMSO/tw} + \mathsf{dp}$-definable, it is not folio representable, and therefore we cannot apply the replacements used in the proof of \zcref{thm_reduceToStars}.

However, checking whether $\cover_{H,\psi}(G,\chi) \leq k$ is still possible when $q = 1.$
For this, we may apply the replacement argument from \zcref{thm_reduceToStars} by using instead the stronger notion of a \emph{$(d,k)$-folio} of a boundaried graph $\mathbf{G} = (G, B, \rho).$
We denote by $(d,k)$-$\folio(\mathbf{G})$ the collection of all $d$-folios of the boundaried graphs derived from $\mathbf{G}$ after removing $\leq k$ vertices.
As shown by Fomin, Lokshtanov, Panolan, Saurabh, and Zehavi in \cite{FedorLPSZ2020Hitting}, one can construct an algorithm that, given $\mathbf{G} = (G, B, \rho),$ outputs, in $\Ocal_{d+k+|B|}(|G|^{4})$ time, a new boundaried graph $\mathbf{G}'$ such that $(d,k)$-$\folio(\mathbf{G}) = (d,k)$-$\folio(\mathbf{G}'),$ and the graph of $\mathbf{G}'$ excludes a clique as a minor, where the size of the excluded clique depends on a function of $d,$ $k,$ and $|B|.$
(This follows by repeatedly applying the algorithm from \cite[Theorem 6]{FedorLPSZ2020Hitting}, which identifies irrelevant vertices.)

This enables a replacement argument that produces an equivalent instance whose graph has bounded Hadwiger number.
The problem can then be solved in $\Ocal_{|H|+k}(|G|^{\Ocal(1)})$ time, again using the decision version of \zcref{prop_HadAnnotated}.

To extend this argument to $q > 1,$ one would need to prove a “colorful” analogue of \cite[Theorem 6]{FedorLPSZ2020Hitting}, which would be able to identify a vertex that is irrelevant for every possible hitting set of size $\leq k.$
While one might reasonably believe such an extension is possible, the details of its construction appear to be quite intricate.

All algorithms discussed in this section have significant complexity, both in their parameter dependence and in their polynomial factors.
Optimizing these algorithms -- either in general or for specific, interesting instantiations of $(H, \psi)$ -- could be a promising direction for future research.

\subsection{Extending graph parameters to colorful ones}\label{sec_extend_to_colorful}

Recently, two main mechanisms have been proposed to extend treewidth to annotated graphs (i.e., $1$-colorful graphs). These are the \textsl{rainbow extension} \cite{ThilikosW2025graphminorsstructure} and the \textsl{torso extension} \cite{JansenS2024SteinerTree,HodorLMR24quick}. We propose here an extension of both mechanisms to $q$-colorful graphs that applies to every minor-monotone parameter.

Let $q\ge 0$ and let $\p\colon \gall\to\Nbbb$ be a minor-monotone graph parameter. We define the \emph{torso $q$-extension} of $\p$ and the \emph{rainbow $q$-extension} of $\p$ as the $q$-colorful graph parameters
$\mathsf{t}\p:\Ccal^{q}_{\mathsf{all}}\to\Nbbb$ and
$\mathsf{a}\p:\Ccal^{q}_{\mathsf{all}}\to\Nbbb$,
where
\begin{eqnarray*}
\!\!\!\!\!\!q\mbox{-}\mathsf{t}\p(G,\chi) \!\!\!\!& = &\!\!\!\! \min \{\p(\torso(G,X))\mid \mbox{$X\subseteq V(G)$ and every component of $G-X$ is $q$-restricted} \}.\\
\!\!\!\!\!\!q\mbox{-}\mathsf{r}\p(G,\chi) \!\!\!\!& = &\!\!\!\! \max \{ \p(H)\mid \mbox{the $q$-colorful rainbow $H$ is a colorful minor of $(G,\chi)$} \}.
\end{eqnarray*}

As mentioned in the introduction, $1\mbox{-}\mathsf{t}\tw$, the \emph{torso treewidth} of an annotated graph, has already been proposed—using different terminology—by Jansen and Swennenhuis in \cite{JansenS2024SteinerTree} and by Hodor, La, Micek, and Rambaud in \cite{HodorLMR24quick}.
Also, it is the main combinatorial engine in the proof of \zcref{the_MSOL}.

On the other hand, $1\mbox{-}\mathsf{r}\tw$, the \emph{rainbow treewidth} of an annotated graph, is equivalent to the bidimensionality $\bidim$, which is $1\mbox{-}\mathsf{r}\bg$, where $\bg(G)$ denotes the maximum $k$ for which $G$ contains the $(k\times k)$-grid as a minor.
This leads us to view $q\mbox{-}\mathsf{r}\bg$ (which, again, is equivalent to $q\mbox{-}\mathsf{r}\tw$) as the correct notion of bidimensionality for colorful graphs.
In this sense, $\bidim$ measures the ``simultaneous presence’’ of all colors across the bidimensional surface of a grid.

\medskip
Recall that, as already proved in \cite{HodorLMR24quick}, $1\mbox{-}\mathsf{t}\tw(G)$ is equivalent to the maximum $k$ for which $(G,X)$ contains the $(1,k)$-segregated grid as a colorful minor.
Moreover, by \zcref{thm_restrictiveTreewidthIntro}, this equivalence extends to colorful graphs with any number $q\ge 1$ of colors.
That is, $q\mbox{-}\mathsf{t}\tw(G)$ is equivalent to the maximum $k$ for which $(G,\chi)$ contains one of the $\lceil \nicefrac{1}{2}(q-1)!\rceil$ parametrically-distinct $(q,k)$-segregated grids as a colorful minor.
This motivates us to interpret this quantity as the \emph{monodimensionality} of a colorful graph, as it measures the ``segregated presence’’ of all colors along the monodimensional boundary of a grid. 

A natural question is when the two aforementioned mechanisms for generating colorful extensions yield equivalent parameters.
Let $\p \colon \gall \to \Nbbb$ be a minor-monotone graph parameter.
It is easy to see that for every $q$-colorful graph $(G,\chi)$, we have
$q\mbox{-}\mathsf{r}\p(G,\chi)\le q\mbox{-}\mathsf{t}\p(G,\chi)$, regardless of the choice of $\p$ or $q$.
However, it is not generally the case that there exists a function $f:\Nbbb\to\Nbbb$ such that
$q\mbox{-}\mathsf{t}\p(G,\chi)\le f(q\mbox{-}\mathsf{r}\p(G,\chi)),$ for all $(G,\chi)$. Thus, the question becomes: For which minor-monotone parameters $\p$ does it hold that $q\mbox{-}\mathsf{r}\p$ and $q\mbox{-}\mathsf{t}\p$ are equivalent?

\subsection{From faces to vortices}

In this paper we presented some algorithmic consequences of two of our three structural theorems, i.e., \zcref{thm_intro_ExcludeRainbowClique}
and \zcref{thm_restrictiveTreewidthIntro}.
Clearly, the structure revealed by \zcref{thm_ExcludingRainbowGridIntro} is “intermediate” with respect to that of \zcref{thm_intro_ExcludeRainbowClique} and \zcref{thm_restrictiveTreewidthIntro}.
It provides a structural decomposition of $q$-colorful graphs
of small bidimensionality, i.e., $q$-colorful graphs excluding a large $q$-colorful rainbow grid as a minor, which we believe should be useful in algorithmic design.
As we already mentioned, the algorithmic consequences of \zcref{thm_ExcludingRainbowGridIntro} go beyond the objectives of this paper.
However, we outline below some directions on the \textsc{Steiner Tree} problem that are inspired by the current research.

In the introduction, we commented on \textsc{Steiner Tree} as one of the most studied problems on annotated graphs.
Using the terminology of this paper, Jansen and Swennenhuis in \cite{JansenS2024SteinerTree} proved that \textsc{Steiner Tree}, with input a $1$-colorful graph $(G,\chi)$ can be solved in $2^{\Ocal(k)} \cdot |G|^{\Ocal(1)}$ time where $k$ is the torso treewidth of $(G,\chi).$
Is torso treewidth the \say{best possible} structural parameter on $(G,\chi)$ that, when bounded, \textsc{Steiner Tree} can be solved in polynomial time?
The answer is negative.
This follows by the fact that a polynomial algorithm has been recently given by Groenland, Nederlof, and Koana in \cite{GroendlandNK2024Polynomial} for the case where $(G, \chi)$ excludes the $1$-colorful rainbow $K_{4}$ as a colorful minor.
The results of \cite{GroendlandNK2024Polynomial} apply to all colorful minors of the $(1,k)$-segregated grid on which torso treewidth may be arbitrarily large due to \zcref{lemma_segregatedLowerbound}.
We believe that the \say{limit} is given by the combinatorial bound of \zcref{thm_ExcludingRainbowGridIntro}.
We conjecture that the bidimensionality of $X=\chi^{-1}(1)$ in $G$ is the structural frontier for the complexity of \textsc{Steiner Tree} problem in classes that are colorful minor-closed.
In particular we conjecture the following complexity dichotomy.

\begin{conjecture}
Unless $\mathsf{P}\neq\NP,$ \textsc{Steiner Tree} can be solved in polynomial time in a {$1$-}colorful minor-closed class $\Ccal$ if and only if $\Ccal$ has bounded bidimensionality.
\end{conjecture}

Note that the negative part of the above conjecture already follows from the fact that the \textsc{Steiner Tree} problem is \NP-complete in planar graphs \cite{GareyJS76}.
We believe that the local structure provided by \zcref{thm_ExcludingRainbowGridIntro} should serve as a starting point for proving the positive part as well.
In particular, we now have that all terminals in $X = \chi^{-1}(1)$ are confined to the apices and the “few” vortices of the near-embedding of $G.$
Motivated by the fact that \textsc{Steiner Tree} can be solved in polynomial time when all terminals lie on “few” faces of a planar graph \cite{DreyfusW1971Steiner, EricksonMV1987SendAndSplit, KisfaludiBakNvL2020NearlyETH}, the key challenge becomes how to generalize such face-based techniques to handle apices and vortices in near-embeddings of the kind provided by \zcref{thm_ExcludingRainbowGridIntro}.
We further anticipate that this methodological approach may lead to algorithmic dichotomies -- such as the one conjectured above -- for a broader class of problems on annotated graphs, where tractability arises from restricting terminals to few faces of the input graph.

\newpage

\bibliographystyle{plainurl}
\bibliography{literature_colorful}

\newpage
\appendix

\section{Some problems where our algorithmic results apply}\label{sec_AppendixA}

In this section, we present a selection of natural optimization problems for which our algorithmic results apply. These problems are definable in $\CMSO$ and satisfy the conditions required by \zcref{the_MSOL} or \zcref{the_MSOL_tw}. 

\subsection{Problems where \zcref{the_MSOL} can be applied}
We now give some examples of optimization problems where \zcref{the_MSOL} is applicable. 
For each of them, there is some formula $\phi$ in $\CMSO$ and a choice of $\bullet\in\{\min,\max\}$ such that if the input is the $q$-colorful graph $(G,\chi),$ then answer is $\p_{\phi}^{\bullet}(G,\chi).$
For all of them it is straightforward to verify that there exists a function $g:\Nbbb\to\Nbbb$ such that $(g(|B|)\mbox{-}\mathsf{frep}, f_{\ref{thm_compression}})$ is a representation of $\funp.$
In particular then the question is about paths or cycles, then $g(|B|)=|B|$ is enough, while when the problem is defined in terms of some graph $H,$ then we may pick $g(b)=|B|+|H|.$

We proceed with the list of problems.
For each of them we may associate a decision problem by introducing an integer $k$ and asking whether the requested maximum/minimum number is at least/at most $k.$
In each case, we give arguments why the corresponding decision version is \NP-hard.

 \paragraph{}
\fbox{\begin{minipage}{15cm}
\noindent\textsc{Spanning $H$-Topological Minor}

\noindent\textsl{Input:} A $1$-colorful graph $(G,\chi).$\\
\noindent\textsl{Question:} What is the maximum number of vertices of $\chi^{-1}(1)$ a subdivision of $H$ in $G$ may span?
\end{minipage}
}
\medskip

Notice that when $H=K_{3},$ then the above problem yields \textsc{Maximum Spanning Cycle} that asks for the maximum number of vertices of $\chi^{-1}(1)$ that can be spanned by a cycle of $G.$ The corresponding decision problem is \NP-hard as, if we additionally set $\chi^{-1}(1)=V(G),$ we obtain the \textsc{Hamiltonian Cycle} problem.

 \paragraph{}
\fbox{\begin{minipage}{15cm}
\noindent\textsc{Longest Spanning Path/Cycle}

\noindent\textsl{Input:} A $1$-colorful graph $(G,\chi).$\\
\noindent\textsl{Question:} What is the maximum path/cycle containing all vertices of $\chi^{-1}(1)$?
\end{minipage}
}
\medskip

The corresponding decision problem is \NP-hard 
as, if we additionally set $\chi^{-1}(1)=V(G),$ we obtain the \textsc{Longest Cycle}/\textsc{Longest Path} problem.

 \paragraph{}

\fbox{\begin{minipage}{15cm}
\noindent\textsc{Minimum/Maximum Number of Spanning Cycles}

\noindent\textsl{Input:} A $1$-colorful graph $(G,\chi).$\\
\noindent\textsl{Question:} What is the {minimum}/{maximum} number of pairwise-disjoint cycles 
spanning $\chi^{-1}(1)$?
\end{minipage}
}
\medskip

The corresponding decision problem of the \textsc{Minimum} variant is \NP-hard because it gives \textsc{Hamiltonian Cycle} when $\chi^{-1}(1)=V(G)$ and $k =1.$
The corresponding decision problem of the \textsc{maximum} variant 
reduces to the \textsc{Disjoint Paths} problem, in which we are given a graph $G$ along 
with a list $(s_{1},t_{1}),\ldots,(s_{r},t_{r})$ of pairs of terminals and asks 
if $G$ contains $r$ pairwise vertex disjoint paths joining them. 
The reduction introduces a new vertex set $X=\{x_{1},\ldots,x_{r}\}$ to $G$ and makes each $x_{i}\in X$
adjacent to $s_{i}$ and $t_{i}.$

 \paragraph{}
\fbox{\begin{minipage}{15cm}
\noindent\textsc{Nailed Packing of $H$-Subdivisions}

\noindent\textsl{Input:} A $1$-colorful graph $(G,\chi).$\\
\noindent\textsl{Question:} what is the maximum 
 $k$ for which $G$ contains the disjoint union of 
 $k$ subdivisions of $H$ as a subgraph, each of which spans at least one vertex of $\chi^{-1}(1)$?
 \end{minipage}
}
\medskip

The corresponding decision problem, for $H=K_{3},$ is \NP-hard because when $\chi^{-1}(1)=V(G),$ this gives the \textsc{Cycle Packing} problem.
This problem may be extended to $q$-colorful graphs where we may add that each cycle spans all $q$ colors.

 \paragraph{}

\fbox{\begin{minipage}{15cm}
\noindent\textsc{Alternating Spanning Cycle/Path}

\noindent\textsl{Input:} A $2$-colorful graph $(G,\chi).$

\noindent\textsl{Question:} What is the maximum length cycle/path spanning $\chi^{-1}(1)\cup \chi^{-1}(2)$ such that no two vertices of the same set in $\{\chi^{-1}(1),\chi^{-1}(2)\}$ appear consecutively along this path/cycle.
 \end{minipage}
}
\medskip

The corresponding decision problem is \NP-hard because 
of the following reduction from \textsc{Hamiltonian Cycle/Path}:
Given an instance $G$ of \textsc{Hamiltonian Cycle/Path} consider the instance $G',\chi^{-1}(1),\chi^{-1}(2),|V(G')|$
of the decision version of \textsc{Alternating Spanning Cycle/Path} where 
$G'$ is the bipartite double cover of $G$ (that is the tensor product $G\times K_{2}$ of $G$ and $K_{2}$) and $\chi^{-1}(1)$ and $\chi^{-1}(2)$ are its two color classes.

We remark that the above problem can be extended for general $q$-colorful graphs and for every pattern of color alternation.

 \paragraph{}

\fbox{\begin{minipage}{15cm}
\noindent\textsc{Minimum Steiner $q$-Forest}

\noindent\textsl{Input:} A $q$-colorful graph $(G,\chi).$

\noindent\textsl{Question:} What is the minimum size of a subgraph 
$H$ of $G$ in
which any two vertices belonging to the same set $\chi^{-1}(i)$ are connected by a path in $H$?
 \end{minipage}
}
\medskip

The corresponding decision problem, for the case $q=2,$ is $\NP$-hard, as proved in \cite{HofPW2009Partitioning}.
\medskip

\subsection{Problems where \zcref{the_MSOL_tw} can be applied} We now give two examples of optimization problems where \zcref{the_MSOL_tw} is applicable.

 \paragraph{Unordered Linkability problems.}
Given a graph $G,$ a set $X\subseteq V(G),$ and a graph $H$ on $h$ vertices, we say that $X$ is \emph{$H$-linkable}
in $G$ if for every subset $S\in \binom{X}{h}$ there is a $\Pcal\in \binom{S}{2}$
and a collection of internally vertex-disjoint 
paths in $G$ between the pairs in $\Pcal$ that, when contracted to single edges, give a graph that is isomorphic to $H.$ We consider the following optimization problem.

 \paragraph{}
\fbox{
\begin{minipage}{15cm}
\noindent \textsc{Unordered Linkability}\\
\noindent{\sl Input}: A $1$-colorful graph $(G,\chi)$ and a graph $H$ on $h$ vertices.\\
\noindent{\sl Question}: what is the biggest $X'\subseteq \chi^{-1}(1)$ so that $X'$ is $H$-linkable in $G$?
\end{minipage}}
\medskip

Notice that the decision version of \textsc{Unordered Linkability} where $H$ is a cycle on $h$ vertices expresses (by asking whether $k=|\chi^{-1}(1)|$) is the \textsc{Cyclability} problem asking, given $(G,\chi),$ whether 
$\chi^{-1}(1)$ is \emph{$h$-cyclable} in $G,$ that is every $h$ vertices of $\chi^{-1}(1)$ belong in some cycle of $G$ (see \cite{Dirac1960vorhandene, WatkinsM1967Cycles, GyoriP2001claw, AldredBHK1999cubic, FlandrinLM2007cycles} for the combinatorial properties of $h$-cyclable sets).
In \cite{GolovachKTM2017cyclability} it was proved that \textsc{Cyclability}, parameterized by $h,$ is fixed parameter tractable (\FPT) for planar graphs, while \textsf{co}-\textsf{W}[1]-hard in general.
By \FPT we mean that it admits an \FPT-algorithm, i.e., an algorithm running in $\Ocal_{k}(|G|^{\Ocal(1)})$ time.

 \paragraph{Ordered Linkability problems.}

We next define the \say{ordered} version of linkability. Let $(G,\chi)$ be a $1$-colorful graph, and 
let $H$ be a graph with $V(H)=\{x_{1},\ldots,x_{h}\}.$ We say that $\chi^{-1}(1)$ is 
 \emph{$H$-linked} in $G$ if, \textsl{for every sequence} $v_{1},\ldots,v_{h}$
of vertices in $\chi^{-1}(1)$ there is a $\Pcal\in \binom{\{v_{1}, \ldots, v_{h}\}}{2}$
and a collection of internally vertex-disjoint 
paths in $G$ between the pairs in $\Pcal$ that, when contracted to single edges, give a graph 
with vertex set $S$ that is isomorphic to $H$ via the isomorphism that maps $x_{i}$ to $v_{i},$ $i\in[h].$ 
The combinatorics of $H$-linked sets
has been studied in \cite{GouldW2007Subdivision,Kostochka2005extremal,EllinghamPDG2012Linkage,GouldKG2006linked,FerraraGTW2006linked}. In order to study the algorithmic properties of $H$-linked sets, we consider the following optimization problem.

 \paragraph{}
\fbox{
\begin{minipage}{15cm}
\noindent \textsc{Ordered Linkability}\\
\noindent{\sl Input}: A $1$-colorful graph $(G,\chi)$ and a graph $H$ on $h$ vertices.\\
\noindent{\sl Question}: what is the biggest $X'\subseteq \chi^{-1}(1)$ so that $X'$ $H$-linked in $G$?
\end{minipage}
}
\medskip

In \cite{GolovachST23Model}, 
a special case of \textsc{Ordered Linkability} problem was examined, namely the \textsc{Disjoint Paths Linkability} problem that is obtained if $H$ is the disjoint union of $h/2$ edges.
It was proven in \cite{GolovachST23Model} that, when parameterized by $h,$ 
this problem is not \FPT, unless $\mathsf{FPT}=\mathsf{W[1]}.$
$H$-linked sets for this particular $H$ where introduced by \cite{Watkins1968arcs} in the late 60s and its graph-theoretical properties 
have been extensively studied in \cite{LarmanΜ1970configurations,Jung1970Verallgemeinerung,RobertsonS1995Graph,BollobasT1996Highly,EgawaFG2000specified,KawarabayashiKG2006linked,ThomasW2005linkages}.

\medskip

As observed in \cite{GolovachST23Model}, the decision versions of \textsc{Unordered Linkability} and \textsc{Ordered Linkability} are $\mathsf{CMSO/tw}\!+\!\mathsf{dp}$-definable, therefore, from \zcref{prop_HadAnnotated}, they both admit an $\Ocal_{h+\hw(G)}(|G|^2)$ time algorithm.
In other words, 
while no \FPT-algorithms are expected when parameterized by $h,$ they both become fixed parameter tractable (\FPT), when they are further restricted to graphs of bounded Hadwiger number.

According to our results, towards achieving \FPT-algorithms, we can do more than bounding the Hadwiger number of the input graph.
For this, it can be easily seen that both these problems are folio-representable, therefore, from \zcref{thm_rcMeta_Intro}, we may deduce that they both admit an $\Ocal_{h+\rh(G,\chi)}(|G|^{\Ocal(1)})$ time algorithm.
This reveals that it is not only the structure of $G$ that determines the parameterized complexity of these problems but the way the vertices in $\chi^{-1}(1)$ are distributed in $G.$

 \paragraph{Packing and covering.}
As commented in the \zcref{sec_concl},
by using \zcref{the_MSOL_tw} we can
construct algorithms running in $\Ocal_{r+|H|+k}(|G|^{\Ocal(1)})$ time for the following decision problem (for every $q\in \Nbbb$).

 \paragraph{}\fbox{\begin{minipage}{15cm}
\noindent\textsc{Packing Minor Models of $(H,\chi)$}

\noindent\textsl{Input:} A $q$-colorful graph $(G,\chi)$ and a $k\in\Nbbb.$\\
\noindent\textsl{Question:} 
Does $(G,\chi)$ contain $k$ subgraphs each containing $(H,\chi)$ as a colorful minor? 
\end{minipage}
}
\medskip

Also, for $q=1,$ as commented in \zcref{subsec_packing_covering}, we have an algorithm running in $\Ocal_{|H|+k}(|G|^{\Ocal(1)})$ time for the following decision problem.

 \paragraph{}\fbox{\begin{minipage}{15cm}
\noindent\textsc{Covering Minor Models of 1-colorful $(H,\chi)$}

\noindent\textsl{Input:} A $q$-colorful graph $(G,\chi)$ and a $k\in\Nbbb.$\\
\noindent\textsl{Question:} 
Does $(G,\chi)$ contain a set $S\subseteq V(G)$ such that $|S|\leq k$ and 
$(G-S,\chi)$ does not contain $(H,\chi)$ as a colorful minor? 
\end{minipage}
}
\medskip

We would like to conclude this section by commenting on how far the simple
arguments we explained in \zcref{subsec_packing_covering} can go towards creating 
constructive algorithms for generalizations of \textsc{Covering Minor Models of 1-colorful $(H,\chi)$}.

Let $\p:\Ccal^{1}_{\textsf{all}}\to\Nbbb$
be any computable colorful minor-monotone 1-colorful graph parameter such that there is a constructive function $f:\Nbbb\to\Nbbb$ where, for every 
1-colorful graph $(G,\chi),$ it holds 
that $\bidim(G,\chi^{-1}(1))\leq f(\p(G,\chi^{-1}(1))).$
Consider now the extension of \textsc{Covering Minor Models of 1-colorful $(H,\chi)$} so that, instead of asking that $|S|\leq k,$ we ask whether 
$\p(G,S)\leq k.$
By using the same arguments as we explained in \zcref{subsec_packing_covering},
one may construct an \FPT-algorithm for this problem.

\medskip
We emphasize that the collection of problems presented in this section is intended to be illustrative rather than exhaustive. Our aim is to showcase the diversity of settings -- spanning connectivity, packing, spanning, and linkability -- where our framework applies.

\end{document}

%% file: formatting.tex
\usepackage[letterpaper, top=25.4mm, bottom=25.4mm, left=25.4mm, right=25.4mm, includefoot]{geometry}

\usepackage[utf8]{inputenc}
\usepackage[T1]{fontenc}
\usepackage{lmodern}
\usepackage{alphabeta}

\usepackage{titling}
\usepackage{xspace}
\usepackage{soul} 
\usepackage{ifthen}

\usepackage{ifthen}

\usepackage{tikz}
\usepackage{xparse}

\newcommand{\colorname}[1]{%
  \ifnum#1=1 red\else%
    \ifnum#1=2 blue\else%
      \ifnum#1=3 green\else%
        \ifnum#1=4 magenta\else%
          \ifnum#1=5 orange\else%
            black%
          \fi%
        \fi%
      \fi%
    \fi%
  \fi%
}

\usepackage[inline]{enumitem}

\usepackage{dsfont}
\usepackage{setspace}

\usepackage{bm}
\usepackage{microtype}
\usepackage{amsmath}
\usepackage{amssymb}
\usepackage{amsfonts}
\usepackage{mathtools}
\usepackage[mathscr]{euscript}

\usepackage[hyphens]{url}
\usepackage{amsthm}

\usepackage{tikz}
\usepackage{xcolor}
\usepackage{subcaption}
\usepackage{graphicx}
\usepackage{wrapfig}

\usepackage{hyperref}
\usepackage{zref-clever}

\usepackage{nicefrac}
\usepackage[textwidth=2.2cm]{todonotes}

\colorlet{myGreen}{green!50!black}
\colorlet{myLightgreen}{green}
\colorlet{myRed}{red!90!black}
\definecolor{myBlue}{rgb}{0.25, 0.0, 1.0}
\definecolor{myLightBlue}{rgb}{0.39, 0.58, 0.93}
\colorlet{myViolet}{myBlue!55!myRed}
\definecolor{myOrange}{rgb}{1.0, 0.66, 0.07}

\definecolor{CornflowerBlue}{rgb}{0.39, 0.58, 0.93}
\definecolor{DarkGoldenrod}{rgb}{0.72, 0.53, 0.04}
\definecolor{BritishRacingGreen}{rgb}{0.0, 0.26, 0.15}
\definecolor{DarkMagenta}{rgb}{0.55, 0.0, 0.55}
\definecolor{AO}{rgb}{0.0, 0.5, 0.0}
\definecolor{BostonUniversityRed}{rgb}{0.8, 0.0, 0.0}
\definecolor{myRed}{rgb}{0.8, 0.0, 0.0}
\definecolor{DarkMidnightBlue}{rgb}{0.0, 0.2, 0.4}
\definecolor{DarkTangerine}{rgb}{1.0, 0.66, 0.07}
\definecolor{AppleGreen}{rgb}{0.55, 0.71, 0.0}
\definecolor{BrightUbe}{rgb}{0.82, 0.62, 0.91}
\definecolor{Amethyst}{rgb}{0.6, 0.4, 0.8}
\definecolor{DarkGray}{rgb}{0.52, 0.52, 0.51}
\definecolor{Gray}{rgb}{0.66, 0.66, 0.66}
\definecolor{BananaYellow}{rgb}{1.0, 0.88, 0.21}
\definecolor{Amber}{rgb}{1.0, 0.75, 0.0}
\definecolor{LightGray}{rgb}{0.83, 0.83, 0.83}
\definecolor{PrincetonOrange}{rgb}{1.0, 0.56, 0.0}
\definecolor{DeepCarrotOrange}{rgb}{0.91, 0.41, 0.17}
\definecolor{CarrotOrange}{rgb}{0.93, 0.57, 0.13}
\definecolor{MidnightBlue}{rgb}{0.1, 0.1, 0.44}
\definecolor{Magenta}{rgb}{0.50, 0.0, 0.50}
\definecolor{BrightPink}{rgb}{1.0, 0.0, 0.5}
\definecolor{BrilliantRose}{rgb}{1.0, 0.33, 0.64}
\definecolor{ChromeYellow}{rgb}{1.0, 0.65, 0.0}
\definecolor{HotMagenta}{rgb}{1.0, 0.11, 0.81}

\definecolor{DarkTangerine}{rgb}{1.0, 0.66, 0.07}
\definecolor{darkyellow}{rgb}{.7, .6, 0.0}
\definecolor{CornflowerBlue}{rgb}{0.39, 0.58, 0.93}
\definecolor{DarkGoldenrod}{rgb}{0.72, 0.53, 0.04}
\definecolor{BritishRacingGreen}{rgb}{0.0, 0.26, 0.15}
\definecolor{AO}{rgb}{0.0, 0.5, 0.0}
\definecolor{MidnightBlack}{rgb}{0.1,0.1,.34}
\definecolor{MidnightBlue}{rgb}{0.1,0.1,0.43}
\definecolor{Black}{rgb}{0,0, 0}
\definecolor{Blue}{rgb}{0, 0 ,1}
\definecolor{Red}{rgb}{1, 0 ,0}
\definecolor{White}{rgb}{1, 1, 1}
\definecolor{DeepMagenta}{rgb}{0.8, 0.0, 0.8}
\definecolor{grey}{rgb}{.6, .6, .6}
\definecolor{darkgrey}{rgb}{.33, .33, .33}
\definecolor{Mygreen}{rgb}{.0, .7, .0}
\definecolor{Yellow}{rgb}{.55,.55,0}
\definecolor{Mustard}{rgb}{1.0, 0.86, 0.35}
\definecolor{applegreen}{rgb}{0.55, 0.71, 0.0}
\definecolor{darkturquoise}{rgb}{0.0, 0.81, 0.82}
\definecolor{celestialblue}{rgb}{0.29, 0.59, 0.82}
\definecolor{green_yellow}{rgb}{0.68, 1.0, 0.18}
\definecolor{crimsonglory}{rgb}{0.75, 0.0, 0.2}
\definecolor{darkmagenta}{rgb}{0.30, 0.0, 0.30}
\definecolor{magenta}{rgb}{0.50, 0.0, 0.50}
\definecolor{internationalorange}{rgb}{1.0, 0.31, 0.0}
\definecolor{darkorange}{rgb}{1.0, 0.55, 0.0}
\definecolor{ao}{rgb}{0.0, 0.5, 0.0}
\definecolor{awesome}{rgb}{1.0, 0.13, 0.32}
\definecolor{darkcyan}{rgb}{0.0, 0.50, 0.50}
\definecolor{violet}{rgb}{0.93, 0.51, 0.93}
\definecolor{brown}{rgb}{0.65, 0.16, 0.16}
\definecolor{orange}{rgb}{1.0, 0.65, 0.0}
\definecolor{DarkGreen}{rgb}{0,.5,0}
\definecolor{BostonUniversityRed}{rgb}{0.8, 0.0, 0.0}

\newcommand{\midnightblue}[1]{{\color{MidnightBlue}#1}}

\vfuzz \hfuzz
\setlist[itemize]{topsep=0pt,partopsep=0pt,itemsep=0pt,parsep=0pt}
\setlist[itemize,1]{label={\small\textbullet}}
\setlist[itemize,2]{label={\tiny\textbullet}}
\setlist[itemize,3]{label=$\cdot$}
\setlist[enumerate]{topsep=0pt,partopsep=0pt,itemsep=0pt,parsep=0pt}
\setlist[enumerate,1]{label=\roman*)}
\setlist[enumerate,2]{label=\alph*)}
\setlist[enumerate,3]{label=\arabic*)}

\hypersetup{
colorlinks=true,
linkcolor=AO!65!black,
citecolor=AO!65!black,
urlcolor=AppleGreen!65!black,
bookmarksopen=true,
bookmarksnumbered,
bookmarksopenlevel=2,
bookmarksdepth=3
}

\theoremstyle{definition}

\newtheorem{environment}{Environment}[section]

\newtheorem{lemma}[environment]{Lemma}
\AddToHook{env/lemma/begin}{\zcsetup{countertype={environment=lemma}}}
\zcRefTypeSetup{lemma}{
Name-sg = Lemma ,
name-sg = Lemma ,
Name-pl = Lemmas ,
name-pl = Lemmas ,
}

\newtheorem*{lemma*}{Lemma}
\AddToHook{env/lemma*/begin}{\zcsetup{countertype={environment=lemma*}}}
\zcRefTypeSetup{lemma*}{
Name-sg = Lemma ,
name-sg = Lemma ,
Name-pl = Lemmas ,
name-pl = Lemmas ,
}

\newtheorem{proposition}[environment]{Proposition}
\AddToHook{env/proposition/begin}{\zcsetup{countertype={environment=proposition}}}
\zcRefTypeSetup{proposition}{
Name-sg = Proposition ,
name-sg = Proposition ,
Name-pl = Propositions ,
name-pl = Propositions ,
}

\newtheorem{corollary}[environment]{Corollary}
\AddToHook{env/corollary/begin}{\zcsetup{countertype={environment=corollary}}}
\zcRefTypeSetup{corollary}{
Name-sg = Corollary ,
name-sg = Corollary ,
Name-pl = Corollaries ,
name-pl = Corollaries ,
}

\newtheorem{theorem}[environment]{Theorem}
\AddToHook{env/theorem/begin}{\zcsetup{countertype={environment=theorem}}}
\zcRefTypeSetup{theorem}{
Name-sg = Theorem ,
name-sg = Theorem ,
Name-pl = Theorems ,
name-pl = Theorems ,
}

\newtheorem*{theorem*}{Theorem}
\AddToHook{env/theorem*/begin}{\zcsetup{countertype={environment=theorem*}}}
\zcRefTypeSetup{theorem*}{
Name-sg = Theorem ,
name-sg = Theorem ,
Name-pl = Theorems ,
name-pl = Theorems ,
}

\newtheorem{conjecture}[environment]{Conjecture}
\AddToHook{env/conjecture/begin}{\zcsetup{countertype={environment=conjecture}}}
\zcRefTypeSetup{conjecture}{
Name-sg = Conjecture ,
name-sg = Conjecture ,
Name-pl = Conjectures ,
name-pl = Conjectures ,
}

\newtheorem*{hypothesis*}{Hypothesis}
\AddToHook{env/hypothesis*/begin}{\zcsetup{countertype={environment=hypothesis*}}}
\zcRefTypeSetup{hypothesis*}{
Name-sg = Hypothesis ,
name-sg = Hypothesis ,
Name-pl = Hypotheses ,
name-pl = Hypotheses ,
}

\newtheorem{observation}[environment]{Observation}
\AddToHook{env/observation/begin}{\zcsetup{countertype={environment=observation}}}
\zcRefTypeSetup{observation}{
Name-sg = Observation ,
name-sg = Observation ,
Name-pl = Observations ,
name-pl = Observations ,
}

\AddToHook{env/example/begin}{\zcsetup{countertype={environment=example}}}
\zcRefTypeSetup{example}{
Name-sg = Example ,
name-sg = Example ,
Name-pl = Examples ,
name-pl = Examples ,
}

\AddToHook{env/remark/begin}{\zcsetup{countertype={environment=remark}}}
\zcRefTypeSetup{remark}{
Name-sg = Remark ,
name-sg = Remark ,
Name-pl = Remarks ,
name-pl = Remarks ,
}

\zcRefTypeSetup{equation}{
Name-sg = Equation ,
name-sg = Equation ,
Name-pl = Equations ,
name-pl = Equations ,
}

\zcRefTypeSetup{chapter}{
Name-sg = Chapter ,
name-sg = Chapter ,
Name-pl = Chapters ,
name-pl = Chapters ,
}

\zcRefTypeSetup{section}{
Name-sg = Section ,
name-sg = Section ,
Name-pl = Sections ,
name-pl = Sections ,
}

\zcRefTypeSetup{algorithm}{
Name-sg = Algorithm ,
name-sg = Algorithm ,
Name-pl = Algorithms ,
name-pl = Algorithms ,
}

\AddToHook{env/notation/begin}{\zcsetup{countertype={environment=notation}}}
\zcRefTypeSetup{notation}{
Name-sg = Notation ,
name-sg = Notation ,
Name-pl = Notations ,
name-pl = Notations ,
}

\AddToHook{env/question/begin}{\zcsetup{countertype={environment=question}}}
\zcRefTypeSetup{question}{
Name-sg = Question ,
name-sg = Question ,
Name-pl = Questions ,
name-pl = Questions ,
}

\AddToHook{env/problem/begin}{\zcsetup{countertype={environment=problem}}}
\zcRefTypeSetup{problem}{
Name-sg = Problem ,
name-sg = Problem ,
Name-pl = Problems ,
name-pl = Problems ,
}

\AddToHook{env/definition/begin}{\zcsetup{countertype={environment=definition}}}
\zcRefTypeSetup{definition}{
Name-sg = Definition ,
name-sg = Definition ,
Name-pl = Definitions ,
name-pl = Definitions ,
}

\zcRefTypeSetup{figure}{
Name-sg = Figure ,
name-sg = Figure ,
Name-pl = Figures ,
name-pl = Figures ,
}

\newtheoremstyle{beautypalour}
{3pt}
{3pt}
{}
{}
{\itshape}
{.}
{.5em}
{}
\theoremstyle{beautypalour}

\AddToHook{env/claim/begin}{\zcsetup{countertype={environment=claim}}}
\zcRefTypeSetup{claim}{
Name-sg = Claim ,
name-sg = Claim ,
Name-pl = Claims ,
name-pl = Claims ,
}

\usetikzlibrary{calc}
\usetikzlibrary{fit}
\usetikzlibrary{decorations}
\usetikzlibrary{decorations.pathmorphing}
\usetikzlibrary{decorations.text}
\usetikzlibrary{shapes,hobby}
\usetikzlibrary{positioning}

\tikzset{
	position/.style args={#1:#2 from #3}{
		at=($(#3)+(#1:#2)$)
	}
}

\tikzset{
  v:main/.style = {draw, circle, scale=0.8, thick,fill=black,inner sep=0.7mm},
  v:ghost/.style = {inner sep=0pt,scale=1},
  >={latex},
  e:marker/.style = {line width=8.5pt,line cap=round,opacity=0.35,color=DarkGoldenrod},
  e:main/.style = {line width=1pt},
}

\newcommand{\Bcal}{\mathcal{B}}
\newcommand{\Ccal}{\mathcal{C}}

\newcommand{\Fcal}{\mathcal{F}}
\newcommand{\Gcal}{\mathcal{G}}

\newcommand{\Lcal}{\mathcal{L}}

\newcommand{\Ncal}{\mathcal{N}}
\newcommand{\Ocal}{\mathcal{O}}
\newcommand{\Pcal}{\mathcal{P}}
\newcommand{\Qcal}{\mathcal{Q}}

\newcommand{\Scal}{\mathcal{S}}
\newcommand{\Tcal}{\mathcal{T}}

\newcommand{\Vcal}{\mathcal{V}}

\newcommand{\Xcal}{\mathcal{X}}

\newcommand{\Nbbb}{\mathbb{N}}

\RequirePackage{stmaryrd}
\usepackage{textcomp}
\DeclareUnicodeCharacter{2286}{\subseteq}
\DeclareUnicodeCharacter{2192}{\ifmmode\to\else\textrightarrow\fi}
\DeclareUnicodeCharacter{2203}{\ensuremath\exists}
\DeclareUnicodeCharacter{183}{\cdot}
\DeclareUnicodeCharacter{2200}{\forall}
\DeclareUnicodeCharacter{2264}{\leq}
\DeclareUnicodeCharacter{2265}{\geq}
\DeclareUnicodeCharacter{8614}{\mathbin{\mapsto}}
\DeclareUnicodeCharacter{8656}{\Leftarrow}
\DeclareUnicodeCharacter{8657}{\Uparrow}
\DeclareUnicodeCharacter{8658}{\Rightarrow}
\DeclareUnicodeCharacter{8659}{\Downarrow}
\DeclareUnicodeCharacter{8669}{\rightsquigarrow}
\newcommand{\eqdef}{\stackrel{{\scriptsize\rm def}}{=}}
\DeclareUnicodeCharacter{8797}{\eqdef}
\DeclareUnicodeCharacter{8870}{\vdash}
\DeclareUnicodeCharacter{8873}{\Vdash}
\DeclareUnicodeCharacter{22A7}{\models}
\DeclareUnicodeCharacter{9121}{\lceil}
\DeclareUnicodeCharacter{9123}{\lfloor}
\DeclareUnicodeCharacter{9124}{\rceil}
\DeclareUnicodeCharacter{2208}{\in}
\DeclareUnicodeCharacter{9126}{\rfloor}
\DeclareUnicodeCharacter{9655}{\triangleright}
\DeclareUnicodeCharacter{9665}{\triangleleft}
\DeclareUnicodeCharacter{9671}{\diamond}
\DeclareUnicodeCharacter{9675}{\circ}
\DeclareUnicodeCharacter{10178}{\bot}
\DeclareUnicodeCharacter{10214}{} 
\DeclareUnicodeCharacter{10215}{} 
\DeclareUnicodeCharacter{10229}{\longleftarrow}
\DeclareUnicodeCharacter{10230}{\longrightarrow}
\DeclareUnicodeCharacter{10231}{\longleftrightarrow}
\DeclareUnicodeCharacter{10232}{\Longleftarrow}
\DeclareUnicodeCharacter{10233}{\Longrightarrow}
\DeclareUnicodeCharacter{10234}{\Longleftrightarrow}
\DeclareUnicodeCharacter{10236}{\longmapsto}
\DeclareUnicodeCharacter{10238}{\Longmapsto} 
\DeclareUnicodeCharacter{10503}{\Mapsto}    
\DeclareUnicodeCharacter{10971}{\mathrel{\not\hspace{-0.2em}\cap}}
\DeclareUnicodeCharacter{65294}{\ldotp}
\DeclareUnicodeCharacter{65372}{\mid}

\newcommand{\tw}{\mathsf{tw}\xspace}
\newcommand{\obs}{\mathsf{obs}\xspace}
\newcommand{\hw}{\mathsf{hw}\xspace}
\newcommand{\torso}{\mathsf{torso}\xspace}

\newcommand{\funp}{\mathsf{f}\xspace}
\newcommand{\poly}{\mathsf{poly}\xspace}
\newcommand{\tfolio}{\mathsf{tfolio}\xspace}
\newcommand{\rep}{\mathsf{rep}\xspace}
\newcommand{\bg}{\mathsf{bg}\xspace}

\newcommand{\MSO}{\mbox{\sf MSO}\xspace}
\newcommand{\CMSO}{\mbox{\sf CMSO}\xspace}

\newcommand{\yes}{{\sf yes}\xspace}

\newcommand{\remove}[1]{}

\newcommand{\bnd}{{\sf bnd}\xspace}

\newcommand{\NP}{{\sf NP}\xspace}
\newcommand{\FPT}{{\sf FPT}\xspace}

\newcommand{\bidim}{{\sf bidim}\xspace}

\newcommand{\p}{\mathsf{p}\xspace}

\newcommand{\gall}{\mathcal{G}_{{\text{\rm  \textsf{all}}}}}
\newcommand{\opt}{\mathsf{opt}\xspace}

\newcommand{\rg}{\midnightblue{\textsl{rg}}\xspace}
\newcommand{\sg}{\midnightblue{\textsl{sg}}\xspace}
\newcommand{\rc}{\midnightblue{\textsl{rc}}\xspace}

\newcommand{\wall}{\midnightblue{\textsl{wall}}\xspace}
\newcommand{\apex}{\midnightblue{\textsl{apex}}\xspace}
\newcommand{\breadth}{\midnightblue{\textsl{breadth}}\xspace}
\newcommand{\link}{\midnightblue{\textsl{link}}\xspace}
\newcommand{\depth}{\midnightblue{\textsl{depth}}\xspace}
\newcommand{\vortex}{\midnightblue{\textsl{vortex}}\xspace}
\newcommand{\hfw}{\midnightblue{\textsl{hfw}}\xspace}
\newcommand{\genus}{\midnightblue{\textsl{genus}}\xspace}

\newcommand{\say}[1]{``#1''}

\newcommand{\adhesion}{\midnightblue{\textsl{adhesion}}\xspace}
\newcommand{\cover}{\mathsf{cover}\xspace}
\newcommand{\folio}{\mathsf{folio}\xspace}
\newcommand{\pack}{\mathsf{pack}\xspace}
\newcommand{\rh}{\mathsf{rh}\xspace}

%% file: literature_colorful.bib
@article {LarmanΜ1970configurations,
    AUTHOR = {Larman, D. G. and Mani, P.},
     TITLE = {On the existence of certain configurations within graphs and
              the {$1$}-skeletons of polytopes},
   JOURNAL = {Proc. London Math. Soc. (3)},
  FJOURNAL = {Proceedings of the London Mathematical Society. Third Series},
    VOLUME = {20},
      YEAR = {1970},
     PAGES = {144--160},
      ISSN = {0024-6115,1460-244X},
   MRCLASS = {05.50},
  MRNUMBER = {263687},
MRREVIEWER = {W.\ Mader},
       DOI = {10.1112/plms/s3-20.1.144},
       URL = {https://doi.org/10.1112/plms/s3-20.1.144},
}

@article{AldredBHK1999cubic,
  author  = {Aldred, R. E. L. and Bau, S. and Holton, D. A. and McKay, Brendan D.},
  title   = {Cycles through {$23$} vertices in {$3$}-connected cubic planar graphs},
  journal = {Graphs and Combinatorics},
  volume  = {15},
  number  = {4},
  pages   = {373--376},
  year    = {1999},
  doi     = {10.1007/s003730050046}
}

@inproceedings{AprileFJKSWY2025Integer,
  author    = {Aprile, Manuel and Fiorini, Samuel and Joret, Gwena{\"e}l and Kober, Stefan and Seweryn, Micha{\l} T. and Weltge, Stefan and Yuditsky, Yelena},
  title     = {Integer programs with nearly totally unimodular matrices: the cographic case},
  booktitle = {Proceedings of the 2025 Annual {ACM}-{SIAM} Symposium on Discrete Algorithms ({SODA})},
  pages     = {2301--2312},
  publisher = {SIAM},
  address   = {Philadelphia, PA},
  year      = {2025},
  doi       = {10.1137/1.9781611978322.76}
}

@article{ArnborgLS1991Easy,
  author  = {Arnborg, Stefan and Lagergren, Jens and Seese, Detlef},
  title   = {Easy problems for tree-decomposable graphs},
  journal = {Journal of Algorithms},
  volume  = {12},
  number  = {2},
  pages   = {308--340},
  year    = {1991},
  doi     = {10.1016/0196-6774(91)90006-K}
}

@article{ArnborgP1989Linear,
  author  = {Arnborg, Stefan and Proskurowski, Andrzej},
  title   = {Linear time algorithms for {NP}-hard problems restricted to partial {$k$}-trees},
  journal = {Discrete Applied Mathematics},
  volume  = {23},
  number  = {1},
  pages   = {11--24},
  year    = {1989},
  doi     = {10.1016/0166-218X(89)90031-0}
}

@article{BodlaenderK2007Combinatorial,
  author  = {Bodlaender, Hans L. and Koster, Arie M. C. A.},
  title   = {Combinatorial Optimization on Graphs of Bounded Treewidth},
  journal = {The Computer Journal},
  volume  = {51},
  number  = {3},
  pages   = {255--269},
  year    = {2008},
  doi     = {10.1093/comjnl/bxm037}
}

@article{BodlaenderCKN2015Deterministic,
  author  = {Bodlaender, Hans L. and Cygan, Marek and Kratsch, Stefan and Nederlof, Jesper},
  title   = {Deterministic single exponential time algorithms for connectivity problems parameterized by treewidth},
  journal = {Information and Computation},
  volume  = {243},
  pages   = {86--111},
  year    = {2015},
  doi     = {10.1016/j.ic.2014.12.008}
}

@article{BodlaenderFLPST2016Meta,
  author  = {Bodlaender, Hans L. and Fomin, Fedor V. and Lokshtanov, Daniel and Penninkx, Eelko and Saurabh, Saket and Thilikos, Dimitrios M.},
  title   = {({M}eta) kernelization},
  journal = {Journal of the ACM},
  volume  = {63},
  number  = {5},
  pages   = {44:1--44:69},
  year    = {2016},
  doi     = {10.1145/2973749}
}

@article{BollobasT1996Highly,
  author  = {Bollob{\'a}s, B{\'e}la and Thomason, Andrew},
  title   = {Highly linked graphs},
  journal = {Combinatorica},
  volume  = {16},
  number  = {3},
  pages   = {313--320},
  year    = {1996},
  doi     = {10.1007/BF01261316}
}

@article{BoriePT1992Automatic,
  author  = {Borie, Richard B. and Parker, R. Gary and Tovey, Craig A.},
  title   = {Automatic generation of linear-time algorithms from predicate calculus descriptions of problems on recursively constructed graph families},
  journal = {Algorithmica},
  volume  = {7},
  number  = {5--6},
  pages   = {555--581},
  year    = {1992},
  doi     = {10.1007/BF01758777}
}

@inproceedings{Chuzhoy2015Improved,
  author    = {Chuzhoy, Julia},
  title     = {Improved bounds for the flat wall theorem},
  booktitle = {Proceedings of the Twenty-Sixth Annual {ACM}-{SIAM} Symposium on Discrete Algorithms},
  pages     = {256--275},
  publisher = {SIAM},
  address   = {Philadelphia, PA},
  year      = {2015},
  doi       = {10.1137/1.9781611973730.20}
}

@article{chuzhoy2021towards,
  author  = {Chuzhoy, Julia and Tan, Zihan},
  title   = {Towards tight(er) bounds for the excluded grid theorem},
  journal = {Journal of Combinatorial Theory, Series B},
  volume  = {146},
  pages   = {219--265},
  year    = {2021},
  doi     = {10.1016/j.jctb.2020.09.010}
}

@article{Courcelle1992Monadic,
  author  = {Courcelle, Bruno},
  title   = {The monadic second-order logic of graphs. {III}. Tree-decompositions, minors and complexity issues},
  journal = {RAIRO Informatique Th{\'e}orique et Applications},
  volume  = {26},
  number  = {3},
  pages   = {257--286},
  year    = {1992},
  doi     = {10.1051/ita/1992260302571}
}

@incollection{Courcelle1997Expression,
  author    = {Courcelle, Bruno},
  title     = {The expression of graph properties and graph transformations in monadic second-order logic},
  booktitle = {Handbook of Graph Grammars and Computing by Graph Transformation, Vol. 1},
  pages     = {313--400},
  publisher = {World Scientific Publishing},
  address   = {River Edge, NJ},
  year      = {1997},
  doi       = {10.1142/9789812384720_0005},
  url       = {https://doi.org/10.1142/9789812384720_0005}
}

@article{Courcelle1990Monadic,
  author  = {Courcelle, Bruno},
  title   = {The monadic second-order logic of graphs. {I}. Recognizable sets of finite graphs},
  journal = {Information and Computation},
  volume  = {85},
  number  = {1},
  pages   = {12--75},
  year    = {1990},
  doi     = {10.1016/0890-5401(90)90043-H}
}

@mastersthesis{Moore2017Rooted,
  author  = {Crump, Iain},
  title   = {Forbidden Minors for {$3$}-Connected Graphs With No Non-Splitting {$5$}-Configurations},
  school  = {Simon Fraser University},
  address = {Burnaby, BC, Canada},
  year    = {2012},
  type    = {Master's thesis},
  note    = {Approved August 9, 2012},
  url     = {https://summit.sfu.ca/item/12373}
}

@inproceedings{Curticapean2016Counting,
  author    = {Curticapean, Radu},
  title     = {Counting Matchings with {$k$} Unmatched Vertices in Planar Graphs},
  booktitle = {24th Annual European Symposium on Algorithms ({ESA} 2016)},
  series    = {Leibniz International Proceedings in Informatics (LIPIcs)},
  volume    = {57},
  pages     = {33:1--33:17},
  editor    = {Sankowski, Piotr and Zaroliagis, Christos},
  publisher = {Schloss Dagstuhl -- Leibniz-Zentrum f{\"u}r Informatik},
  address   = {Dagstuhl, Germany},
  year      = {2016},
  doi       = {10.4230/LIPIcs.ESA.2016.33},
  url       = {https://drops.dagstuhl.de/entities/document/10.4230/LIPIcs.ESA.2016.33}
}

@inproceedings{Curticapean2016DefectMatching,
  author    = {Curticapean, Radu},
  title     = {Counting Matchings with {$k$} Unmatched Vertices in Planar Graphs},
  booktitle = {24th Annual European Symposium on Algorithms ({ESA} 2016)},
  series    = {Leibniz International Proceedings in Informatics (LIPIcs)},
  volume    = {57},
  pages     = {33:1--33:17},
  editor    = {Sankowski, Piotr and Zaroliagis, Christos},
  publisher = {Schloss Dagstuhl -- Leibniz-Zentrum f{\"u}r Informatik},
  address   = {Dagstuhl, Germany},
  year      = {2016},
  doi       = {10.4230/LIPIcs.ESA.2016.33},
  url       = {https://drops.dagstuhl.de/entities/document/10.4230/LIPIcs.ESA.2016.33},
  note      = {Duplicate key retained for compatibility with existing citations}
}

@book{CyganFKLMPPS2015Parameterized,
  author    = {Cygan, Marek and Fomin, Fedor V. and Kowalik, {\L}ukasz and Lokshtanov, Daniel and Marx, D{\'a}niel and Pilipczuk, Marcin and Pilipczuk, Micha{\l} and Saurabh, Saket},
  title     = {Parameterized Algorithms},
  publisher = {Springer},
  address   = {Cham},
  year      = {2015},
  edition   = {1st},
  doi       = {10.1007/978-3-319-21275-3}
}

@inproceedings{DahlhausEtAl1992,
  author    = {Dahlhaus, Elias and Johnson, David S. and Papadimitriou, Christos H. and Seymour, Paul D. and Yannakakis, Mihalis},
  title     = {The complexity of multiway cuts},
  booktitle = {Proceedings of the 24th Annual {ACM} Symposium on Theory of Computing ({STOC})},
  pages     = {241--251},
  publisher = {ACM},
  year      = {1992},
  doi       = {10.1145/129712.129736}
}

@article{Dirac1960vorhandene,
  author  = {Dirac, Gabriel Andrew},
  title   = {In abstrakten Graphen vorhandene vollst{\"a}ndige {$4$}-Graphen und ihre Unterteilungen},
  journal = {Mathematische Nachrichten},
  volume  = {22},
  pages   = {61--85},
  year    = {1960},
  doi     = {10.1002/mana.19600220107}
}

@article{DreyfusW1971Steiner,
  author  = {Dreyfus, S. E. and Wagner, R. A.},
  title   = {The Steiner problem in graphs},
  journal = {Networks},
  volume  = {1},
  pages   = {195--207},
  year    = {1971/72},
  doi     = {10.1002/net.3230010302}
}

@article{EgawaFG2000specified,
  author  = {Egawa, Yoshimi and Faudree, Ralph J. and Gy{\"o}ri, Ervin and Ishigami, Yoshiyasu and Schelp, Richard H. and Wang, Hong},
  title   = {Vertex-disjoint cycles containing specified edges},
  journal = {Graphs and Combinatorics},
  volume  = {16},
  number  = {1},
  pages   = {81--92},
  year    = {2000},
  doi     = {10.1007/s003730050005}
}

@article{EllinghamPDG2012Linkage,
  author  = {Ellingham, M. N. and Plummer, Michael D. and Yu, Gexin},
  title   = {Linkage for the diamond and the path with four vertices},
  journal = {Journal of Graph Theory},
  volume  = {70},
  number  = {3},
  pages   = {241--261},
  year    = {2012},
  doi     = {10.1002/jgt.20612}
}

@article{EricksonMV1987SendAndSplit,
  author  = {Erickson, Ranel E. and Monma, Clyde L. and Veinott, Arthur F., Jr.},
  title   = {Send-and-split method for minimum-concave-cost network flows},
  journal = {Mathematics of Operations Research},
  volume  = {12},
  number  = {4},
  pages   = {634--664},
  year    = {1987},
  doi     = {10.1287/moor.12.4.634}
}

@article{Fabila-MonroyW2013Rooted,
  author  = {Fabila-Monroy, Ruy and Wood, David R.},
  title   = {Rooted {$K_4$}-minors},
  journal = {Electronic Journal of Combinatorics},
  volume  = {20},
  number  = {2},
  pages   = {Paper 64, 19 pp.},
  year    = {2013},
  doi     = {10.37236/3476}
}

@article{FellowsL1987Nonconstructive,
  author  = {Fellows, Michael R. and Langston, Michael A.},
  title   = {Nonconstructive advances in polynomial-time complexity},
  journal = {Information Processing Letters},
  volume  = {26},
  number  = {3},
  pages   = {157--162},
  year    = {1987},
  doi     = {10.1016/0020-0190(87)90054-8}
}

@article{FellowsL1988Nonconstructive,
  author  = {Fellows, Michael R. and Langston, Michael A.},
  title   = {Nonconstructive tools for proving polynomial-time decidability},
  journal = {Journal of the ACM},
  volume  = {35},
  number  = {3},
  pages   = {727--739},
  year    = {1988},
  doi     = {10.1145/44483.44491}
}

@article{FerraraGTW2006linked,
  author  = {Ferrara, Michael and Gould, Ronald and Tansey, Gerard and Whalen, Thor},
  title   = {On {$H$}-linked graphs},
  journal = {Graphs and Combinatorics},
  volume  = {22},
  number  = {2},
  pages   = {217--224},
  year    = {2006},
  doi     = {10.1007/s00373-006-0651-6}
}

@misc{Fiorini2025Face,
  author        = {Fiorini, Samuel and Kober, Stefan and Seweryn, Micha{\l} T. and Shantanam, Abhinav and Yuditsky, Yelena},
  title         = {Face covers and rooted minors in bounded genus graphs},
  year          = {2025},
  eprint        = {2503.09230},
  archivePrefix = {arXiv},
  primaryClass  = {math.CO},
  url           = {https://arxiv.org/abs/2503.09230}
}

@article{FlandrinLM2007cycles,
  author  = {Flandrin, Evelyne and Li, Hao and Marczyk, Antoni and Wo{\'z}niak, Mariusz},
  title   = {A generalization of Dirac's theorem on cycles through {$k$} vertices in {$k$}-connected graphs},
  journal = {Discrete Mathematics},
  volume  = {307},
  number  = {7--8},
  pages   = {878--884},
  year    = {2007},
  doi     = {10.1016/j.disc.2005.11.052}
}

@inproceedings{FedorLPSZ2020Hitting,
  author    = {Fomin, Fedor V. and Lokshtanov, Daniel and Panolan, Fahad and Saurabh, Saket and Zehavi, Meirav},
  title     = {Hitting topological minors is {FPT}},
  booktitle = {{STOC} '20---Proceedings of the 52nd Annual {ACM} {SIGACT} Symposium on Theory of Computing},
  pages     = {1317--1326},
  publisher = {ACM},
  address   = {New York},
  year      = {2020},
  doi       = {10.1145/3357713.3384318}
}

@article{Frederickson1991Planar,
  author  = {Frederickson, Greg N.},
  title   = {Planar graph decomposition and all pairs shortest paths},
  journal = {Journal of the ACM},
  volume  = {38},
  number  = {1},
  pages   = {162--204},
  year    = {1991},
  doi     = {10.1145/102782.102788}
}

@article{GareyJ1977Rectilinear,
  author  = {Garey, M. R. and Johnson, D. S.},
  title   = {The rectilinear Steiner tree problem is {NP}-complete},
  journal = {SIAM Journal on Applied Mathematics},
  volume  = {32},
  number  = {4},
  pages   = {826--834},
  year    = {1977},
  doi     = {10.1137/0132071}
}

@article{GareyJS76,
  author  = {Garey, Michael R. and Johnson, David S. and Stockmeyer, Larry},
  title   = {Some simplified {NP}-complete graph problems},
  journal = {Theoretical Computer Science},
  volume  = {1},
  number  = {3},
  pages   = {237--267},
  year    = {1976},
  doi     = {10.1016/0304-3975(76)90059-X}
}

@article{GolovachKTM2017cyclability,
  author  = {Golovach, Petr A. and Kami{\'n}ski, Marcin and Maniatis, Spyridon and Thilikos, Dimitrios M.},
  title   = {The parameterized complexity of graph cyclability},
  journal = {SIAM Journal on Discrete Mathematics},
  volume  = {31},
  number  = {1},
  pages   = {511--541},
  year    = {2017},
  doi     = {10.1137/141000014}
}

@inproceedings{GolovachST23Model,
  author    = {Golovach, Petr A. and Stamoulis, Giannos and Thilikos, Dimitrios M.},
  title     = {Model-checking for first-order logic with disjoint paths predicates in proper minor-closed graph classes},
  booktitle = {Proceedings of the 2023 {ACM-SIAM} Symposium on Discrete Algorithms, {SODA} 2023},
  editor    = {Bansal, Nikhil and Nagarajan, Viswanath},
  pages     = {3684--3699},
  publisher = {SIAM},
  year      = {2023},
  doi       = {10.1137/1.9781611977554.CH141},
  url       = {https://doi.org/10.1137/1.9781611977554.CH141}
}

@misc{GorskySW2025Polynomial,
  author        = {Gorsky, Maximilian and Seweryn, Micha{\l} T. and Wiederrecht, Sebastian},
  title         = {Polynomial bounds for the graph minor structure theorem},
  year          = {2025},
  eprint        = {2504.02532},
  archivePrefix = {arXiv},
  primaryClass  = {math.CO},
  url           = {https://arxiv.org/abs/2504.02532}
}

@article{GouldW2007Subdivision,
  author  = {Gould, Ronald and Whalen, Thor},
  title   = {Subdivision extendibility},
  journal = {Graphs and Combinatorics},
  volume  = {23},
  number  = {2},
  pages   = {165--182},
  year    = {2007},
  doi     = {10.1007/s00373-006-0665-0}
}

@article{GouldKG2006linked,
  author  = {Gould, Ronald J. and Kostochka, Alexandr and Yu, Gexin},
  title   = {On minimum degree implying that a graph is {$H$}-linked},
  journal = {SIAM Journal on Discrete Mathematics},
  volume  = {20},
  number  = {4},
  pages   = {829--840},
  year    = {2006},
  doi     = {10.1137/050624662}
}

@inproceedings{GroendlandNK2024Polynomial,
  author    = {Groenland, Carla and Nederlof, Jesper and Koana, Tomohiro},
  title     = {A polynomial time algorithm for Steiner tree when terminals avoid a rooted {$K_4$}-minor},
  booktitle = {19th International Symposium on Parameterized and Exact Computation},
  series    = {Leibniz International Proceedings in Informatics (LIPIcs)},
  volume    = {321},
  pages     = {12:1--12:17},
  publisher = {Schloss Dagstuhl -- Leibniz-Zentrum f{\"u}r Informatik},
  address   = {Dagstuhl, Germany},
  year      = {2024},
  doi       = {10.4230/LIPIcs.IPEC.2024.12}
}

@inproceedings{GroheKMW11Finding,
  author    = {Grohe, Martin and Kawarabayashi, Ken-ichi and Marx, D{\'a}niel and Wollan, Paul},
  title     = {Finding topological subgraphs is fixed-parameter tractable},
  booktitle = {Proceedings of the 43rd {ACM} Symposium on Theory of Computing ({STOC} 2011)},
  editor    = {Fortnow, Lance and Vadhan, Salil P.},
  pages     = {479--488},
  publisher = {ACM},
  year      = {2011},
  doi       = {10.1145/1993636.1993700}
}

@article{GyoriP2001claw,
  author  = {Gy{\H{o}}ri, Ervin and Plummer, Michael D.},
  title   = {A nine vertex theorem for {$3$}-connected claw-free graphs},
  journal = {Studia Scientiarum Mathematicarum Hungarica},
  volume  = {38},
  pages   = {233--244},
  year    = {2001},
  doi     = {10.1556/SScMath.38.2001.1-4.16}
}

@misc{Hodor2024Quickly,
  author        = {Hodor, J{\k{e}}drzej and La, Hoang and Micek, Piotr and Rambaud, Cl{\'e}ment},
  title         = {Quickly excluding an apex-forest},
  year          = {2024},
  eprint        = {2404.17306},
  archivePrefix = {arXiv},
  primaryClass  = {math.CO},
  url           = {https://arxiv.org/abs/2404.17306}
}

@misc{HodorLMR24quick,
  author        = {Hodor, J{\k{e}}drzej and La, Hoang and Micek, Piotr and Rambaud, Cl{\'e}ment},
  title         = {Quickly excluding an apex-forest},
  year          = {2024},
  eprint        = {2404.17306},
  archivePrefix = {arXiv},
  primaryClass  = {math.CO},
  url           = {https://arxiv.org/abs/2404.17306}
}

@inproceedings{JansenS2024SteinerTree,
  author    = {Jansen, Bart M. P. and Swennenhuis, C{\'e}line M. F.},
  title     = {Steiner tree parameterized by multiway cut and even less},
  booktitle = {32nd Annual European Symposium on Algorithms},
  series    = {Leibniz International Proceedings in Informatics (LIPIcs)},
  volume    = {308},
  pages     = {76:1--76:16},
  publisher = {Schloss Dagstuhl -- Leibniz-Zentrum f{\"u}r Informatik},
  address   = {Dagstuhl, Germany},
  year      = {2024},
  doi       = {10.4230/LIPIcs.ESA.2024.76}
}

@inproceedings{JohnsonMPPSL2024multiway,
  author    = {Johnson, Matthew and Martin, Barnaby and Pandey, Sukanya and Paulusma, Dani{\"e}l and Smith, Siani and van Leeuwen, Erik Jan},
  title     = {Edge multiway cut and node multiway cut are hard for planar subcubic graphs},
  booktitle = {19th Scandinavian Symposium on Algorithm Theory},
  series    = {Leibniz International Proceedings in Informatics (LIPIcs)},
  volume    = {294},
  pages     = {29:1--29:17},
  publisher = {Schloss Dagstuhl -- Leibniz-Zentrum f{\"u}r Informatik},
  address   = {Dagstuhl, Germany},
  year      = {2024},
  doi       = {10.4230/LIPIcs.SWAT.2024.29}
}

@article{Jung1970Verallgemeinerung,
  author  = {Jung, H. A.},
  title   = {Eine Verallgemeinerung des {$n$}-fachen Zusammenhangs f{\"u}r Graphen},
  journal = {Mathematische Annalen},
  volume  = {187},
  pages   = {95--103},
  year    = {1970},
  doi     = {10.1007/BF01350174}
}

@article{Kawarabayashi2004Rooted,
  author  = {Kawarabayashi, Ken-ichi},
  title   = {Rooted minor problems in highly connected graphs},
  journal = {Discrete Mathematics},
  volume  = {287},
  number  = {1--3},
  pages   = {121--123},
  year    = {2004},
  doi     = {10.1016/j.disc.2004.07.007}
}

@article{KawarabayashiKG2006linked,
  author  = {Kawarabayashi, Ken-ichi and Kostochka, Alexandr and Yu, Gexin},
  title   = {On sufficient degree conditions for a graph to be {$k$}-linked},
  journal = {Combinatorics, Probability and Computing},
  volume  = {15},
  number  = {5},
  pages   = {685--694},
  year    = {2006},
  doi     = {10.1017/S0963548305007479}
}

@article{KwarabayashiTW2018NewProof,
  author  = {Kawarabayashi, Ken-ichi and Thomas, Robin and Wollan, Paul},
  title   = {A new proof of the flat wall theorem},
  journal = {Journal of Combinatorial Theory, Series B},
  volume  = {129},
  pages   = {204--238},
  year    = {2018},
  doi     = {10.1016/j.jctb.2017.09.006}
}

@misc{KawarabayashiTW2020Quickly,
  author        = {Kawarabayashi, Ken-ichi and Thomas, Robin and Wollan, Paul},
  title         = {Quickly excluding a non-planar graph},
  year          = {2020},
  eprint        = {2010.12397},
  archivePrefix = {arXiv},
  primaryClass  = {math.CO},
  url           = {https://arxiv.org/abs/2010.12397}
}

@inproceedings{KawarabayashiW2010Shorter,
  author    = {Kawarabayashi, Ken-ichi and Wollan, Paul},
  title     = {A shorter proof of the graph minor algorithm---the unique linkage theorem---[extended abstract]},
  booktitle = {{STOC}'10---Proceedings of the 2010 {ACM} International Symposium on Theory of Computing},
  pages     = {687--694},
  publisher = {ACM},
  address   = {New York},
  year      = {2010},
  doi       = {10.1145/1806689.1806784}
}

@article{KisfaludiBakNvL2020NearlyETH,
  author  = {Kisfaludi-Bak, S{\'a}ndor and Nederlof, Jesper and van Leeuwen, Erik Jan},
  title   = {Nearly {ETH}-tight algorithms for planar Steiner tree with terminals on few faces},
  journal = {ACM Transactions on Algorithms},
  volume  = {16},
  number  = {3},
  pages   = {28:1--28:30},
  year    = {2020},
  doi     = {10.1145/3371389}
}

@inproceedings{KorhonenPS2024Minor,
  author    = {Korhonen, Tuukka and Pilipczuk, Micha{\l} and Stamoulis, Giannos},
  title     = {Minor containment and disjoint paths in almost-linear time},
  booktitle = {2024 {IEEE} 65th Annual Symposium on Foundations of Computer Science ({FOCS} 2024)},
  pages     = {53--61},
  publisher = {IEEE Computer Society},
  address   = {Los Alamitos, CA},
  year      = {2024},
  doi       = {10.1109/FOCS61266.2024.00014},
  eprint    = {2404.03958},
  archivePrefix = {arXiv},
  url       = {https://arxiv.org/abs/2404.03958}
}

@article{Kostochka2005extremal,
  author  = {Kostochka, Alexandr and Yu, Gexin},
  title   = {An extremal problem for {$H$}-linked graphs},
  journal = {Journal of Graph Theory},
  volume  = {50},
  number  = {4},
  pages   = {321--339},
  year    = {2005},
  doi     = {10.1002/jgt.20115}
}

@inproceedings{KrauthgamerLRr2019FlowCut,
  author    = {Krauthgamer, Robert and Lee, James R. and Rika, Havana},
  title     = {Flow-cut gaps and face covers in planar graphs},
  booktitle = {Proceedings of the Thirtieth Annual {ACM}-{SIAM} Symposium on Discrete Algorithms},
  pages     = {525--534},
  publisher = {SIAM},
  address   = {Philadelphia, PA},
  year      = {2019},
  doi       = {10.1137/1.9781611975482.33}
}

@incollection{Marx12,
  author    = {Marx, D{\'a}niel},
  title     = {A tight lower bound for planar multiway cut with fixed number of terminals},
  booktitle = {Automata, Languages, and Programming. Part {I}},
  series    = {Lecture Notes in Computer Science},
  volume    = {7391},
  pages     = {677--688},
  publisher = {Springer},
  address   = {Heidelberg},
  year      = {2012},
  doi       = {10.1007/978-3-642-31594-7_57},
  url       = {https://doi.org/10.1007/978-3-642-31594-7_57}
}

@article{MarxSW2017Rooted,
  author  = {Marx, D{\'a}niel and Seymour, Paul D. and Wollan, Paul},
  title   = {Rooted grid minors},
  journal = {Journal of Combinatorial Theory, Series B},
  volume  = {122},
  pages   = {428--437},
  year    = {2017},
  doi     = {10.1016/j.jctb.2016.07.003}
}

@article{OkamuraSeymour1981,
  author  = {Okamura, Hiroshi and Seymour, Paul D.},
  title   = {Multicommodity flows in planar graphs},
  journal = {Journal of Combinatorial Theory, Series B},
  volume  = {31},
  number  = {1},
  pages   = {75--81},
  year    = {1981},
  doi     = {10.1016/0095-8956(81)90005-7}
}

@misc{PandeyL2025Planar,
  author        = {Pandey, Sukanya and van Leeuwen, Erik Jan},
  title         = {Planar multiway cut with terminals on few faces},
  year          = {2025},
  eprint        = {2506.23399},
  archivePrefix = {arXiv},
  primaryClass  = {cs.DS},
  url           = {https://arxiv.org/abs/2506.23399}
}

@misc{PaulPTW2024obstructionsArXiV,
  author        = {Paul, Christophe and Protopapas, Evangelos and Thilikos, Dimitrios M. and Wiederrecht, Sebastian},
  title         = {Obstructions to Erd{\H{o}}s-{P}{\'o}sa dualities for minors},
  year          = {2024},
  eprint        = {2407.09671},
  archivePrefix = {arXiv},
  primaryClass  = {math.CO},
  url           = {https://arxiv.org/abs/2407.09671}
}

@inproceedings{PaulPTW2024Obstructions,
  author    = {Paul, Christophe and Protopapas, Evangelos and Thilikos, Dimitrios M. and Wiederrecht, Sebastian},
  title     = {Obstructions to Erd{\H{o}}s-{P}{\'o}sa dualities for minors},
  booktitle = {2024 {IEEE} 65th Annual Symposium on Foundations of Computer Science ({FOCS} 2024)},
  pages     = {31--52},
  publisher = {IEEE Computer Society},
  address   = {Los Alamitos, CA},
  year      = {2024},
  doi       = {10.1109/FOCS61266.2024.00013},
  eprint    = {2407.09671},
  archivePrefix = {arXiv},
  url       = {https://arxiv.org/abs/2407.09671}
}

@misc{PaulPTS2025LocalIndex,
  author        = {Paul, Christophe and Protopapas, Evangelos and Thilikos, Dimitrios M. and Wiederrecht, Sebastian},
  title         = {The local structure theorem for graph minors with finite index},
  year          = {2025},
  eprint        = {2507.02769},
  archivePrefix = {arXiv},
  primaryClass  = {math.CO},
  url           = {https://arxiv.org/abs/2507.02769}
}

@inproceedings{Reed1992Finding,
  author    = {Reed, Bruce A.},
  title     = {Finding approximate separators and computing tree width quickly},
  booktitle = {Proceedings of the Twenty-Fourth Annual {ACM} Symposium on Theory of Computing ({STOC} '92)},
  pages     = {221--228},
  publisher = {ACM Press},
  address   = {Victoria, British Columbia, Canada},
  year      = {1992},
  doi       = {10.1145/129712.129734}
}

@incollection{RobertsonS1984Grapha,
  author    = {Robertson, Neil and Seymour, Paul D.},
  title     = {Graph Width and Well-Quasi-Ordering: A Survey},
  booktitle = {Progress in Graph Theory},
  volume    = {2},
  pages     = {399--406},
  publisher = {Academic Press},
  address   = {Toronto, Orlando},
  year      = {1984}
}

@article{RobertsonS1986Graphb,
  author  = {Robertson, Neil and Seymour, Paul D.},
  title   = {Graph Minors. {II}. Algorithmic aspects of tree-width},
  journal = {Journal of Algorithms},
  volume  = {7},
  number  = {3},
  pages   = {309--322},
  year    = {1986},
  doi     = {10.1016/0196-6774(86)90023-4}
}

@article{RobertsonS1986Grapha,
  author  = {Robertson, Neil and Seymour, Paul D.},
  title   = {Graph Minors. {V}. Excluding a planar graph},
  journal = {Journal of Combinatorial Theory, Series B},
  volume  = {41},
  number  = {1},
  pages   = {92--114},
  year    = {1986},
  doi     = {10.1016/0095-8956(86)90030-4}
}

@article{RobertsonS1990Grapha,
  author  = {Robertson, Neil and Seymour, Paul D.},
  title   = {Graph Minors. {IV}. Tree-width and well-quasi-ordering},
  journal = {Journal of Combinatorial Theory, Series B},
  volume  = {48},
  number  = {2},
  pages   = {227--254},
  year    = {1990},
  doi     = {10.1016/0095-8956(90)90120-O}
}

@article{RobertsonS1991Graph,
  author  = {Robertson, Neil and Seymour, Paul D.},
  title   = {Graph Minors. {X}. Obstructions to tree-decomposition},
  journal = {Journal of Combinatorial Theory, Series B},
  volume  = {52},
  number  = {2},
  pages   = {153--190},
  year    = {1991},
  doi     = {10.1016/0095-8956(91)90061-N}
}

@article{RobertsonS1995Graph,
  author  = {Robertson, Neil and Seymour, Paul D.},
  title   = {Graph Minors. {XIII}. The disjoint paths problem},
  journal = {Journal of Combinatorial Theory, Series B},
  volume  = {63},
  number  = {1},
  pages   = {65--110},
  year    = {1995},
  doi     = {10.1006/jctb.1995.1006}
}

@article{RobertsonS2003GraphMinorsXVI,
  author  = {Robertson, Neil and Seymour, Paul D.},
  title   = {Graph Minors. {XVI}. Excluding a non-planar graph},
  journal = {Journal of Combinatorial Theory, Series B},
  volume  = {89},
  number  = {1},
  pages   = {43--76},
  year    = {2003},
  doi     = {10.1016/S0095-8956(03)00042-X}
}

@article{GraphMinorsXXIII,
  author  = {Robertson, Neil and Seymour, Paul D.},
  title   = {Graph minors {XXIII}. Nash-Williams' immersion conjecture},
  journal = {Journal of Combinatorial Theory, Series B},
  volume  = {100},
  number  = {2},
  pages   = {181--205},
  year    = {2010},
  doi     = {10.1016/j.jctb.2009.07.003}
}

@inproceedings{SauST25Parameterizing,
  author    = {Sau, Ignasi and Stamoulis, Giannos and Thilikos, Dimitrios M.},
  title     = {Parameterizing the quantification of {CMSO}: model checking on minor-closed graph classes},
  booktitle = {Proceedings of the 2025 Annual {ACM-SIAM} Symposium on Discrete Algorithms, {SODA} 2025},
  editor    = {Azar, Yossi and Panigrahi, Debmalya},
  pages     = {3728--3742},
  publisher = {SIAM},
  year      = {2025},
  doi       = {10.1137/1.9781611978322.124},
  eprint    = {2406.18465},
  archivePrefix = {arXiv},
  url       = {https://arxiv.org/abs/2406.18465}
}

@book{Schrijver2003,
  author    = {Schrijver, Alexander},
  title     = {Combinatorial Optimization: Polyhedra and Efficiency},
  volume    = {B},
  publisher = {Springer},
  year      = {2003}
}

@article{ThilikosW2024Killing,
  author  = {Thilikos, Dimitrios M. and Wiederrecht, Sebastian},
  title   = {Killing a vortex},
  journal = {Journal of the ACM},
  volume  = {71},
  number  = {4},
  pages   = {27:1--27:56},
  year    = {2024},
  doi     = {10.1145/3664648}
}

@misc{ThilikosW2025graphminorsstructure,
  author        = {Thilikos, Dimitrios M. and Wiederrecht, Sebastian},
  title         = {The graph minors structure theorem through bidimensionality},
  year          = {2025},
  eprint        = {2306.01724},
  archivePrefix = {arXiv},
  primaryClass  = {math.CO},
  url           = {https://arxiv.org/abs/2306.01724}
}

@misc{ThilikosW2024Excluding,
  author        = {Thilikos, Dimitrios M. and Wiederrecht, Sebastian},
  title         = {Excluding surfaces as minors in graphs},
  year          = {2026},
  eprint        = {2601.19230},
  archivePrefix = {arXiv},
  primaryClass  = {math.CO},
  url           = {https://arxiv.org/abs/2601.19230}
}

@article{ThomasW2005linkages,
  author  = {Thomas, Robin and Wollan, Paul},
  title   = {An improved linear edge bound for graph linkages},
  journal = {European Journal of Combinatorics},
  volume  = {26},
  number  = {3--4},
  pages   = {309--324},
  year    = {2005},
  doi     = {10.1016/j.ejc.2004.02.013}
}

@article{HofPW2009Partitioning,
  author  = {van {'t} Hof, Pim and Paulusma, Dani{\"e}l and Woeginger, Gerhard J.},
  title   = {Partitioning graphs into connected parts},
  journal = {Theoretical Computer Science},
  volume  = {410},
  number  = {47--49},
  pages   = {4834--4843},
  year    = {2009},
  doi     = {10.1016/j.tcs.2009.06.028}
}

@article{WatkinsM1967Cycles,
  author  = {Watkins, M. E. and Mesner, D. M.},
  title   = {Cycles and connectivity in graphs},
  journal = {Canadian Journal of Mathematics},
  volume  = {19},
  pages   = {1319--1328},
  year    = {1967},
  doi     = {10.4153/CJM-1967-121-2}
}

@article{Watkins1968arcs,
  author  = {Watkins, Mark E.},
  title   = {On the existence of certain disjoint arcs in graphs},
  journal = {Duke Mathematical Journal},
  volume  = {35},
  pages   = {231--246},
  year    = {1968},
  url     = {http://projecteuclid.org/euclid.dmj/1077377609}
}

@article{Wollan2008Extremal,
  author  = {Wollan, Paul},
  title   = {Extremal functions for rooted minors},
  journal = {Journal of Graph Theory},
  volume  = {58},
  number  = {2},
  pages   = {159--178},
  year    = {2008},
  doi     = {10.1002/jgt.20301}
}

@article{Hananya1994disjoint,
  author  = {Yinnone, Hananya},
  title   = {Maximum number of disjoint paths connecting specified terminals in a graph},
  journal = {Discrete Applied Mathematics},
  volume  = {55},
  number  = {2},
  pages   = {183--195},
  year    = {1994},
  doi     = {10.1016/0166-218X(94)90007-8}
}
